\documentclass[10pt]{article}
\usepackage{latexsym}
\usepackage{times}
\usepackage{mathrsfs}
\usepackage{stmaryrd}
\usepackage{amsopn}
\usepackage{amsmath}
\usepackage{amssymb}
\usepackage{amsfonts}
\usepackage{amsbsy}
\usepackage{amscd}
\newcommand{\cb}{C_{\rm b}}
\newcommand{\cbu}{C_{\rm b}^{\rm u}}
\newcommand{\cu}{C^{\rm u}}
\newcommand{\co}{{C_0}}
\newcommand{\cc}{C_{\rm c}}
\newcommand{\C}{\mathbb{C}}

\newcommand{\N}{\mathbb{N}}
\newcommand{\R}{\mathbb{R}}

\newcommand{\Z}{\mathbb{Z}}
\newcommand{\Ac}{\mathcal{A}}
\newcommand{\Bc}{\mathcal{B}}
\newcommand{\Cc}{\mathcal{C}}
\newcommand{\Dc}{\mathcal{D}}
\newcommand{\Kc}{\mathcal{K}}
\newcommand{\Fc}{\mathcal{F}}
\newcommand{\Mc}{\mathcal{M}}
\newcommand{\Nc}{\mathcal{N}}

\newcommand{\Br}{\rond{B}}

\newcommand{\Er}{\rond{E}}
\newcommand{\Fr}{\rond{F}}
\newcommand{\Gr}{\rond{G}}
\newcommand{\Hr}{\rond{H}}

\newcommand{\Kr}{\rond{K}}
\newcommand{\Lr}{\rond{L}}
\newcommand{\Mr}{\rond{M}}
\newcommand{\Nr}{\rond{N}}

\newcommand{\Sr}{\rond{S}}

\font\tenrsf=rsfs10 at 11pt
\font\sevenrsf=rsfs7 at 8pt
\font\fiversf=rsfs5 at 6pt
\newfam\rsffam
\textfont\rsffam=\tenrsf
\scriptfont\rsffam=\sevenrsf
\scriptscriptfont\rsffam=\fiversf
\def\rond#1{{\tenrsf\fam\rsffam#1}}

\newcommand{\proof}{\noindent{\bf Proof: }}

\def\cchi{\raisebox{.45 ex}{$\chi$}} 
\def\ra{\rangle}
\def\se{\sigma_{\rm ess}}

\def\w{{\rm w}}
\def\d{{\rm d}}
\def\crm{{\rm c}}
\def\veps{\varepsilon}
\def\la{\langle}
\def\rarrow{\rightarrow}
\def\vphi{\varphi}
\def\what{\widehat}
\def\wtilde{\widetilde}

\def\slim{\mbox{\rm s-}\!\lim}

\def\supp{\mbox{\rm supp }}
\def\Re{\mbox{\rm Re }}

\def\qed{\hfill \vrule width 8pt height 9pt depth-1pt \medskip}

\def\build#1_#2^#3{\mathrel{\mathop{\kern 0pt#1}\limits_{#2}^{#3}}}
\def\bol{\Bc_0^{\,l}}
\def\bor{\Bc_0^{\,r}}
\def\bo{\Bc_0}
\def\bool{\Bc_{00}^{\,l}}

\def\bql{\Bc_q^{\,l}}
\def\bqr{\Bc_q^{\,r}}
\def\bq{\Bc_q}

\newtheorem{theorem}{\rm \bf Theorem}[section]
\newtheorem{lemma}[theorem]{Lemma}
\newtheorem{proposition}[theorem]{Proposition}
\newtheorem{corollary}[theorem]{Corollary}

\newtheorem{definition}[theorem]{Definition}
\newtheorem{example}[theorem]{Example}
\newtheorem{remark}[theorem]{Remark}
\newtheorem{remarks}[theorem]{Remarks}
\font\teneuf=eufm10 at 12pt
\font\seveneuf=eufm7 at 8pt
\font\fiveeuf=eufm5 at 6pt
\newfam\euffam
\textfont\euffam=\teneuf
\scriptfont\euffam=\seveneuf
\scriptscriptfont\euffam=\fiveeuf

\newfont{\secgoth}{eufm10 at 16pt}
\begin{document}
\title{Decay Preserving Operators and Stability of the Essential
  Spectrum} 
\author{ V. Georgescu and S. Gol\'enia\\
CNRS and  D\'epartement de Math\'ematiques\\
Universit\'e de Cergy-Pontoise\\
2 avenue Adolphe Chauvin \\
95302 Cergy-Pontoise Cedex France}
\date{}
\maketitle
\begin{center}
Revised Version
\end{center} 
\begin{abstract}
  We establish criteria for the stability of the essential
  spectrum for unbounded operators acting in Banach modules. The
  applications cover operators acting on sections of vector fiber
  bundles over non-smooth manifolds or locally compact abelian
  groups, in particular differential operators of any order with
  complex measurable coefficients on $\R^n$, singular Dirac
  operators, and Laplace-Beltrami operators on Riemannian manifolds
  with measurable metrics.
\end{abstract}
\tableofcontents
\section{Introduction}
The main purpose of this paper is to establish criteria which ensure
that the difference  of the resolvents of two  operators is compact. 
In  order  to  simplify  later  statements,  we  use  the  following
definition (our notations are  quite standard; we recall however the
most important ones at the end of this section).

\begin{definition}\label{d:int} 
  Let $A$ and $B$ be two closed operators acting in a Banach space
  $\Hr$.  We say that {\em $B$ is a compact perturbation of $A$} if
  there is $z\in\rho(A)\cap\rho(B)$ such that
  $(A-z)^{-1}-(B-z)^{-1}$ is a compact operator.
\end{definition} 

Under the  conditions of this definition the difference
$(A-z)^{-1}-(B-z)^{-1}$ is a compact operator for all
$z\in\rho(A)\cap\rho(B)$. In particular,
{\em if $B$ is a compact perturbation of $A$, then $A$ and $B$
have the same essential spectrum}, and this for any reasonable
definition of the essential spectrum, see \cite{GW}.
To be precise, in this paper we define the essential spectrum
of $A$ as the set of points $\lambda\in\C$ such that $A-\lambda$
is not Fredholm.

We shall describe now a standard and simple, although quite
powerful, method of proving that $B$ is a compact perturbation of
$A$.  Note that we are interested in situations where $A$ and $B$
are differential (or pseudo-differential) operators with complex
measurable coefficients which differ little on a neighborhood of
infinity. An important point in such situations is that one has not
much information about the domains of the operators. However, one
often knows explicitly a generalized version of the ``quadratic form
domain'' of the operator.  Since we want to consider operators of
any order (in particular Dirac operators) we shall work in the
following framework, which goes beyond the theory of accretive
forms.

Let $\Gr,\Hr,\Kr$ be reflexive Banach spaces such that
$\Gr\subset\Hr\subset\Kr$ continuously and densely. We are
interested in operators in $\Hr$ constructed according to the
following procedure: let $A_0,B_0$ be continuous bijective maps
$\Gr\rarrow\Kr$ and let $A,B$ be their restrictions to $A_0^{-1}\Hr$
and $B_0^{-1}\Hr$. These are closed densely defined operators in
$\Hr$ and $z=0\in\rho(A)\cap\rho(B)$.  Then in
$\Bc(\Kr,\Gr)$ we have
\begin{equation}\label{e:!} 
A_0^{-1}-B_0^{-1}=A_0^{-1}(B_0-A_0)B_0^{-1}.
\end{equation} 
In particular, we get in $\Bc(\Hr)$
\begin{equation}\label{e:!!}
A^{-1}-B^{-1}=A_0^{-1}(B_0-A_0)B^{-1}.
\end{equation} 
We get the simplest compactness criterion: if
$A_0-B_0:\Gr\rarrow\Kr$ is compact, then $B$ is a compact
perturbation of $A$. But in this case we have more: the operator
$A_0^{-1}-B_0^{-1}:\Kr\rarrow\Gr$ is also compact, and this can not
happen if $A_0,B_0$ are differential operators with distinct
principal part (cf.\ below). This also excludes singular lower order
perturbations, e.g. Coulomb potentials in the Dirac case.

The advantage of the preceding criterion is that no knowledge of the
domains $\Dc(A),\Dc(B)$ is needed. To avoid the mentioned
disadvantages, one may assume that one of the operators is more
regular than the second one, so that the functions in its domain
are, at least locally, slightly better than those from $\Gr$.  Note
that $\Dc(B)$ when equipped with the graph topology is such that
$\Dc(B)\subset\Gr$ continuously and densely and we get a second
compactness criterion by asking that $A_0-B_0:\Dc(B)\rarrow\Kr$ be
compact. This time again we get more than needed, because not only
$B$ is a compact perturbation of $A$, but also
$A_0^{-1}-B_0^{-1}:\Hr\rarrow\Gr$ is compact. However, perturbations
of the principal part of a differential operator are allowed and
also much more singular perturbations of the lower order terms, cf.\ 
\cite{N1} for the Dirac case.

In this paper we are interested in situations where we have really
no information concerning the domains of $A$ and $B$ (besides the
fact that they are subspaces of $\Gr$). The case when $A,B$ are
second order elliptic operators with measurable complex coefficients
acting in $\Hr=L^2(\R^n)$ has been studied by Ouhabaz and Stollmann
in \cite{OS} and, as far as we know, this is the only paper where
the ``unperturbed'' operator is not smooth.  Their approach consists
in proving that the difference $A^{-k}-B^{-k}$ is compact for some
$k\geq 2$ (which implies the compactness of $A^{-1}-B^{-1}$). In
order to prove this, they take advantage of the fact that $\Dc(A^k)$
is a subset of the Sobolev space $W^{1,p}$ for some $p>2$, which
means that we have a certain gain of local regularity. Of course,
$L^p$ techniques from the theory of partial differential equations
are required for their methods to work.

We shall explain now in the most elementary situation the main ideas
of our approach to these questions. Let $\Hr= L^2(\R)$ and
$P=-i\frac{d}{dx}$.  We consider operators of the form $A_0=PaP+V$
and $B_0=PbP+W$ where $a,b$ are bounded operators on $\Hr$ such that
$\Re a$ and $\Re b$ are bounded below by strictly positive numbers.
$V$ and $W$ are assumed to be continuous operators
$\Hr^{1}\rarrow\Hr^{-1}$, where $\Hr^{s}$ are Sobolev spaces
associated to $\Hr$.  Then $A_0,B_0\in\Bc(\Hr^{1},\Hr^{-1})$ and we
put some conditions on $V,W$ which ensure that $A_0,B_0$ are
invertible (e.g.\ we could include the constant $z$ in them). Thus
we are in the preceding abstract framework with $\Gr=\Hr^{1}$ and
$\Kr=\Hr^{-1}\equiv\Gr^*$.  Then from (\ref{e:!!}) we get
\begin{equation}\label{e:!!!}
A^{-1}-B^{-1}=A_0^{-1}P(b-a)PB^{-1}+A_0^{-1}(W-V)B^{-1}.
\end{equation}
Let $R$ be the first term on the right hand side and let us see how
we could prove that it is a compact operator on $\Hr$. Note that the
second term should be easier to treat since we expect $V$ and $W$ to
be operators of order less than $2$.

We have $R\Hr\subset\Hr^1$, so we can write $R=\psi(P)R_1$ for some
$\psi\in B_0(\R)$ (bounded Borel function which tends to zero at
infinity) and $R_1\in\Bc(\Hr)$.  This is just half of the conditions
needed for compactness, in fact $R$ will be compact if and only if
one can also find $\vphi\in B_0(\R)$ and $R_2\in\Bc(\Hr)$ such that
$R=\vphi(Q)R_2$, where $\vphi(Q)$ is the operator of multiplication
by $\vphi$. Of course, the only factor which can help to get such a
decay is $b-a$. So let us suppose that we can write $b-a=\xi(Q)U$
for some $\xi\in B_0(\R)$ and a bounded operator $U$ on $\Hr$. We
denote $S=A_0^{-1}P$ and note that this is a bounded operator on
$\Hr$, because $P:\Hr\rarrow\Hr^{-1}$ and
$A_0^{-1}:\Hr^{-1}\rarrow\Hr^{1}$ are bounded. Then
$R=S\xi(Q)UPB^{-1}$ and $UPB^{-1}\in\Bc(\Hr)$, hence $R$ will be
compact if the operator $S\in\Bc(\Hr)$ has the following property:
for each $\xi\in B_0(\R)$ there are $\vphi\in B_0(\R)$ and
$T\in\Bc(\Hr)$ such that $S\xi(Q)=\vphi(Q)T$.

An operator $S$ with the property specified above will be called
\emph{decay preserving}.  Thus we see that the compactness of $R$
follows from the fact that $S$ preserves decay and our main point is
that it is easy to check this property under very general
assumptions on $A$, cf.\ Corollary \ref{ex27} and Proposition
\ref{prop28} for abstract criteria, Lemmas \ref{l:ab}, \ref{l:40}
and \ref{l:analytic} and Theorems \ref{t:fquasi} and
\ref{th:fquasi} for more concrete examples. Note that the
perturbative technique described in Proposition \ref{pr:pco}
shows that in many cases it suffices to prove the decay preserving
property only for operators with smooth coefficients (cf.\ Lemma
\ref{l:analytic}). 

An abstract formulation of the ideas described above (see
Proposition \ref{p:comp}) allows one to treat situations of a very
general nature, like pseudo-differential operators on finite
dimensional vector spaces over a local\footnote{ See \cite{Sa,Ta}
  for the corresponding pseudo-differential calculus.  } (for
example $p$-adic) field, in particular differential operators of
arbitrary order with irregular coefficients on $\R^n$, the Laplace
operator on manifolds with locally $L^\infty$ Riemannian metrics,
and operators acting on sections of vector bundles over locally
compact spaces. Sections \ref{lcg}, \ref{div}, \ref{wvp} and
\ref{lrm} are devoted to such applications.  We stress once again
that, in the applications to differential operators, we are
interested only in situations where the coefficients are not smooth
and the lower order terms are singular.

\medskip

\noindent{\bf Plan of the paper:}
In Section \ref{pr} we introduce an algebraic formalism which allows
us to treat in a unified and simple way operators which have an
algebraically complicated structure, e.g.\ operators acting on
sections of vector fiber bundles over a locally compact space. The
class of decay improving (or vanishing at infinity) operators is
defined through an a priori given algebra of operators on a Banach
space $\Hr$, that we call \emph{multiplier algebra of} $\Hr$, and
this allows us to define the notion of decay preserving operator in
a natural and general context, that of Banach modules.  Several
examples of multiplier algebras are given Subsections \ref{bm},
\ref{lcg} and \ref{wvpr}. We stress that Section \ref{pr} is only an
accumulation of definitions and straightforward consequences.

We mention that this algebraic framework allows one to study
differential operators in $L^p$ or more general Banach spaces. Since
these extensions are rather obvious and the examples are not
particularly interesting, we shall not consider explicitly such
situations.

Section \ref{mr} contains several abstract compactness criteria which
formalize in the context of Banach modules the ideas involved in the
example discussed above.

In Subsection \ref{anyord} we give our first concrete applications
of the abstract theory: we consider ``hypoelliptic'' operators on
abelian groups and treat as an example the Dirac operators on
$\R^n$. In Section \ref{div} we discuss operators in divergence
form on $\R^n$, hence of order $2m$ with $m\geq1$ integer, with
coefficients of a rather general form (they do not have to be
functions, for example).

In Section \ref{wvp} we present several results concerning the case
when the coefficients of the operator $A-B$ vanish at infinity only
in some weak sense. This question has been studied before, for
example in \cite{He,LV,OS,We}.  We present the notion of
\emph{weakly vanishing at infinity} functions in terms of filters
finer than the Fr\'echet filter, a natural idea in our context being
to extend the standard notion of neighborhood of infinity.

If $X$ is a locally compact space, it is usual to define the filter
of neighborhoods of infinity as the family of subsets of $X$ with
relatively compact complement; we shall call this the
\emph{Fr\'echet filter}. If $\Fc$ is a filter on $X$ finer than the
Fr\'echet filter then a function $\vphi:X\rarrow\C$ such that
$\lim_\Fc\vphi=0$ can naturally be thought as convergent to zero at
infinity in a generalized sense (recall that $\lim_\Fc\vphi=0$ means
that for each $\veps>0$ the set of points $x$ such that
$|\vphi(x)|<\veps$ belongs to $\Fc$ ).  In Subsection \ref{wvpr} we
consider three such filters and describe corresponding classes of
decay preserving operators in Theorem \ref{th:meas}, Proposition
\ref{p:Lql} and Theorem \ref{t:fquasi}.

Theorem \ref{t:fquasi} is a consequence of a factorization theorem
that we prove in Section \ref{mft} and which involves interesting
tools from the modern theory of Banach spaces. In fact, Theorem
\ref{t:maurey}, the main result of Section \ref{mft}, is a version
of the ``strong factorization theorem'' of B.\ Maurey (see Theorem
\ref{th:bmsft}) which does not seem to be covered by the results
existing in the literature. We also use Maurey's theorem directly to
prove some of our main results, for example Theorems
\ref{th:wmetric} and \ref{th:wmetric*} which depend on Theorem
\ref{th:meas}.

Theorem \ref{t:wos} is one of the main applications of our
formalism: we prove a compactness result for operators of order $2m$
in divergence form assuming that the difference between their
coefficients vanishes at infinity in a weak sense. Such results were
known before only in the case $m=1$, see especially Theorem 2.1 in
\cite{OS}. We assume that the coefficients of the higher order terms
are bounded, thus their Theorem 3.1 is not covered unless we add an
implicit assumption, as is done in \cite{OS} (or in our Theorems
\ref{th:wmetric} and \ref{th:wmetric*}).  In fact, our main abstract
compactness result Theorem \ref{t:main} is stated such as to apply
to situations when the coefficients of the principal part of the
operators are locally unbounded, as in \cite{Ba1,Ba2}, but we have
not developed this idea here.

Perturbations of the Laplace operator on a Riemannian manifold with
locally $L^\infty$ metric are considered in Section \ref{lrm}. We
introduce and study an abstract model of this situation which fits
very naturally in our algebraic framework and covers the case of
Lipschitz manifolds with measurable metrics. We consider in more
detail the case when the manifold is $C^1$ (but the metric is only
locally $L^\infty$) and establish stability of the essential
spectrum under certain perturbations of the metric, see Theorems
\ref{t:metric}, \ref{th:wmetric} and \ref{th:wmetric*}. We also
consider, in an abstract setting and without going into technical
details, the Laplace operator acting on differential forms.

In an Appendix we collect some general facts concerning operators
acting in scales of spaces which are often used without comment in
the rest of the paper.

\medskip
\noindent{\bf Notations:}
If $\Gr$ and $\Hr$ are Banach spaces then $\Bc(\Gr,\Hr)$ is the
space of bounded linear operators $\Gr\rightarrow\Hr$, the subspace
of compact operators is denoted $\Kc(\Gr,\Hr)$, and we set
$\Bc(\Hr)=\Bc(\Hr,\Hr)$ and $\Kc(\Hr)=\Kc(\Hr,\Hr)$. The domain and
the resolvent set of an operator $S$ will be denoted by $\mathcal
D(S)$ and $\rho(S)$ respectively.  The norm of a Banach space $\Gr$
is denoted by $\|\cdot\|_\Gr$ and we omit the index if the space
plays a central r\^ole.  The adjoint space (space of antilinear
continuous forms) of a Banach space $\Gr$ is denoted $\Gr^*$ and if
$u\in\Gr$ and $v\in\Gr^*$ then we set $v(u)=\langle u,v\rangle$. The
embedding $\Gr\subset\Gr^{**}$ is realized by defining $\langle
v,u\rangle=\overline{\langle u,v\rangle}$.

If $\Gr,\Hr,\Kr$ are Banach spaces such that $\Gr\subset\Hr$
continuously and densely and $\Hr\subset\Kr$ continuously then we
we have a natural continuous embedding
$\Bc(\Hr)\hookrightarrow\Bc(\Gr,\Kr)$ that will be used without
comment later on.

A {\em Friedrichs couple} $(\Gr,\Hr)$ is a pair of Hilbert spaces
$\Gr,\Hr$ together with a continuous dense embedding
$\Gr\subset\Hr$. The {\em Gelfand triplet} associated to it is
obtained by identifying $\Hr=\Hr^*$ with the help of the Riesz
isomorphism and then taking the adjoint of the inclusion map
$\Gr\rarrow\Hr$.  Thus we get $\Gr\subset\Hr\subset\Gr^*$ with
continuous and dense embeddings.  Now if $u\in\Gr$ and
$v\in\Hr\subset\Gr^*$ then $\langle u,v\rangle$ is the scalar
product in $\Hr$ of $u$ and $v$ and also the action of the
functional $v$ on $u$. As noted above, we have
$\Bc(\Hr)\subset\Bc(\Gr,\Gr^*)$.

If $X$ is a locally compact topological space then $B(X)$ is the
$C^*$-algebra of bounded Borel complex functions on $X$, with norm
$\sup_{x\in X}|\varphi(x)|$, and $B_0(X)$ is the subalgebra
consisting of functions which tend to zero at infinity. Then $C(X)$,
$\cb(X)$, $\co(X)$ and $\cc(X)$ are the spaces of complex functions
on $X$ which are continuous, continuous and bounded, continuous and
convergent to zero at infinity, and continuous with compact support
respectively. We denote $\cchi_S$ the characteristic function of a
set $S\subset X$.

\medskip

\noindent{\bf Acknowledgments:} 
We would like to thank Fran\c{c}oise Piquard: several discussions
with her on factorization theorems for Banach space operators have
been very helpful in the context of Section \ref{wvp}. We are also
indebted to Francis Nier for a critical reading of the first version
of this text and for several useful suggestions and to Thierry
Coulhon for a discussion concerning the regularity assumptions from
Theorems \ref{th:wmetric} and \ref{th:wmetric*} and for the
references \cite{AC,ACDH}.

\section{Banach modules and decay preserving operators} \label{pr}
\subsection{Banach modules}\label{bm}
We use the terminology of \cite{FD} but with some abbreviations,
e.g.\ a {\em morphism} is a linear multiplicative map between two
algebras, and a {\em $*$-morphism} is a morphism between two
$*$-algebras which commutes with the involutions. We recall that an
\emph{approximate unit} in a Banach algebra $\Mc$ is a net
$\{J_\alpha\}$ in $\Mc$ such that $\|J_\alpha\|\leq C$ for some
constant $C$ and all $\alpha$ and $\lim_\alpha \|J_\alpha M -
M\|=\lim_\alpha\|MJ_\alpha -M\|=0$ for all $M\in\Mc$.  An
approximate unit exists if and only if there is a number $C$ such
that for each $\varepsilon>0$ and for each finite set $\mathcal
F\subset\Mc$ there is $J\in \Mc$ with $\|J\|\leq C$ and $\|JM
-M\|\leq\varepsilon$, $\|MJ- M\|\leq \varepsilon$ for all
$M\in\mathcal F$.  It is well known that any $C^*$-algebra has an
approximate unit. If $\Hr$ is a Banach space, we shall say that a
Banach subalgebra $\Mc$ of $\Bc(\Hr)$ is \emph{non-degenerate} if
the linear subspace of $\Hr$ generated by the elements $M u$, with
$M\in\Mc$ and $u\in\Hr$, is dense in $\Hr$.

In view of its importance in our paper, we state below the
Cohen-Hewitt factorization theorem \cite[Ch.\ V--9.2]{FD}.
\begin{theorem}\label{t:ch}
  Let $\Cc$ be a Banach algebra with an approximate unit, let $\Er$
  be a Banach space, and let $Q:\Cc\rarrow\Bc(\Er)$ be a continuous
  morphism.  Denote $\Er_0$ the closed linear subspace of $\Er$
  generated by the elements of the form $Q(\varphi)v$ with
  $\varphi\in\Cc$ and $v\in\Er$. Then for each $u\in\Er_0$ there are
  $\varphi\in\Cc$ and $v\in\Er$ such that $u=Q(\varphi)v$.
\end{theorem}

Now we introduce the framework in which we shall work.
\begin{definition}\label{d:cmod}
  A \emph{Banach module} is a couple $(\Hr,\Mc)$ consisting of a
  Banach space $\Hr$ and a non-degenerate Banach subalgebra $\Mc$ of
  $\Bc(\Hr)$ which has an approximate unit.  If $\Hr$ is a Hilbert
  space and $\Mc$ is a $C^*$-algebra of operators on $\Hr$, we say
  that $\Hr$ is a \emph{Hilbert module}.
\end{definition}
We shall adopt the usual \emph{abus de language} and say that $\Hr$
is a Banach module (over $\Mc$). The distinguished subalgebra $\Mc$
will be called \emph{multiplier algebra of $\Hr$} and, when required
by the clarity of the presentation, we shall denote it $\Mc(\Hr)$.
We are only interested in the case when $\Mc$ does not have a unit:
the operators from $\Mc$ are the prototype of decay improving (or
vanishing at infinity) operators, and the identity operator cannot
have such a property.  Note that it is implicit in Definition
\ref{d:cmod} that \emph{if $\Hr$ is a Hilbert module then its
  adjoint space $\Hr^*$ is identified with $\Hr$} with the help of
the Riesz isomorphism.

If $\{J_\alpha\}$ is an approximate unit of $\Mc$, then the density
in $\Hr$ of the linear subspace generated by the elements $Mu$ is
equivalent to
\begin{equation}\label{e:au}
\lim_\alpha\|J_\alpha u-u\|=0 \quad\textrm{for all}\quad u\in\Hr.
\end{equation}
But much more is true:
\begin{equation}\label{e:ch}
u\in\Hr\Rightarrow
u=Mv \textrm{ for some } M\in\Mc\textrm{ and }v\in\Hr.
\end{equation}
This follows from the Cohen-Hewitt theorem, see Theorem \ref{t:ch}.
By using (\ref{e:au}) we could avoid any reference to this result in
our later arguments; this would make them more elementary but less
simple. From Theorem \ref{t:ch} we also get:

\begin{lemma}\label{lm:ch}
Assume that $\Ac$ is a Banach algebra with approximate unit and that
a morphism $\Phi:\Ac\rightarrow\Mc(\Hr)$ with dense image is given.
Then  each $u\in\Hr$ can be
written as $u=Av$ where $A\in\Phi(\Ac)$ and $v\in\Hr$.
\end{lemma}

\begin{example}\label{ex:xmod}{\rm
The simplest example of Banach module is the following. Let $X$
be a locally compact non-compact topological space and let $\Hr$ be
a Banach space. We say that $\Hr$ is a \emph{Banach $X$-module} if a
continuous morphism $Q:\co(X)\rarrow\Bc(\Hr)$ has been given such
that the linear subspace generated by the vectors of the form
$Q(\varphi)u$, with $\varphi\in \co(X)$ and $u\in\Hr$, is dense in
$\Hr$. If $\Hr$ is a Hilbert space and $Q$ is a $*$-morphism, we say
that $\Hr$ is a \emph{Hilbert $X$-module}.
We shall use the notation $\varphi(Q)\equiv Q(\varphi)$.  The Banach
module structure on $\Hr$ is defined by the closure $\Mc$ in
$\Bc(\Hr)$ of the set of operators of the form $\varphi(Q)$ with
$\varphi\in \co(X)$. In the case of a Hilbert $X$-module the closure
is not needed and we get a Hilbert module structure (because a
$*$-morphism between two $C^*$-algebras is continuous and its range
is a $C^*$-algebra). Banach $X$-modules appear
naturally in differential geometry as spaces of sections of vector
fiber bundles over a manifold $X$, and this is the point of interest
for us.  
}\end{example}

\begin{remark}\label{re:xmod}{\rm
In the case of a Banach $X$-module, Lemma \ref{lm:ch} gives:
each $u\in\Hr$ can be written as $u=\psi(Q)v$ with
$\psi\in\co(X)$ and $v\in\Hr$. In particular, we deduce
that \emph{the morphism $Q$ has an extension, also denoted
$Q$, to a unital continuous morphism of $\cb(X)$ into $\Bc(\Hr)$}
which is uniquely determined by the following strong continuity
property: if $\{\varphi_n\}$ is a bounded sequence in $\cb(X)$ such
that $\varphi_n\rarrow\varphi$ locally uniformly, then
$\varphi_n(Q)\rarrow\varphi(Q)$ strongly on $\Hr$. Indeed, 
we can define $\varphi(Q)u=(\varphi\psi)(Q)v$ for each
$\varphi\in\cb(X)$; then if $e_\alpha$ is an approximate unit for
$\co(X)$ with $\|e_\alpha\|\leq1$ we get $\varphi(Q)u=\lim(\varphi
e_\alpha)(Q)u$ hence the definition is independent
of the factorization of $u$ and $\|\varphi(Q)\|\leq
\|Q\|\sup|\varphi|$. 
}\end{remark}

\begin{remark}\label{re:xmod2}{\rm
If $\Hr$ is a Hilbert $X$-module one can extend the morphism even
further: \emph{$Q$ canonically extends to a $*$-morphism
$\varphi\mapsto\varphi(Q)$   of $B(X)$ into $\Bc(\Hr)$
such that\footnote{ If $X$ is
second countable then this property determines uniquely the
extension. In general, uniqueness is assured by the property:
if $U\subset X$ is open then
$\cchi_U(Q)=\sup_{\varphi}\varphi(Q)$, where $\varphi$ runs over the
set of continuous functions with compact support such that $0\leq
\varphi\leq\cchi_U$.
}:
if $\{\varphi_n\}$ is a bounded sequence in $B(X)$ and
$\lim_{n\rightarrow\infty}\varphi_n(x)=\varphi(x)$ for all $x\in X$,
then $\slim_n\varphi_n(Q)=\varphi(Q)$.}
This follows from standard integration theory see \cite{Be,Loo}.
In particular, a separable Hilbert $X$-module is essentially a direct
integral of Hilbert spaces over $X$, see \cite[II.6.2]{Dix}, but we
shall not need this fact.
}\end{remark}

The class of $X$-modules is more general than it appears at first
sight. Indeed, if $\Cc$ is an abelian $C^*$-algebra then one has a
canonical identification $\Cc\equiv\co(\mathcal X)$ where $\mathcal
X$ is the spectrum of $\Cc$.  However, the space $\mathcal X$ is in
general rather complicated so it is not really useful to take it
into account. In particular, this happens in the following class of
examples of interest in applications (see Section \ref{wvp}).

\begin{example}\label{ex:xmodf}
{\rm 
Let $X$ be a set and $\Fc$ a filter on $X$. 
Let us choose a $C^*$-algebra $\Cc$ of bounded complex functions on
$X$ (with the sup norm) and then let $\Cc_0$ be the set of
$\varphi\in\mathcal C$ such that $\lim_\Fc \varphi=0$.
Then $\Cc_0$ is a $C^*$-algebra and its spectrum $\mathcal X$
contains $X$ but is much larger than $X$ in general. 
}\end{example}

Let us say that a Banach module structure defined by a Banach
algebra $\Nc$ on $\Hr$ is \emph{finer} than that defined by $\Mc$ if
$\Mc\subset\Nc$.  In the next example we show that, by using the
same idea as in Example \ref{ex:xmodf}, one can define on each
$X$-module new Banach module structures finer than the initial one.
In Section \ref{wvp} we shall consider the question of the stability
of the essential spectrum in situations of this type, when the
perturbation vanishes at infinity in a weak sense.

\begin{example}\label{ex:fxmod}{\rm
    Let $\Hr$ be a Hilbert $X$-module over a locally compact
    non-compact topological space $X$ and let $\Fc$ be a filter on
    $X$ finer than the Fr\'echet filter. We extend the morphism
    $Q$ to all of $B(X)$ as explained in Remark \ref{re:xmod2} and
    observe that we get a finer Hilbert module structure on $\Hr$
    by taking $\{\vphi(Q)\mid \lim_\Fc\vphi=0\}$ as multiplier
    algebra. One can proceed similarly in the case of a Banach
    $X$-module, it suffices to replace $B(X)$ by $\cb(X)$.
}\end{example}

We now give an example of a non-topological nature.

\begin{example}\label{ex:meas}{\rm
    Let $(X,\mu)$ be a measure space with $\mu(X)=\infty$.  We
    define the class of functions which ``vanish at infinity'' as
    follows.  Let us say that a set $F\subset X$ is of
    \emph{cofinite measure} if its complement $F^\crm$ is of finite
    (exterior) measure.  The family of sets of cofinite measure is
    clearly a filter $\Fc_\mu$. If $\vphi$ is a function on $X$ then
    $\lim_{\Fc_\mu}\vphi=0$ means that for each $\veps>0$ the set
    where $|\vphi(x)|\geq\veps$ is of finite measure. We denote
    $B_\mu(X)$ the $C^*$-subalgebra of $L^\infty(X)$ consisting of
    functions such that $\lim_{\Fc_\mu}\vphi=0$.  Let $\Nc_\mu$ be
    the set of (equivalence classes of) Borel subsets of finite
    measure of $X$.  Then $\{\cchi_N\}_{N\in\Nc_\mu}$ is an
    approximate unit of $B_\mu(X)$ because for each $\vphi\in
    B_\mu(X)$ and each $\varepsilon>0$ we have $N=\{x\mid
    |\vphi(x)|\geq\varepsilon\}\in\Nc_\mu$ and
    $\mbox{ess-sup\,}|\vphi-\cchi_N\vphi|\leq\veps$. Now it is clear
    that $L^2(X)$ and, more generally, any direct integral of
    Hilbert spaces over $X$, has a natural Hilbert module structure
    with $B_\mu(X)$ as multiplier algebra.  }\end{example}

If $\Hr$ is a Banach module and the Banach space $\Hr$ is reflexive
we say that $\Hr$ is a \emph{reflexive Banach module}. In this case
the adjoint Banach space $\Hr^*$ \emph{is equipped with a canonical
Banach module structure}, its multiplier algebra being
$\Mc(\Hr^*):=\{A^*\mid A\in\Mc(\Hr)\}$.  This is a closed subalgebra
of $\Bc(\Hr^*)$ which clearly has an approximate unit and the linear
subspace generated by the elements of the form $A^*v$, with
$A\in\Mc(\Hr)$ and $v\in\Hr^*$, is weak$^*$-dense, hence dense, in
$\Hr^*$.  Indeed, if $u\in\Hr$ and $\langle u,A^*v\rangle=0$ for all
such $A,v$ then $Au=0$ for all $A\in\Mc(\Hr)$ hence $u=0$ because of
(\ref{e:au}).
\begin{example}\label{ex:sob}{\rm
    For each real number $s$ let $\Hr^s:=\Hr^s(\R^n)$ be the Hilbert
    space of distributions $u$ on $\R^n$ such that
    $\|u\|_s^2:=\int(1+|k|^2)^s|\what u(k)|^2 d k<\infty$, where
    $\what u$ is the Fourier transform of $u$. This is the usual
    Sobolev space of order $s$ on $\R^n$. The algebra $\Sr$ of
    Schwartz test functions on $\R^n$ is naturally embedded in
    $\Bc(\Hr^s)$, a function $\varphi\in\Sr$ being identified with
    the operator of multiplication by $\vphi$ on $\Hr^s$.  If we
    denote by $\Mc^s$ the closure of $\Sr$ in $\Bc(\Hr^s)$, then
    clearly $(\Hr^s,\Mc^s)$ is a Banach module and this Banach
    module is a Hilbert module if and only if $s=0$.  The module
    adjoint to $(\Hr^s,\Mc^s)$ is identified with
    $(\Hr^{-s},\Mc^{-s})$. Note that $\Mc^s$ can be realized as a
    subalgebra of $\Mc^0=C_0(\R^n)$, namely $\Mc^s$ is the
    completion of $\Sr$ for the norm
    $\|\varphi\|_{\Mc^s}:=\sup_{\|u\|_s=1}\|\varphi u\|_s$, and then
    we have $\Mc^s=\Mc^{-s}$ isometrically and $\Mc^s\subset\Mc^t$
    if $s\geq t\geq 0$ (by interpolation).  }\end{example}
\begin{definition}\label{d:1}
  A couple $(\Gr,\Hr)$ consisting of a Hilbert module $\Hr$ and a
  Hilbert space $\Gr$ such that $\Gr\subset\Hr$ continuously and
  densely will be called a {\em Friedrichs module}.  If $\Hr$ is a
  Hilbert $X$-module over a locally compact space $X$, we say that
  $(\Gr,\Hr)$ is a \emph{Friedrichs $X$-module}.  If
  $\Mc(\Hr)\subset \Kc(\Gr,\Hr)$, we say that $(\Gr,\Hr)$ is a
  {\em compact Friedrichs module}.
\end{definition}
In the situation of this definition we always identify $\Hr$ with
its adjoint space, which gives us a Gelfand triplet
$\Gr\subset\Hr\subset\Gr^*$.  If $(\Gr,\Hr)$ is a compact Friedrichs
module then each operator $M$ from $\Mc(\Hr)$ extends to a compact
operator $M:\Hr\rightarrow\Gr^*$ (this is the adjoint of the compact
operator $M^*:\Gr\rightarrow\Hr$).  Thus we shall have
$\Mc(\Hr)\subset\Kc(\Gr,\Hr)\cap\Kc(\Hr,\Gr^*)$.

\begin{example}\label{ex:sob*}{\rm
    With the notations of Example \ref{ex:sob}, if we set
    $\Hr=\Hr^0$ and take $s>0$, then $(\Hr^s,\Hr)$ is a compact
    Friedrichs module and the associated Gelfand triplet is
    $\Hr^s\subset\Hr\subset\Hr^{-s}$. Indeed, if $\varphi\in
    C_0(\R^n)$ then the operator of multiplication by $\vphi$ is a
    compact operator $\Hr^s\rarrow\Hr$.  }\end{example}

\subsection{Decay improving operators}

Let $\Hr$ and $\Kr$ be Banach spaces. If $\Kr$ is a Banach module
then we shall denote by $\bol(\Hr,\Kr)$ the norm closed linear
subspace generated by the operators $MT$, with $T\in\Bc(\Hr,\Kr)$
and $M\in\Mc(\Kr)$.  We say that an operator in $\bol(\Hr,\Kr)$ is
{\em decay improving}, or {\em left vanishes at infinity} (with
respect to $\Mc(\Kr)$, if this is not obvious from the context). If
$J_\alpha$ is an approximate unit for $\Mc(\Kr)$, then for an
operator $S\in\Bc(\Hr,\Kr)$ we have:
\begin{eqnarray}\label{e:lv}
S\in\bol(\Hr,\Kr) & \Leftrightarrow &
\lim_\alpha\|J_\alpha S-S\|=0 \\
& \Leftrightarrow &
S=MT\textrm{ for some }M\in\Mc(\Kr) \textrm{ and }T\in \Bc(\Hr, \Kr).
\nonumber
\end{eqnarray}
The second equivalence  follows from the Cohen-Hewitt theorem
(Theorem \ref{t:ch}).

If $\Hr$ is a Banach module then one can similarly define
$\Bc_0^{\,r}(\Hr,\Kr)$ as the norm closed linear subspace generated
by the operators $TM$ with $T\in\Bc(\Hr,\Kr)$ and $M\in\Mc(\Hr)$.
We say that the elements of $\bor(\Hr,\Kr)$ \emph{right vanish at
infinity}.  If both $\Hr$ and $\Kr$ are Banach modules then both
spaces $\bol(\Hr,\Kr)$ and $\bor(\Hr,\Kr)$ make sense and we set
$\Bc_0(\Hr,\Kr)=\bol(\Hr,\Kr)\cap\bor(\Hr,\Kr)$.

Some simple properties of these spaces are described below.
\begin{proposition}\label{p:ref}
If $\Kr$ is a reflexive Banach module and $S\in\Bc_0^{\,l}(\Hr,\Kr)$
then $S^*$ belongs to $\Bc_0^{\,r}(\Kr^*,\Hr^*)$.
\end{proposition}
\proof We have $S=MT$ with $M\in\Mc(\Kr)$ and $T\in\Bc(\Hr,\Kr)$ by
(\ref{e:lv}), which implies $S^*=T^*M^*$ and we have
$M^*\in\Mc(\Kr^*)$ by definition.
\qed
\begin{proposition}
If $\Hr$ is a Hilbert module then $\Bc_0(\Hr)$ is a
$C^*$-algebra and  an operator $S\in\Bc(\Hr)$ belongs to it if and
only if one can write $S=MTN$ with $M,N\in\Mc(\Hr)$ and
$T\in\Bc(\Hr)$. 
\end{proposition}
\proof $\Bc_0(\Hr)$ is clearly a $C^*$-algebra, so if
$S\in\Bc_0(\Hr)$ then $S=S_1S_2$ for some operators
$S_1,S_2\in\Bc_0(\Hr)$. 
Thus $S_1=MT_1$ and $S_2=T_2N$ for some $M,N\in\Mc(\Hr)$ and
$T_1,T_2\in\Bc(\Hr)$, hence $S=MT_1T_2N$.
\qed
\begin{proposition}\label{p:cvv}
  If $\Kr$ is a Banach module then
  $\Kc(\Hr,\Kr)\subset\Bc_0^{\,l}(\Hr,\Kr)$.  If $\Hr$ is a
  reflexive Banach module, then
  $\Kc(\Hr,\Kr)\subset\Bc_0^{\,r}(\Hr,\Kr)$.
\end{proposition}
\proof If $\{J_\alpha\}$ is an approximate unit for $\Mc(\Kr)$ then
$\slim_\alpha J_\alpha u=u$ uniformly in $u$ if $u$ belongs to a
compact subset of $\Kr$. Hence if $S\in\Kc(\Hr,\Kr)$ then
$\lim_\alpha\|J_\alpha S-S\|=0$ and thus $S\in\bol(\Hr,\Kr)$ by
(\ref{e:lv}). To prove the second part of the proposition, note that
if $S\in\Kc(\Hr,\Kr)$ then $S^*\in\Kc(\Kr^*,\Hr^*)$, so
$S^*\in\bol(\Kr^*,\Hr^*)$ by what we just proved, hence
$S^{**}\in\bor(\Hr,\Kr^{**})$ by Proposition \ref{p:ref}. Thus we
get $\lim_\alpha\|S^{**}J_\alpha -S^{**}\|=0$ if $\{J_\alpha\}$ is
an approximate unit for $\Mc(\Hr)$. But clearly $S=S^{**}$, hence
$S\in\bor(\Hr,\Kr)$.  \qed
\begin{proposition}\label{p:comp}
Let $\Hr$ be a Banach module and $\Gr$ a Banach space continuously
embedded in $\Hr$ and such that for each $M\in\Mc(\Hr)$ the
restriction of $M$ to $\Gr$ is a compact operator
$\Gr\rarrow\Hr$. If $R\in\bol(\Hr)$ and $R\Hr\subset\Gr$, then 
$R\in\Kc(\Hr)$. 
\end{proposition}
\proof According to (\ref{e:lv}) we have $R=\lim_\alpha J_\alpha R$,
the limit being taken in norm. But $R\in\Bc(\Hr,\Gr)$ by the closed
graph theorem and $J_\alpha\in\Kc(\Gr,\Hr)$ by hypothesis, so that
$J_\alpha R\in\Kc(\Hr)$.  \qed

\subsection{Decay preserving operators}\label{ql}
\begin{definition}\label{d:ql}
  Let $\Hr$, $\Kr$ be Banach modules and let $S\in\Bc(\Hr,\Kr)$.  We
  say that $S$ is \emph{left decay preserving } if for each
  $M\in\Mc(\Hr)$ we have $SM\in\Bc_0^{\,l}(\Hr, \Kr)$.  We say that
  $S$ is \emph{right decay preserving } if for each $M\in\Mc(\Kr)$
  we have $MS\in\Bc_0^{\,r}(\Hr, \Kr)$.  If $S$ is left and right
  decay preserving, we say that $S$ is \emph{decay preserving}.
\end{definition}
We denote $\Bc_q^{\,l}(\Hr, \Kr)$, $\Bc_q^{\,r}(\Hr, \Kr)$ and
$\Bc_q(\Hr, \Kr)$ these classes of operators (the index q comes from
\emph{quasilocal}, a terminology which is sometimes more convenient
than ``decay preserving''). These are closed subspaces of
$\Bc(\Hr,\Gr)$. The next result is obvious; a similar assertion
holds in the right decay preserving case.
\begin{proposition}\label{p:dql}
  Let $\{J_\alpha\}$ be an approximate unit for $\Mc(\Hr)$ and let
  $S$ be an operator in $\Bc(\Hr, \Kr)$. Then $S$ is left decay
  preserving   if and only if one of the following conditions is
  satisfied:\\ 
  {\rm(1)} $SJ_\alpha \in\Bc_0^{\,l}(\Hr,\Kr)$  for all $\alpha$.\\
  {\rm(2)} for each $M\in\Mc(\Hr)$ there are $T\in\Bc(\Hr, \Kr)$ and
  $N\in\Mc(\Kr)$ such that $SM=NT$.
\end{proposition}

The next proposition, which says that the decay preserving property
is stable under the usual algebraic operations, is an immediate
consequence of Proposition \ref{p:dql}.  There is, of course, a
similar statement with ``left'' and ``right'' interchanged.  We
denote by $\Gr,\Hr$ and $\Kr$ Banach modules.
\begin{proposition}\label{p:aql}
{\rm (1)} $S\in\bql(\Hr,\Kr)$ and $T\in\bql(\Gr,\Hr)$ $\Rightarrow$
$ST\in\bql(\Gr,\Kr)$.\\
{\rm (2)} If $\Hr,\Kr$ are reflexive and $S\in\bql(\Hr,\Kr)$, then
$S^*\in\bqr(\Kr^*,\Hr^*)$.\\
{\rm (3)} If $\Hr$ is a Hilbert module then $\Bc_q(\Hr)$ is a unital
$C^*$-subalgebra of $\Bc(\Hr)$.
\end{proposition}

Obviously $\bol(\Hr,\Kr)\subset\bql(\Hr,\Kr)$ and
$\bor(\Hr,\Kr)\subset\bqr(\Hr,\Kr)$ but this fact is of no interest.
The main results of this paper depend on finding other, more
interesting examples of decay preserving operators. We shall give in
the rest of this subsection some elementary examples of such
operators and in Subsections \ref{lcg} and \ref{wvpr} more subtle
ones.

\medskip

From now on in this subsection $X$ will be a locally compact
non-compact topological space.  The \emph{support} $\supp u\subset
X$ of an element $u$ of a Banach $X$-module $\Hr$ is defined as the
smallest closed set such that its complement $U$ has the property
$\varphi(Q)u=0$ if $\varphi\in\cc(U)$.  Clearly, the set $\Hr_{\rm
  c}$ of elements $u\in\Hr$ such that $\supp u$ is compact is a
dense subspace of $\Hr$.

Let $\Hr,\Kr$ be Banach $X$-modules, let $S\in \Bc(\Hr, \Kr)$, and
let $\varphi, \psi\in C(X)$, not necessarily bounded.  We say that
$\varphi(Q) S \psi(Q)$ is a bounded operator if there is a constant
$C$ such that
$$\|\xi(Q)\varphi(Q)S\psi(Q)\eta(Q)\|\leq C\sup|\xi|\sup|\eta|$$
for
all $\xi,\eta\in\cc(X)$. The lower bound of the admissible constants
$C$ in this estimate is denoted $\|\varphi(Q)S\psi(Q)\|$.  If $\Kr$
is a reflexive Banach $X$-module, then the product
$\varphi(Q)S\psi(Q)$ is well defined as sesquilinear form on the
dense subspace $\Kr^*_{\rm c}\times\Hr_{\rm c}$ of $\Kr^*\times\Hr$
and the preceding boundedness notion is equivalent to the continuity
of this form for the topology induced by $\Kr^*\times\Hr$.  We
similarly define the boundedness of the commutator $[S,
\varphi(Q)]$.
\begin{proposition}\label{prop26}
  Assume that $S\in \Bc(\Hr, \Kr)$ and let $\theta : X\rightarrow[1,
  \infty[$ be a continuous function such that
  $\lim_{x\rightarrow\infty}\theta(x)=\infty$.  If
  $\theta(Q)S\theta^{-1}(Q)$ is a bounded operator then $S$ is left
  decay preserving.  If $\theta^{-1}(Q)S\theta(Q)$ is a bounded
  operator then $S$ is right decay preserving.
\end{proposition}
\proof Let $K\subset X$ be compact, let $U\subset X$ be a
neighbourhood of infinity in $X$, and let $\vphi,\psi\in\cb(X)$ such
that $\supp\vphi\subset K,\ \supp\psi\subset U$ and
$|\vphi|\leq1,|\psi|\leq1$.  Then $\theta\vphi$ is a bounded
function and $\psi\theta^{-1}$ is bounded and can be made as small
as we wish by choosing $U$ conveniently. Thus given $\varepsilon>0$
we have
$$
\|\psi(Q)S\vphi(Q)\|\leq\|\psi\theta^{-1}\|\cdot
\|\theta(Q)S\theta^{-1}(Q)\|\cdot
\|\theta\vphi\|\leq \varepsilon
$$
if $U$ is a sufficiently small neighbourhood of infinity. Then the
result follows from Proposition \ref{p:dql}(1) and relation
(\ref{e:lv}). 
\qed

The boundedness of $\theta(Q)S\theta^{-1}(Q)$ can be checked by
estimating the commutator $[S,\theta(Q)]$; we give an example for
the case of metric spaces.  Note that on metric spaces one has a
natural class of regular functions, namely the Lipschitz functions,
for example the functions which give the distance to subsets:
$\rho_K(x)=\inf_{y\in K}\rho(x,y)$ for $K\subset X$.

We say that a locally compact metric space $(X,\rho)$ is
\emph{proper} if the metric $\rho$ has the property
$\lim_{y\rightarrow\infty} \rho(x,y)=\infty$ for some (hence for
all) points $x\in X$.  Equivalently, if $X$ is not compact but the
closed balls are compact.

\begin{corollary}\label{ex27}
  Let $(X,\rho)$ be a proper locally compact metric space.  If $S$
  belongs to $\Bc(\Hr, \Kr)$ and if $[S,\theta(Q)]$ is bounded for
  each positive Lipschitz function $\theta$, then $S$ is decay
  preserving. 
\end{corollary}
\proof Indeed, by taking $\theta=1+\rho_K$ and by using the
notations of the proof of Proposition \ref{prop26}, we easily get
the following estimate: there is $C<\infty$ depending only on $K$
such that
\begin{equation*}
\|\vphi(Q)S\psi(Q)\|\leq C(1+\rho(K,U))^{-1},
\end{equation*}
where $\rho(K,U)$ is the distance from $K$ to $U$.
Since $S^*$ has the same properties as $S$, this proves that $S$ is
decay preserving. Note that the boundedness of
$[S, \rho_x(Q)]$ for some $x\in X$ suffices in this argument.
\qed

\section{Compact perturbations in Banach modules} \label{mr}

In this section $(\Gr,\Hr)$ will be a \emph{compact
Friedrichs module} in the sense of Definition \ref{d:1}. As usual, we
associate to it a Gelfand triplet $\Gr\subset\Hr\subset\Gr^*$ and we
set $\|\cdot\|=\|\cdot\|_\Hr$. We are interested in criteria which
ensure that an operator $B$ is a compact perturbation of an operator
$A$, both operators being unbounded operators in $\Hr$ obtained as
restrictions of some bounded operators $\Gr\rarrow\Gr^*$.  More
precisely, the following is a general assumption (suggested by the
statement of Theorem 2.1 in \cite{OS}) which will always be
fulfilled: \label{AB}

$$
\leqno{\mbox{\bf($AB$)}}\hspace{1mm}
\left\{
\begin{array}{ll}
A, B \mbox{ are  closed densely defined operators in } \Hr
\mbox{ with }
\rho(A)\cap\rho(B)\neq\emptyset \\
\mbox{ and having the following properties: }
\Dc(A)\subset\Gr \mbox{ densely, } \Dc(A^*)\subset\Gr, \\
\Dc(B)\subset\Gr \mbox{ and } A, B
\mbox{ extend to continuous operators }
\widetilde{A}, \widetilde{B}\in \Bc(\Gr, \Gr^*).
\end{array}
\right.
$$
\begin{example}\label{e:ab}{\rm
    One can construct interesting classes of operators with the
    properties required in (AB) as follows.  Let $\Gr_a$, $\Gr_b$ be
    Hilbert spaces such that $\Gr\subset\Gr_a\subset\Hr$ and
    $\Gr\subset\Gr_b\subset\Hr$ continuously and densely. Thus we
    have two scales
\begin{eqnarray*}
\Gr\subset\Gr_a\subset \!\!\! & \Hr & \!\!
\!\subset\Gr_a^*\subset\Gr^*,\\ 
\Gr\subset\Gr_b\subset \!\!\! & \Hr & \!\!
\!\subset\Gr_b^*\subset\Gr^*. 
\end{eqnarray*}
Then let $A_0\in \Bc(\Gr_a, \Gr_a^*)$ and $B_0\in \Bc(\Gr_b,
\Gr_b^*)$ such that $A_0-z:\Gr_a\rightarrow\Gr_a^*$ and $B_0-z:
\Gr_b\rightarrow\Gr_b^*$ are bijective for some number $z$.
According to Lemma \ref{prop33} we can associate to $A_0,B_0$ closed
densely defined operators $A=\widehat{A_0},B=\widehat{B_0}$ in
$\Hr$, such that the domains $\Dc(A)$ and $\Dc(A^*)$ are dense
subspaces of $\Gr_a$ and the domains $\Dc(B)$ and $\Dc(B^*)$ are
dense subspaces of $\Gr_b$.  If we also have $\Dc(A)\subset\Gr$
densely, $\Dc(A^*)\subset\Gr$ and $\Dc(B)\subset\Gr$, then all the
conditions of the assumption (AB) are fulfilled with $\wtilde
A=A_0|\Gr$ and $\wtilde B=B_0|\Gr$. Such a construction will
be used in Corollary \ref{c:main-aa}. }\end{example}

The r\^ole of the assumption (AB) is to allow us to give a rigorous
meaning to the formal relation, where $z\in\rho(A)\cap\rho(B)$,
\begin{equation}\label{33}
(A-z)^{-1}-(B-z)^{-1}=(A-z)^{-1}(B-A)
(B-z)^{-1}.
\end{equation}
Recall that $z\in\rho(A)$ if and only if $\bar{z}\in\rho(A^*)$ and
then $(A^*-\bar z)^{-1}=(A-z)^{-1*}$.  Thus we have
$(A-z)^{-1*}\Hr\subset\Gr$ by the assumption (AB) and this allows
one to deduce that ($A-z)^{-1}$ extends to a unique continuous
operator $\Gr^*\rarrow\Hr$, that we shall denote for the moment by
$R_z$. From $R_z(A-z)u=u$ for $u\in \Dc(A)$ we get, by density of
$\Dc(A)$ in $\Gr$ and continuity, $R_z(\widetilde{A}-z)u=u$ for
$u\in\Gr$, in particular
\begin{equation*}
(B-z)^{-1}= R_z(\widetilde{A}-z)(B-z)^{-1}.
\end{equation*}
On the other hand, the identity
\begin{equation*}
(A-z)^{-1}= (A-z)^{-1}(B-z)(B-z)^{-1}=
R_z(\widetilde{B}-z)(B-z)^{-1}
\end{equation*}
is trivial. Subtracting the last two relations we get
\begin{equation*}
(A-z)^{-1}-(B-z)^{-1}=R_z(\widetilde{B}-\widetilde{A})(B-z)^{-1}
\end{equation*}
Since $R_z$ is uniquely determined as extension of $(A-z)^{-1}$ to a
continuous map $\Gr^*\rarrow\Hr$, we shall keep the notation
$(A-z)^{-1}$ for it. With this convention, the rigorous version of
(\ref{33}) that we shall use is:
\begin{equation}\label{34}
(A-z)^{-1}-(B-z)^{-1}=(A-z)^{-1}(\wtilde B-\wtilde A)(B-z)^{-1}.
\end{equation}

\begin{theorem}\label{t:main}
  Let $A,B$ satisfy assumption (AB) and let us assume that there are
  a Banach module $\Kr$ and operators $S\in \Bc(\Kr,\Gr^*)$ and
  $T\in\bol(\Gr,\Kr)$ such that $\wtilde B-\wtilde A=ST$ and
  $(A-z)^{-1}S\in\bql(\Kr,\Hr)$ for some $z\in\rho(A)\cap\rho(B)$.
  Then the operator $B$ is a compact perturbation of the operator
  $A$ and $\se(B)=\se(A)$.
\end{theorem}
\proof It suffices to show that
$R\equiv(A-z)^{-1}-(B-z)^{-1}\in\bol(\Hr)$, because the domains of
$A$ and $B$ are included in $\Gr$, hence $R\Hr\subset\Gr$, which
finishes the proof because of Proposition \ref{p:comp}.  Now due to
(\ref{34}) and to the factorization assumption, we can write $R$ as
a product $R=[(A-z)^{-1}S][T(B-z)^{-1}]$ where the first factor is
in $\bql(\Kr,\Hr)$ and the second in $\bol(\Hr,\Kr)$, so the product
is in $\bol(\Hr)$.  \qed

\begin{remarks}\label{r:1}{\rm
    (1) We could have stated the assumptions of Theorem \ref{t:main}
    in an apparently more general form, namely
    $B-A=\sum_{k=1}^nS_kT_k$ with operators $S_k\in\Bc(\Kr_k,\Gr^*)$
    and $T_k\in \Bc(\Gr,\Kr_k)$. But we are reduced to the stated
    version of the assumption by considering the Hilbert
    module $\Kr=\oplus\Kr_k$ and $S=\oplus S_k,T=\oplus T_k$.\\
    (2) If $V\in\Kc(\Gr,\Gr^*)$ and if $\Kr$ is an infinite
    dimensional module, then there are operators
    $S\in\Bc(\Kr,\Gr^*)$ and $T\in\Kc(\Gr,\Kr)$ such that $V=ST$
    (the proof is an easy exercise). This and the preceding remark
    show that compact contributions to $\wtilde B-\wtilde A$ are
    trivially covered by the factorization assumption.
  }\end{remarks}

If $A$ is self-adjoint then the conditions on $A$ in assumption (AB)
are satisfied if $\Dc(A)\subset\Gr\subset \Dc(|A|^{1/2})$ densely
(see the Appendix). Moreover, if $A$ is semibounded, then this
condition is also necessary. In particular, we have:
\begin{corollary}\label{c:main-aa}
Let $A, B$ be self-adjoint operators on $\Hr$ such that
\begin{equation*}
\Dc(A)\subset\Gr\subset\Dc(|A|^{1/2}) \mbox{\rm\ \  and\ \ }
\Dc(B)\subset\Gr\subset \Dc(|B|^{1/2}) \mbox{\rm\ \  densely.}
\end{equation*}
Let $\widetilde{A}, \widetilde{B}$ be the unique extensions of $A,
B$ to operators in $\Br(\Gr,\Gr^* )$. Assume that there is a Hilbert
module $\Kr$ and that $\widetilde{B}-\widetilde{A}=S^*T$ for some
$S\in \Bc(\Gr,\Kr)$ and $T\in\bol(\Gr,\Kr)$ such that
$S(A-z)^{-1}\in\bqr(\Hr,\Kr)$ for some $z\in\rho(A)\cap\rho(B)$.
Then $B$ is a compact perturbation of $A$ and
$\sigma_{\rm ess}(B)=\sigma_{\rm ess}(A)$.
\end{corollary}

The next theorem is convenient for applications to differential
operators in divergence form. Observe that if $(\Er,\Kr)$ is a
Friedrichs module then $\Bc(\Kr)\subset\Bc(\Er,\Er^*)$ hence we can
define
\begin{equation}\label{e:ool}
\Bc_{00}^{\,l}(\Er,\Er^*)=
\textrm{ norm closure of }
\Bc_0^{\,l}(\Kr)\textrm{ in } \Bc(\Er,\Er^*).
\end{equation}
We shall use the terminology and the facts established in the
Appendix, in particular Lemma \ref{prop33}: the operators $D^*aD$
and $D^*bD$ considered below belong to $\Bc(\Gr,\Gr^*)$ and we
denote by $\Delta_a$ and $\Delta_b$ the operators on $\Hr$
associated to them.
\begin{theorem}\label{t:main-2}
  Let $(\Er,\Kr)$ be an arbitrary Friedrichs module and let
  $D\in\Bc(\Gr,\Er)$, $a,b\in\Bc(\Er,\Er^*)$ and $z\in\C$ such
  that:\\ 
  {\rm (1)} The operators $D^*aD-z$ and $D^*bD-z$ are bijective maps
  $\Gr\rightarrow\Gr^*$,\\
  {\rm (2)}
  $a-b\in\Bc_{00}^{\,l}(\Er,\Er^*)$,\\
  {\rm (3)}
  $D(\Delta_a^*-\bar z)^{-1}\in\bqr(\Hr,\Kr)$.\\
  Then $\Delta_b$ is a compact perturbation of $\Delta_a$.
\end{theorem}
\proof We give a proof independent of Theorem \ref{t:main}, although
we could apply this theorem. Clearly $\Delta_a-z$ and
$\Delta_b-z$ extend to bijections $\Gr\rightarrow\Gr^*$ and the
identity
$$
R:=(\Delta_a-z)^{-1}-(\Delta_b-z)^{-1}=(\Delta_a-z)^{-1}
D^*(b-a)D(\Delta_b-z)^{-1}
$$
holds in $\Bc(\Gr^*,\Gr)$, hence in $\Bc(\Hr)$.  Since the
domains of $\Delta_a$ and $\Delta_b$ are included in $\Gr$, we have
$R\Hr\subset\Gr$. Thus, according to Proposition \ref{p:comp}, it
suffices to show that $R\in\Bc^{\,l}_0(\Hr)$.  Since the space
$\Bc_0^{\,l}(\Hr)$ is norm closed and since by hypothesis we can
approach $b-a$ in norm in $\Bc(\Er,\Er^*)$ by operators in
$\Bc_0^{\,l}(\Kr)$, it suffices to show that
$$
(D(\Delta_a^*-\bar z)^{-1})^*cD(\Delta_b-z)^{-1}\in\Bc_0^{\,l}(\Hr)
$$
if $c\in\Bc_0^{\,l}(\Kr)$. But this is clear because
$cD(\Delta_b-z)^{-1}\in\bol(\Hr,\Kr)$ and
$(D(\Delta_a^*-\bar z)^{-1})^*\in\bql(\Kr,\Hr)$
by Proposition \ref{p:aql}.
\qed

The spaces $\Bc_{00}^{\,r}(\Er,\Er^*)$ and $\Bc_{00}(\Er,\Er^*)$ are
defined in an obvious way and we have
\begin{equation}\label{e:coo}
\Kc(\Er,\Er^*)\subset
\Bc_{00}(\Er,\Er^*)
\end{equation}
because $\Kc(\Kr)$ is a dense subset of $\Kc(\Er,\Er^*)$ and
$\Kc(\Kr)\subset\Bc_0(\Kr)$. So we could assume
$a-b\in\Kc(\Er,\Er^*)$, but this case is trivial from the point of
view of this paper.  Although the space $\bool(\Er,\Er^*)$ is much
larger than $\Kc(\Er,\Er^*)$, it is not satisfactory in some
applications, cf.\ Example \ref{ex:sob**} below and Remark
\ref{r:00-0}.  However, we can allow still more general
perturbations and obtain more explicit results if we impose more
structure on the modules.  In Subsection \ref{cpmag} we describe
such improvements for a class of Banach modules over abelian groups.

\begin{example}\label{ex:sob**}{\rm
    In the context of Example \ref{ex:sob*} we may consider the two
    classes of operators $\Bc_{00}^{\,l}(\Hr^s,\Hr^{-s})$ and
    $\bol(\Hr^s,\Hr^{-s})$.  The first space is included in the
    second one and the inclusion is strict, for example
    $\Bc_{00}^{\,l}(\Hr^s,\Hr^{-s})$ does not contain operators of
    order $2s$, while $\bol(\Hr^s,\Hr^{-s})$ contains such
    operators.
}\end{example}

The only condition of Theorem \ref{t:main-2} which, in some concrete
situations, is not easy to check is condition (3). We now give a
perturbative method for checking it.

For the rest of this section we fix two Friedrichs modules
$(\Gr,\Hr)$ and $(\Er,\Kr)$ and a continuous operator
$D:\Gr\rarrow\Er$. Let $a\in\Bc(\Er,\Er^*)$ such that the operator
$D^*aD$ is coercive (see the Appendix), more precisely we have
\begin{equation}\label{eq:peco}
\Re\langle Du,aDu\rangle\geq\mu\|u\|_{\Gr}^2-\nu\|u\|^2_\Hr
\end{equation}
for some strictly positive constants $\mu,\nu$ and all $u\in\Gr$
Then, as explained in the Appendix, if $\Re z\leq-\nu$ the
operator $D^*aD-z$ is a bijective map $\Gr\rarrow\Gr^*$ and
\begin{equation}\label{eq:peco1}
\|(D^*aD-z)^{-1}\|_{\Bc(\Gr^*,\Gr)}\leq\mu^{-1}.
\end{equation}
Note that $a^*$ has all these properties too so the closed densely
defined operators $\Delta_a$ and $\Delta_{a^*}$ in $\Hr$ are well
defined, their domains are dense subsets of $\Gr$, and we have
$\Delta_a^*=\Delta_{a^*}$. It is easy to check that
$\|(\Delta_a-z)^{-1}\|_{\Bc(\Hr)}\leq|\Re z+\nu|^{-1}$ if
$\Re z+\nu<0$. Since $a$ and $a^*$ play a symmetric role, it will
suffice to consider $\Delta_a-z$ in place of $\Delta^*_a-\bar z$
in condition (3) of Theorem \ref{t:main-2}.

Now let $c$ be a second operator with the same properties as
$a$. We assume, without loss of generality, that it satisfies an
estimate like (\ref{eq:peco}) with the same constants $\mu,\nu$.

\begin{proposition}\label{pr:pco}
Assume that
$$
D(\Delta_c-z)^{-1}\in\bqr(\Hr,\Kr)
\hspace{2mm} \mbox{\rm and }\hspace{2mm}
D(D^*cD-z)^{-1}D^*\in\bqr(\Kr)
$$
for some $z$ with $\Re z\leq-\nu$.
If $a-c\in\bqr(\Kr)$ then
$$
D(\Delta_a-z)^{-1}\in\bqr(\Hr,\Kr)
\hspace{2mm} \mbox{\rm and }\hspace{2mm}
D(D^*aD-z)^{-1}D^*\in\bqr(\Kr).
$$
A similar assertion holds for the spaces $\bql$.
\end{proposition}
\proof Let $V=D^*(a-c)D$ and $L_t=(1-t)D^*cD+tD^*aD=D^*cD+tV$.
For $z$ as in the statement of the proposition we have
$\Re\langle u,(L_t-z)u\rangle\geq\mu\|u\|_{\Gr}^2$
if $0\leq t \leq 1$. Hence there is $\varepsilon>0$ such that
$\Re\langle u,(L_t-z)u\rangle\geq\mu/2\|u\|_{\Gr}^2$
if $-\varepsilon \leq t \leq 1+\varepsilon $, in particular
$\|(L_t-z)^{-1}\|_{\Bc(\Gr^*,\Gr)}\leq 2/\mu$ for all such $t$.
If $-\varepsilon \leq s \leq 1+\varepsilon $ and
$|t-s|\|V(L_t-z)^{-1}\|_{\Bc(\Gr^*,\Gr)}<1$ we get a norm convergent
expansion in $\Bc(\Gr^*,\Gr)$
$$
(L_t-z)^{-1}=(L_s-z-(s-t)V)^{-1}=\sum_{k\geq0}(s-t)^k
(L_s-z)^{-1}[V(L_s-z)^{-1}]^k
$$
so the map $t\mapsto (L_t-z)^{-1}\in\Bc(\Gr^*,\Gr)$ is real
analytic on the interval $]-\varepsilon, 1+\varepsilon[$. Let us
denote $\Delta_t$ the operator in$\Hr$ associated to $L_t$ then we
see that the maps $t\mapsto D(\Delta_t-z)^{-1}\in\Bc(\Hr,\Kr)$ and
$t\mapsto D(L_t-z)^{-1}D^*\in\Bc(\Kr)$ are real analytic on the same
interval.  The set of decay preserving operators is a closed
subspace of the Banach space $\Bc(\Hr,\Kr)$ and an analytic function
which on a nonempty open set takes values in a closed subspace
remains in that subspace for ever. Thus it suffices to show that
$D(\Delta_t-z)^{-1}\in\bqr(\Hr,\Kr)$ for small positive values of
$t$. Similarly, we need to prove $D(L_t-z)^{-1}D^*\in\Bc(\Kr)$ only
for small $t$. To prove the first assertion for example, we take
$s=0$ above and get a norm convergent series in $\Bc(\Hr,\Kr)$:
$$
D(L_t-z)^{-1}=\sum_{k\geq0}(-t)^k
D(D^*cD-z)^{-1}[D^*(a-c)D(D^*cD-z)^{-1}]^k.
$$
It is clear that each term belongs to $\bqr(\Hr,\Kr)$.
\qed

\section{Banach modules over abelian groups}\label{lcg}
\subsection{$X$-modules over locally compact abelian groups}

Since a locally compact abelian group $X$ is a locally
compact space, we can consider $X$-modules in the sense of Example
\ref{ex:xmod}. However, the group structure of $X$ allows us to
associate to it more interesting classes of Banach modules
that we shall also call $X$-modules.  Whenever necessary in order to
avoid ambiguities we shall speak of \emph{$X$-module over the
  topological space $X$} if we have in mind the context of Example
\ref{ex:xmod} and of \emph{$X$-module over the group $X$} when we
refer to the structure introduced in the next Definition \ref{d:xm}.

In this section we fix a locally compact non-compact abelian group
$X$ with the group operation denoted additively.  For example, $X$
could be $\R^n,\Z^n$, or a finite dimensional vector space over a
local field, e.g.\ over the field of p-adic numbers.  Let $X^*$ be
the abelian locally compact group dual to $X$.
\begin{definition}\label{d:xm}
A \emph{Banach $X$-module} over the group $X$ is a Banach space
$\Hr$ equipped with a strongly continuous representation $\{V_k\}$
of $X^*$ on $\Hr$.
\end{definition}
Note that we shall use the same notation $V_k$ for the
representations of $X^*$ in different spaces $\Hr$ whenever this
does not lead to ambiguities.

Such \emph{a Banach $X$-module has a canonical structure of Banach
module} that we now define.  We choose Haar measures $dx$ and $dk$
on $X$ and $X^*$ normalized by the following condition: if the
Fourier transform of a function $\varphi$ on $X$ is given by
$(\mathcal F\varphi)(k)\equiv\what\varphi(k)=
\int_{X}\overline{k(x)}\varphi(x)dx$ then $
\varphi(x)=\int_{X^*}k(x)\widehat{\varphi}(k)dk$.  Recall that
$X^{**}=X$. Let $C^{(a)}(X):=\Fc L^1_{\rm c}(X^*)$ be the set of
Fourier transforms of integrable functions with compact support on
$X^*$. It is easy to see that $C^{(a)}(X)$ is a $*$-algebra for the
usual algebraic operations; more precisely, it is a dense subalgebra
of $\co(X)$ stable under conjugation.  For $\varphi \in C^{(a)}(X)$
we set
\begin{equation}\label{e:group1}
\varphi(Q)=\int_{X^*}V_k\what\varphi(k) dk.
\end{equation}
This definition is determined by the formal requirement $k(Q)=V_k$.
Then
\begin{equation}\label{e:group2}
\Mc:=\mbox{ norm closure of }\{\varphi(Q)\mid \varphi\in C^{(a)}(X)\}
\mbox{ in } \Bc(\Hr)
\end{equation}
is a Banach subalgebra of $\Bc(\Hr)$. By using the next lemma we see 
that the couple $(\Hr,\Mc)$ satisfies the conditions of Definition
\ref{d:cmod}, which gives us the canonical Banach module structure
on $\Hr$.
\begin{lemma}\label{l:group3}
The algebra $\Mc$ has an approximate unit consisting of elements of
the form $e_\alpha(Q)$ with $e_\alpha\in C^{(a)}(X)$.
\end{lemma}
\proof Let us fix a compact neighborhood $K$ of the identity
in $X^*$. The set of compact neighborhoods of the identity $\alpha$
such that $\alpha\subset K$ is ordered by $\alpha_1\geq\alpha_2
\Leftrightarrow \alpha_1\subset\alpha_2$. For each such $\alpha$
define $e_\alpha$ by $\what e_\alpha=\cchi_\alpha/|\alpha|$, where
$|\alpha|$ is the Haar measure of $\alpha$. Then
$\|e_\alpha(Q)\|\leq \sup_{k\in K}\|V_k\|< \infty$, from which it is
easy to infer that $\lim_\alpha
\|e_\alpha(Q)\varphi(Q)-\varphi(Q)\|=0$ for all $\varphi\in
C^{(a)}(X)$.  \qed

\begin{example}\label{ex:sob&}{\rm
    Let $X=\R^n$ with the additive group structure and let $\Hr$ be
    the Sobolev space $\Hr^s(X)$ for some real number $s$.  We
    identify as usual $X^*$ with $X$ by setting $k(x)=\exp(i\la
    x,k\ra)$, where $\la x,k\ra$ is the scalar product in $X$.  Then
    we get a Banach $X$-module structure on $\Hr$ by setting
    $(V_ku)(x)=\exp(i\la x,k\ra)u(x)$, where $\la x,k\ra$ is the
    scalar product. Note that $V_k\Hr^s\subset\Hr^s$ and
    $\|V_k\|\leq C(1+|k|)^s$. It is easy to see that the Banach
    module structure associated to this $X$-module structure
    coincides with that defined in Example \ref{ex:sob}.
}\end{example}

\begin{remark}\label{r:wmod}{\rm
    Algebras $\Ac$ as in Lemma \ref{lm:ch} can be easily constructed
    in this context. Indeed, if $\omega$ is a sub-multiplicative
    function on $X^*$, i.e.\ a Borel map $X^*\rightarrow [1,\infty[$
    satisfying $\omega(k'k'')\leq\omega(k')\omega(k'')$ (hence
    $\omega$ is locally bounded), let $C^{(\omega)}(X)$ be the set
    of functions $\varphi$ whose Fourier transform
    $\widehat{\varphi}$ satisfies
\begin{equation}\label{e:fn}
\|\varphi\|_{C^{(\omega)}}:=
\int_{X^*}|\widehat{\varphi}(k)|\omega(k)dk< \infty.
\end{equation}
Then $C^{(\omega)}(X)$ is a subalgebra of $C_0(X)$ and is a Banach
algebra for the norm (\ref{e:fn}). Moreover, $C^{(a)}(X)\subset
C^{(\omega)}(X)$ densely and the net $\{e_\alpha\}$ defined in the
proof of Lemma \ref{l:group3} is an approximate unit of
$C^{(\omega)}(X)$.  If $\|V_k\|_{\Bc(\Hr)}\leq c\omega(k)$ for some
number $c>0$ then $\varphi(Q)$ is well defined for each $\varphi\in
C^{(\omega)}(X)$ by the relation (\ref{e:group1}) and
$\Phi(\varphi)=\varphi(Q)$ is a continuous morphism with dense range
of $C^{(\omega)}(X)$ into $\Mc(\Hr)$. We could take
$\omega(k)=\sup(1,\|V_k\|_{\Bc(\Hr)})$ but if a second Banach
$X$-module $\Kr$ is given then it is more convenient to take
$\omega(k)=\sup\{1, \|V_k\|_{\Bc(\Hr)}, \|V_k\|_{\Bc(\Kr)}\}$.
}\end{remark}

The adjoint of a reflexive Banach $X$-module has a natural structure
of Banach $X$-module.  Indeed, a weakly continuous representation is
strongly continuous, so we can equip the adjoint space $\Hr^*$ with
the Banach $X$-module structure defined by the representation
$k\mapsto (V_{\bar k})^*$, where $\bar k=k^{-1}$ is the complex
conjugate of $k$.

\medskip

The group $X$ is, in particular, a locally compact topological
space, hence the notion of Banach $X$-module in the sense of Example
\ref{ex:xmod} makes sense. But this is in fact a particular case of
that of Banach $X$-module in the sense of Definition \ref{d:xm}.
Indeed, according to Remark \ref{re:xmod}, we get a strongly
continuous representation of $X^*$ on $\Hr$ by setting $V_k=k(Q)$.
In the case of Hilbert $X$-modules we have a more precise fact.
\begin{lemma}\label{l:bhm}
  Let $\Hr$ be a Hilbert space. Then giving a Hilbert $X$-module
  structure on $\Hr$ is equivalent with giving on $\Hr$ a Banach
  $X$-module structure over the group $X$ such that the
  representation $\{V_k\}_{k\in X^*}$ is unitary. The relation
  between the two structures is determined by the condition
  $V_k=k(Q)$.
\end{lemma}
\proof If $\Hr$ is a Hilbert $X$-module then we can define
$V_k=k(Q)\in\Bc(\Hr)$ and check that $\{V_k\}_{k\in X^*}$ is a
strongly continuous unitary representation of $X^*$ on $\Hr$ with
the help of Remark \ref{re:xmod2}.  Reciprocally, it is well known
that such a representation allows one to equip $\Hr$ with a Hilbert
$X$-module structure. The main point is that the estimate
$\|\varphi(Q)\|\leq\sup|\varphi|$ holds, see \cite{Loo}.  \qed

Banach $X$-modules over the group $X$ which are not Hilbert
$X$-modules often appear in the following context (cf.\ Example
\ref{ex:sob&} in the case $s>0$).

\begin{definition}\label{d:5.2}
  A \emph{stable Friedrichs $X$-module} over the group $X$ is a
  Friedrichs $X$-module $(\Gr,\Hr)$ satisfying $V_k\Gr\subset\Gr$
  for all $k\in X^*$ and such that if $u\in\Gr$ and if $K\subset
  X^*$ is compact then $\sup_{k\in K}\|V_ku\|_{\Gr}<\infty$.
\end{definition}
Here $V_k=k(Q)$. It is clear that $V_k\Gr\subset\Gr$ implies
$V_k\in\Bc(\Gr)$ and that the local boundedness condition implies
that the map $k\mapsto V_k\in\Bc(\Gr)$ is a weakly, hence strongly,
continuous representation of $X^*$ on $\Gr$ (not unitary in
general).  The local boundedness condition is automatically
satisfied if $X^*$ is second countable.

Thus, \emph{if $(\Gr,\Hr)$ is a stable Friedrichs $X$-module,
then $\Gr$ is equipped with a canonical Banach $X$-module
structure}. Then, by taking adjoints, we get
a natural Banach $X$-module structure on $\Gr^*$ too. Our
definitions are such that after the identifications
$\Gr\subset\Hr\subset\Gr^*$ the restriction to $\Hr$ of the
operator $V_k$ acting in $\Gr^*$ is just the initial $V_k$.
Indeed, we have $V_k^*=V_k^{-1}=V_{\bar k}$ in $\Hr$.
Thus there is no ambiguity in using the same notation $V_k$ for the
representation of $X^*$ in the spaces $\Gr,\Hr$ and $\Gr^*$.

\begin{proposition}\label{p:sfxm}
  If $\Kr$ is a Banach space then
  $\bol(\Kr,\Gr)\subset\bol(\Kr,\Hr)$, and if $\Kr$ is a Banach
  module then $\bql(\Kr,\Gr)\subset\bql(\Kr,\Hr)$.
\end{proposition}    
\proof If $S\in\bol(\Kr,\Gr)$ then $S=\vphi(Q)T$ for some
$\vphi\in C^{(\omega)}(X)$ with
$\omega(k)=\sup(1,\|V_k\|_{\Bc(\Gr)})$ and some $T\in\Bc(\Kr,\Gr)$
(see Remark \ref{r:wmod}).  But clearly such a $\vphi(Q)$ belongs to
the multiplier algebra of $\Hr$ and $T\in\Bc(\Kr,\Hr)$.  \qed

\subsection{Regular operators are decay preserving}

We show now that, in the case of Banach $X$-modules over groups, the
decay preserving property is related to regularity in the sense of
the next definition.
\begin{definition}
  Let $\Hr$ and $\Kr$ be Banach $X$-modules.  We say that a
  continuous operator $S:\Hr\rarrow\Kr$ is of class $\cu(Q)$, and we
  write $S\in \cu(Q;\Hr,\Kr)$, if the map $k\mapsto V_k^{-1} S
  V_k\in\Bc(\Hr, \Kr)$ is norm continuous.
\end{definition}
Note that norm continuity at the origin implies norm continuity
everywhere.  The class of regular operators is stable under
algebraic operations:
\begin{proposition}\label{p:apr}
Let $\Gr,\Hr,\Kr$ be Banach $X$-modules.\\
(i) If $S\in \cu(Q;\Hr,\Kr)$ and $T\in \cu(Q;\Gr,\Hr)$ then
$ST\in \cu(Q;\Gr,\Kr)$.\\
(ii) If $S\in \cu(Q;\Hr,\Kr)$ is bijective,
then $S^{-1}\in \cu(Q;\Kr,\Hr)$.\\
(iii) If $S\in \cu(Q;\Hr,\Kr)$ and $\Hr,\Gr$ are reflexive, then
$S^*\in \cu(Q;\Kr^*,\Hr^*)$.
\end{proposition}
\proof We prove only (ii). If we set $S_k=V_k^{-1} S V_k$ then
$V_k^{-1} S^{-1} V_k=S_k^{-1}$, hence
$$
\|V_k^{-1} S^{-1} V_k-S^{-1}\|=
\|S_k^{-1}-S^{-1}\|=\|S_k^{-1}(S-S_k)S^{-1}\|
\leq C\|S-S_k\|
$$
if $k$ is in a compact set, and this tends to zero as $k\rarrow0$.
\qed

\begin{proposition}\label{prop28}
If $T\in \cu(Q;\Hr,\Kr)$ then $T$ is  decay preserving.
\end{proposition}
\proof
We show that $\varphi(Q)T\in \bor(\Hr, \Kr)$ if
$\varphi\in C^{(a)}(X)$. A similar argument gives
$T\varphi(Q)\in \bol(\Hr, \Kr)$. Set $T_k=V_kT V_k^{-1}$, then
\begin{equation*}
\varphi(Q)T=\int_{X^*}\widehat{\varphi}(k)V_k T dk=\int_{X^*} T_k
\widehat{\varphi}(k)V_k dk.
\end{equation*}
Since $k\mapsto T_k$ is norm continuous on the compact support of
$\widehat{\varphi}$, for each $\varepsilon>0$ we can construct, with
the help of a partition of unity, functions $\theta_i\in \cc(X^*)$
and operators $S_i\in \Bc(\Hr, \Kr)$, such that $\|T_k-\sum_{i=1}^n
\theta_i(k)S_i\|<\varepsilon$ if $\widehat\varphi(k)\neq0$. Thus
\begin{equation*}
\|\varphi(Q)T-\sum_{i=1}^n\int_{X^*}
\theta_i(k)S_i\widehat{\varphi}(k)
V_k dk\|\leq \varepsilon\sum_{i=1}^n\int_{X^*}|\widehat{\varphi}(k)|
\|V_k\|_{\Bc(\Hr)}dk.
\end{equation*}
Now, since $\bor(\Hr, \Kr)$ is a norm closed subspace, it suffices
to show that the operator
$\int_{X^*}\theta_i(k)S_i\widehat{\varphi}(k)V_k dk$ belongs to
$\Bc_0^{\,r}(\Hr, \Kr)$ for each $i$.  But if $\psi_i$ is the
inverse Fourier transform of $\theta_i\widehat{\varphi}$ then this
is $S_i\psi_i(Q)$ and $\psi_i\in C^{(a)}(X)$.\qed

Let $\Hr,\Kr$ be Hilbert $X$-modules over the group $X$.  We say
that an operator $S\in\Bc(\Hr,\Kr)$ is of \emph{finite
  range}\footnote{ If $X$ is an euclidean space and
  $\Hr=\Kr=L^2(X)$, the next condition means that there is
  $r<\infty$ such that the distribution kernel of $S$ satisfies
  $S(x,y)=0$ for $|x-y|>r$.}  if there is a compact neighborhood
$\Lambda$ of the origin such that for any compact sets $H,K\subset
X$ with $(H-K)\cap\Lambda=\emptyset$ we have
$\cchi_H(Q)S\cchi_K(Q)=0$. From Remark \ref{re:xmod2} we get that
this is equivalent to $S\cchi_K(Q)=\cchi_{K+\Lambda}(Q)S\cchi_K(Q)$
for any Borel set $K$.  A finite range operator is clearly decay
preserving. Moreover, the set of finite range operators is stable
under sums and products, and the adjoint of such an operator is also
of finite range.
\begin{proposition}\label{p:frange} 
  If $\Hr,\Kr$ are Hilbert $X$-modules over the group $X$, then each
  operator of class $\cu(Q)$ is a norm limit of a sequence
  of finite range operators.
\end{proposition}
\proof We fix a Haar measure $dk$ on $X^*$ and if $S\in\Bc(\Hr,
\Kr)$ and $\theta\in L^1(X^*)$ we define
\begin{equation}\label{eq21}
S_\theta=\int_{X^*}V_k^* S V_k \theta(k) dk.
\end{equation}
The integral is well defined because $k\mapsto V_k^*SV_k\in\Bc(\Hr)$
is a bounded strongly continuous map. In order to explain the main
idea of the proof we shall make a formal computation involving the
spectral measure $E(A)=\cchi_A(Q)$, see Remark \ref{re:xmod2} and
Lemma \ref{l:bhm} (we shall use the same notation for the spectral
measures in $\Hr$ and $\Kr$).  We have for $k\in X^*$ and
$\varphi(Q)\in B(X)$
\begin{eqnarray*}
\varphi(Q)V_k^*=\varphi(Q)k(Q)^*=(\varphi\overline{k})(Q)
=\int \varphi(x)\overline{k}(x)E(dx).
\end{eqnarray*}
Note also that for $x,y\in X$ we have $\overline{k}(x)k(y)
=k(-x)k(y)=k(y-x)$. Let $\widehat{\theta}(x) =\int
\overline{k(x)}\theta(k)dk$ be the Fourier transform of $\theta$.
Then if $\varphi, \psi\in B(X)$:
\begin{eqnarray}
\varphi(Q)S_\theta\psi(Q)&=&\int_{X^*}\theta(k)dk\int_X\int_X
\varphi(x)\overline{k}(x) k(y)\psi(y) E(dx) S E(dy)\nonumber\\
&=&\label{eq22}\int_X\int_X \widehat{\theta}(x-y)
\varphi(x)\psi(y) E(dx) S E(dy).
\end{eqnarray}
This clearly implies the following: 
$$
\leqno{\mbox{\bf($*$)}}\hspace{1mm}
\left\{
\begin{array}{ll}
\mbox{\rm If the support of }\widehat{\theta}
\mbox{\rm\, is a compact set } \Lambda
\mbox{\rm\  and if } \supp\varphi\cap(\Lambda+\supp\psi)=\varnothing
\\
\mbox{\rm then } \varphi(Q) S_\theta\psi(Q)=0.
\end{array}
\right.
$$
We shall note give a rigorous justification of (\ref{eq22}) but we
shall prove the preceding assertion, which suffices for our
purposes.  Observe that if ($*$) holds for a certain set of
operators $S$ then it also holds for the strongly closed linear
subspace of $\Bc(\Hr, \Kr)$ generated by it. So it suffices to prove
($*$) for $S$ an operator of rank one $S f=v\langle u,f\rangle$ with
some fixed $u\in\Hr$ and $v\in\Kr$. Now the computation giving
(\ref{eq22}) obviously makes sense in the weak topology and gives
for $f\in\Hr$ and $g\in\Kr$:
\begin{eqnarray*}
\langle g, \varphi(Q)S_\theta\psi(Q) f\rangle
=\int_X\int_X\widehat{\theta}(x-y)\varphi(x)\psi(y)
\langle g, E(dx) u\rangle \langle u, E(dy) f\rangle,
\end{eqnarray*}
hence ($*$) holds for such $S$.

Finally, note that if $S\in \cu(Q)$ then $S$ is norm limit of
operators of the form $S_\theta$. For this it suffices to take
$\theta=|K|^{-1}\cchi_{K}$ where $K$ runs over the set of open
relatively compact neighbourhoods of the neutral element of $X^*$,
$|K|$ being the Haar measure of $K$.  Then, by approximating
conveniently $\theta$ in $L^1$ norm, one shows that $S$ is norm
limit of operators $S_\theta$ such that $\widehat{\theta}$ has
compact support.  \qed
\begin{proposition}\label{pr:frange*}
  Assume that $X$ is a disjoint union $X=\cup_{a\in A} X_a$ of Borel
  sets $X_a$ such that: 1) there is a compact set $K$ such that each
  $X_a$ is a translate of a subset of $K$, and 2) for each compact
  neighborhood $\Lambda$ of the origin, the number of sets
  $X_b+\Lambda$ which intersects a given $X_a$ is bounded by a
  constant independent of $a$. Then, if $\Hr,\Kr$ are Hilbert
  $X$-modules over the group $X$, a finite range operator is of
  class $\cu(Q)$.
\end{proposition}
\proof Let $S$ be a finite range operator and let $\Lambda$ be such
that $\cchi_H(Q)S\cchi_K(Q)=0$ if $H,K$ are compact sets with
$(H-K)\cap\Lambda=\emptyset$.  Let $\cchi_a$ be the characteristic
function of $X_a$ and $\varphi_a$ that of $Y_a=X_a+\Lambda$. We can
assume that $A\subset X$ and that $X_a=a+K_a $ for some
$K_a\subset\Lambda$.  We shall abbreviate $\cchi_a=\cchi_a(Q)$ and
$\varphi_a=\varphi_a(Q)$.  We have $\sum_a\cchi_a=1$ strongly on
$\Hr$, cf.\ Remark \ref{re:xmod2}, and
$[V_k,S]\cchi_a=\varphi_a[V_k,S]\cchi_a$ because $V_k=k(Q)$. Thus
there is a constant $C$, depending only on an upper bound for the
number of $Y_b$ which intersects a fixed $X_a$, such that for $u\in
\Hr$ with compact support:
\begin{eqnarray*}
\|[V_k,S]u\|^2 & \leq & C\sum\|\varphi_a[V_k,S]\cchi_au\|^2\\
& = &
C\sum\|\varphi_a[V_k-k(a),S]\cchi_a\cdot\cchi_au\|^2\\
& \leq & 2CL_k\sum \|\cchi_au\|^2=2CL_k\|u\|^2
\end{eqnarray*} 
where $L_k=\sup_a\|(V_k-k(a))\vphi_a\|$. But
$$
\|(V_k-k(a))\vphi_a\|\leq \sup_{y\in Y_a}|k(y)-k(a)|=
\sup_{y\in Y_a}|k(y-a)-1|\leq \sup_{x\in L}|k(x)-1|
$$
where $L=\Lambda+\Lambda$ is a compact set. Thus $L_k\rarrow0$ if
$k\rarrow0$ in $X^*$.
\qed

If $X$ is an abelian locally compact group then there is enough
structure in order to develop a rich pseudo-differential calculus in
$L^2(X)$ and Proposition \ref{prop28} shows that many
pseudo-differential operators are decay preserving. We give a simple
example below. If $\varphi$ and $\psi$ are Borel functions on $X$
and $X^*$ respectively then, following standard quantum mechanical
conventions, we denote by $\varphi(Q)$ the operator of
multiplication by $\varphi$ in $L^2(X)$ and we set \label{w(p)}
$\psi(P)=\mathcal F^{-1}M_\psi\mathcal F$, where $M_\psi$ is the
operator of multiplication by $\psi$ in $L^2(X^*)$.

Let $\cbu(X)$ and $\cbu(X^*)$ be the algebras of bounded uniformly
continuous functions on $X$ and $X^*$ respectively. Below the space
$L^2(X)$ is equipped with its natural Hilbert $X$-module structure.
\begin{proposition}\label{prop29}
  The $C^*$-algebra generated by the operators $\varphi(Q)$ and
  $\psi(P)$, with $\varphi\in \cbu(X)$ and $\psi\in \cbu(X^*)$,
  consists of decay preserving operators.
\end{proposition}
\proof By Proposition \ref{p:aql}, $\Bc_q(L^2(X))$ is a
$C^*$-algebra, hence it suffices to show that each $\varphi(Q)$ and
$\psi(P)$ is decay preserving. For $\varphi(Q)$ the assertion is
trivial while for $\psi(P)$ we apply Proposition \ref{prop28}.\qed

\subsection{Compact perturbations in modules over abelian groups}
\label{cpmag}
In the present context it is possible to improve the results of
Section \ref{mr}.
\begin{lemma}\label{l:ab}
  Let $(\Gr,\Hr)$ and $(\Er,\Kr)$ be stable Friedrichs $X$-modules
  over the group $X$.  Let $D\in\Bc(\Gr,\Er)$ and
  $a\in\Bc(\Er,\Er^*)$ be operators of class $\cu(Q)$ such that
  $D^*aD-z:\Gr\rightarrow\Gr^*$ is bijective for some complex number
  $z$ and let $\Delta_a$ the operator on $\Hr$ associated to
  $D^*aD$.  Then the operator $D(\Delta_a-z)^{-1}\in\Bc(\Hr,\Er)$ is
  decay preserving.
\end{lemma}
\proof The lemma is an easy consequence of Propositions \ref{p:apr}
and \ref{prop28}.  Indeed, due to Proposition \ref{prop28}, it
suffices to show that the operator $D(\Delta_a-z)^{-1}$ is of class
$\cu(Q;\Hr,\Er)$.  We shall prove more, namely that
$D(D^*aD-z)^{-1}$ is of class $\cu(Q;\Gr^*,\Er)$.  Since $D$ is of
class $\cu(Q;\Gr,\Er)$, and due to (i) of Proposition \ref{p:apr},
it suffices to show that $(D^*aD-z)^{-1}$ is of class
$\cu(Q;\Gr^*,\Gr)$. But $D^*aD-z$ is of class $\cu(Q;\Gr,\Gr^*)$ by
(i) and (iii) of Proposition \ref{p:apr} and is a bijective map
$\Gr\rarrow\Gr^*$, so the result follows from (ii) of Proposition
\ref{p:apr}.  \qed

\begin{theorem}\label{t:main-3}
  Let $X$ be an abelian locally compact group and let $(\Gr,\Hr)$ be
  a compact stable Friedrichs $X$-module and $(\Er,\Kr)$ a stable
  Friedrichs $X$-module.  Assume that $D\in\Bc(\Gr,\Er)$ and
  $a,b\in\Bc(\Er,\Er^*)$ are operators of class $\cu(Q)$ such that
  the operators $D^*aD-z$ and $D^*bD-z$ are bijective maps
  $\Gr\rightarrow\Gr^*$ for some complex number $z$.  If
  $a-b\in\bol(\Er,\Er^*)$ then $\Delta_b$ is a compact perturbation
  of $\Delta_a$.
\end{theorem}
\proof The proof is a repetition of that of Theorem \ref{t:main-2}.
The only difference is that we write directly
$$
R=(D(\Delta_a^*-\bar z)^{-1})^*(b-a)D(\Delta_b-z)^{-1}
$$
and observe that $(b-a)D(\Delta_b-z)^{-1}\in\bol(\Hr,\Er^*)$ and
that $(D(\Delta_a^*-\bar z)^{-1})^*$ as an operator
$\Er^*\rarrow\Hr$ is decay preserving by (2) of Proposition
\ref{p:aql} and because the operator $D(\Delta_a^*-\bar
z)^{-1}:\Hr\rarrow\Er$ is decay preserving by Lemma \ref{l:ab}.
\qed

We finish with a simple corollary of Theorem \ref{t:main} which
nevertheless covers interesting examples of differential operators
of any order.
\begin{theorem}\label{t:main-ab}
  Assume that $(\Gr,\Hr)$ is a compact stable Friedrichs $X$-module
  over the group $X$ and that condition (AB) from page \pageref{AB}
  is satisfied. Let us also assume that $\wtilde
  A-z:\Gr\rarrow\Gr^*$ is bijective for some
  $z\in\rho(A)\cap\rho(B)$ and that $\wtilde A\in\cu(Q;\Gr,\Gr^*)$.
  If $\wtilde B-\wtilde A\in\bol(\Gr,\Gr^*)$, then $B$ is a compact
  perturbation of $A$.
\end{theorem}
\proof We apply Theorem \ref{t:main} with $\Kr=\Gr^*$, $S$ the
identity operator and $T=\wtilde B-\wtilde A$. Then $(\wtilde
A-z)^{-1}$ is of class $\cu(Q;\Gr^*,\Gr)$ by (ii) of Proposition
\ref{p:apr}, hence $(\wtilde A-z)^{-1}\in \bq(\Gr^*,\Gr)$ by
Proposition \ref{prop28}. But this is stronger than $(\wtilde
A-z)^{-1}\in \bql(\Gr^*,\Hr)$, as follows from Proposition
\ref{p:sfxm}.  \qed

\subsection{A class of hypoelliptic operators on abelian groups}
\label{anyord}
In this subsection we assume that $X$ is non-discrete, so $X^*$ is
non-compact. We also fix a finite dimensional complex Hilbert space
$E$ and take $\Hr=L^2(X;E)$ equipped with its natural Hilbert
$X$-module structure. Note that, according to our conventions, the
unitary representation of $X^*$ is given by the multiplication
operators $V_k=k(Q)$.

Let $w:X^*\rarrow[1,\infty[$ be a continuous function
satisfying $w(k)\rarrow\infty$ as $k\rarrow\infty$ and
such that $w(k'k)\leq \omega(k')w(k)$ holds for some function
$\omega$ and all $k',k$. We shall assume that $\omega$ is the
smallest function satisfying the preceding estimate. It is
clear then that $\omega$ is sub-multiplicative in the sense
defined in Remark \ref{r:wmod} (see \cite[Section 10.1]{Ho} for
this construction).

Then $w(P)$ is a self-adjoint operator on $\Hr$ with $w(P)\geq1$
(see page \pageref{w(p)} for this notation). We denote
$\Hr^w=\Dc(w(P))$ and equip it with the Banach $X$-module structure
given by the norm $\|u\|_w=\|w(P)u\|$ and the representation
$V_k|\Hr^w$.  Obviously, this space is a generalization of the usual
notion of Sobolev spaces.

\begin{lemma}\label{l:omega}
$(\Hr^w,\Hr)$ is a compact stable Friedrichs $X$-module.
\end{lemma}
\proof If $\varphi\in\co(X)$ then $\varphi(Q)w(P)^{-1}$ is a compact
operator because $w^{-1}$ belongs to $C_0(X)$, hence
$\varphi(Q)\in\Kc(\Hr^w,\Hr)$. Then observe that
$V_k^{-1}w(P)V_k=w(kP)$ and $w(kP)\leq \omega(k)w(P)$.  Thus $V_k$
leaves stable $\Hr^w$ and we have the estimate
$\|V_k\|_{\Bc(\Hr^w)}\leq \omega(k)$.  \qed

We call \emph{uniformly hypoelliptic} an operator $A$ on $\Hr$ such
that there are $w$ as above and an operator $\wtilde
A\in\Bc(\Hr^w,\Hr^{w*})$ such that $\wtilde
A-z:\Hr^w\rarrow\Hr^{w*}$ is bijective for some complex $z$ and such
that $A$ is the operator induced by $\wtilde A$ in $\Hr$ (see the
Appendix).  For example, the constant coefficients case with $E=\C$
corresponds to the choice $A=h(P)$ with $h:X^*\rarrow\C$ a Borel
function such that $c'w^2\leq1+|h|\leq c''w^2$ and such that the
range of $h$ is not dense in $\C$. We shall justify our terminology
in the remark at the end of this subsection.

Theorem \ref{t:main-ab} is quite well adapted to show the stability
of the essential spectrum of such operators under perturbations
which are small at infinity.  We stress that the differential
operators covered by these results can be of any order and that in
the usual case when the coefficients are complex measurable
functions a condition of the type $\wtilde
A\in\cu(Q;\Hr^w,\Hr^{w*})$ is very general, if not automatically
satisfied (see the remark at the end of this subsection).  Hence the
only condition really relevant in this context is $\wtilde B-\wtilde
A\in\bol(\Hr^w,\Hr^{w*})$ and the main point is that it allows
perturbations of the higher order coefficients even in the
non-smooth case.

It is clear that these results can be used to establish the
stability of the essential spectrum of pseudo-differential operators
on finite dimensional vector spaces over local fields
(see \cite{Sa,Ta}) under perturbations of the same order.

\medskip

We shall give an application of physical interest to Dirac
operators. Let $X=\R^n$ and let
$\alpha_0\equiv\beta,\alpha_1,\dots,\alpha_n$ be symmetric operators
on $E$ such that $\alpha_j\alpha_k+\alpha_k\alpha_j=\delta_{jk}$.
Then the free Dirac operator is $D=\sum_{k=1}^n\alpha_kP_k+m\beta$
for some real number $m$. The natural compact stable Friedrichs
$X$-module in this context is $(\Hr^{1/2},\Hr)$. Note that we use
the same notation $\Hr^s$ for Sobolev spaces of $E$-valued
functions.

\begin{proposition}\label{p:dirac}
  Let $V,W$ be measurable functions on $X$ with values symmetric
  operators on $E$ and such that the operators of multiplication by
  $V$ and $W$ define continuous maps $\Hr^{1/2}\rarrow\Hr^{-1/2}$
  and $V-W\in\bo(\Hr^{1/2},\Hr^{-1/2})$.  Assume that $D+V+i$ and
  $D+W+i$ are bijective maps $\Hr^{1/2}\rarrow\Hr^{-1/2}$.  Then
  $D+V$ and $D+W$ induce self-adjoint operators $A$ and $B$ in
  $\Hr$, $B$ is a compact perturbation of $A$, and $\se(B)=\se(A)$.
\end{proposition}

This follows immediately from Theorem \ref{t:main-ab}.  We stress
that the main new feature of this result is that the ``unperturbed''
operator $A$ is locally as singular as the ``perturbed'' one $B$.
The assumptions imposed on $V,W$ are quite general, compare with
\cite{Ar,AY,Kl,N1,N2}.

\medskip

\noindent{\bf Remark:}
In order to clarify the relation between the notion of uniform
hypoellipticity introduced above and the original notion of
hypoellipticity due to H\"ormander, we shall consider the case of
differential operators on $\R^n$ (which is identified with its dual
group in the standard way).  Assume first that $h$ is a polynomial
on $\R^n$ and that $A=h(P)$. Then the function defined by
$w(k)^4=\sum_\alpha|h^{(\alpha)}(k)|^2$ satisfies
$w(k'+k)\leq(1+c|k'|)^{m/2}w(k)$, where $c$ is a number and $m$ is
the order of $h$, see \cite[Example 10.1.3]{Ho}.  Now the ``form
domain'' of the operator $h(P)$ in $L^2(\R^n)$ is the space
$\Gr=\Dc(|h(P)|^{1/2})$ and this domain is stable under $V_k=\exp
i\la k,Q\ra$ if and only the function $w$ satisfies $w^2\leq
c(1+|h|)$, see Lemma 7.6.7 in \cite{ABG}.  On the other hand,
Definition 11.1.2 and Theorem 11.1.3 from \cite{Ho} show that $A$ is
hypoelliptic if and only if $h^{(\alpha)}(k)/h(k)\rarrow0$ when
$k\rarrow\infty$, for all $\alpha\neq0$. So in this case we have
$c'w^2\leq1+|h|\leq c''w^2$ and the operator $h(P)$ is uniformly
hypoelliptic in our sense if $h(\R^n)$ is not dense in $\C$. If
$n=2$ then $h(k)=k_1^4+k_2^2$ is a simple example of polynomial
which satisfies all these conditions but is not elliptic.  See
\cite[Subsections 2.7-2.10]{GM} for the case of matrix valued
functions $h$.

In the variable coefficient case the notion of hypoellipticity
defined in \cite[Definition 13.4.3]{Ho} is a local one and one may
consider different global versions. For instance,
\cite[Theorem 13.4.4]{Ho} suggests that the
notion we introduced above is natural for operators of uniform
constant strength. But the uniform constant strength condition
is not satisfied by the operators with polynomial coefficients,
for example, hence such operators are not uniformly hypoelliptic
in our sense in general.

\section{Operators in divergence form on Euclidean
spaces}\label{div}
The results of this section are corollaries of Theorem
\ref{t:main-3}. We shall take $X=\R^n$, we fix a finite dimensional
Hilbert space $E$, and choose $\Hr=L^2(X;E)$ with the obvious Hilbert
$X$-module structure.  If $s\in\R$ then $\Hr^s$ is the usual
Sobolev space of $E$ valued functions. Then for each $s>0$ the
couple $(\Hr^s, \Hr)$ is a compact stable Friedrichs $X$-module,
cf.\ Examples \ref{ex:sob}, \ref{ex:sob*} and \ref{ex:sob&}.

Let us describe the objects which appear in Theorem \ref{t:main-3}
in the present context.
We fix  an integer $m\geq 1$ and take $\Gr=\Hr^m$.
Let $\Kr=\bigoplus_{|\alpha|\leq m} \Hr_\alpha$, where
$\Hr_{\alpha}\equiv\Hr$,
with the natural direct sum Hilbert $X$-module structure.
Here $\alpha$ are multi-indices $\alpha\in\N^n$ and
$|\alpha|=\alpha_1+\dots+\alpha_n$. Then we define
\begin{equation*}
\Er=\bigoplus_{|\alpha|\leq m} \Hr^{m-|\alpha|}=
\{(u_\alpha)_{|\alpha|\leq m}\in \Kr \mid
u_\alpha\in\Hr^{m-|\alpha|} \}
\end{equation*}
equipped with the Hilbert direct sum structure. It is obvious that
$(\Er,\Kr)$ is a stable Friedrichs $X$-module (but not compact).

We set $P_k=-i\partial_k$, where $\partial_k$ is the derivative with
respect to the $k$-th variable, and $P^\alpha=P_1^{\alpha_1}\ldots
P_n^{\alpha_n}$ if $\alpha\in\N^n$. Then for $u\in\Gr$ let
$Du=(P^\alpha u)_{|\alpha|\leq m}\in \Kr$. Since
\begin{equation*}
\|Du\|^2=\sum_{|\alpha|\leq m}\|P^\alpha u\|^2=\|u\|^2_{\Hr^m}
\end{equation*}
we see that $D:\Gr\rightarrow\Kr$ is a linear isometry.
Moreover, we have  defined $\Er$ such as to have
$D\Gr\subset\Er$, hence  $D\in \Bc(\Gr, \Er)$.
We have $D\in\cu(Q;\Gr,\Er)$ because 
$$
V_k^{-1}DV_k=(V_k^{-1}P^\alpha V_k)_{|\alpha|\leq m}=
((P+k)^\alpha)_{|\alpha|\leq m}
$$
and this a polynomial in $k$ with coefficients in $\Bc(\Gr,\Er)$.

We shall identify $\Hr^*=\Hr$ and $\Kr^*=\Kr$, which implies
$\Gr^*=\Hr^{-m}$ and
\begin{equation*}
\Er^*=\oplus_{|\alpha|\leq m}\Hr^{|\alpha|-m}.
\end{equation*}
The operator $D^*\in\Bc(\Er^*, \Gr^*)$ acts as follows:
\begin{equation*}
D^*(u_\alpha)_{|\alpha|\leq m}=\sum_{|\alpha|\leq m}P^\alpha u_\alpha
\in \Hr^{-m},
\end{equation*}
because $u_\alpha\in \Hr^{|\alpha|-m}$.

By taking into account the given expressions for $\Er$ and $\Er^*$
we see that we can identify an operator $a\in\Bc(\Er, \Er^*)$ with a
matrix of operators $a=(a_{\alpha\beta})_{|\alpha|,|\beta|\leq m}$,
where $a_{\alpha\beta}\in\Bc(\Hr^{m-|\beta|}, \Hr^{|\alpha|-m})$ and
\begin{equation*}
a(u_\beta)_{|\beta|\leq m} = 
\big(\sum_{|\beta|\leq m} a_{\alpha\beta}
u_\beta\big)_{|\alpha|\leq m}.
\end{equation*}
Then we clearly have
\begin{equation}\label{e:div}
D^* aD=\sum_{|\alpha|, |\beta|\leq m} 
P^\alpha a_{\alpha\beta} P^\beta.
\end{equation}
which is a general version of a differential operator in divergence
form.  We must, however, emphasize that our $a_{\alpha\beta}$ are
not necessarily ($B(E)$ valued) functions, they could be
pseudo-differential or more general operators.

In view of the statement of the next theorem, we note that, since
the Sobolev spaces are Banach $X$-modules over the group $X$, the
class of regularity $\cu(Q;\Hr^s,\Hr^t)$ is well defined for all
real $s,t$. A bounded operator $S:\Hr^s\rarrow\Hr^t$ belongs to this
class if and only if the map $k\mapsto
V_{-k}SV_k\in\Bc(\Hr^s,\Hr^t)$ is norm continuous.  In particular,
this condition is trivially satisfied if $S$ is the operator of
multiplication by a function, because then $V_k$ commutes with
$S$. Since the coefficients $a_{\alpha\beta}$ of the differential
expression (\ref{e:div}) are usually assumed to be functions, this
is a quite weak restriction in the setting of the next theorem.  The
condition $S\in\bol(\Hr^s,\Hr^t)$ is also well defined and it is
easily seen that it is equivalent to
\begin{equation}\label{e:bst}
\lim_{r\rarrow\infty}\|\theta(Q/r)S\|_{\Hr^s\rarrow\Hr^t}=0
\end{equation}
where $\theta$ is a $C^\infty$ function on $X$ equal to zero on a
neighborhood of the origin and equal to one on a neighborhood of
infinity.  Now we can state the following immediate consequence of
Theorem \ref{t:main-3}.
\begin{proposition}\label{p:applic1}
  Let $a_{\alpha\beta}$ and $b_{\alpha\beta}$ be operators of class
  $\cu(\Hr^{m-|\beta|}, \Hr^{|\alpha|-m})$ and such that the
  operators $D^*aD-z$ and $D^*bD-z$ are bijective maps
  $\Hr^m\rarrow\Hr^{-m}$ for some complex $z$. Let $\Delta_a$ and
  $\Delta_b$ be the operators in $\Hr$ associated to $D^*aD$ and
  $D^*bD$ respectively.  Assume that
\begin{equation}\label{e:sd8} 
\lim_{r\rarrow\infty}
\|\theta(Q/r)(a_{\alpha\beta}-b_{\alpha\beta})\|_{\Hr^{m-|\beta|}
\rarrow\Hr^{|\alpha|-m}}=0
\end{equation} 
for each $\alpha,\beta$, where $\theta$ is a function as above.
Then $\Delta_b$ is a compact perturbation of $\Delta_a$ and the
operators $\Delta_a$ and $\Delta_b$ have the same essential
spectrum.
\end{proposition}

\noindent{\bf Example:}
In the simplest case the coefficients $a_{\alpha\beta}$ and
$b_{\alpha\beta}$ of the principal parts (i.e.\ 
$|\alpha|=|\beta|=m$) are functions. Then the conditions become:
$a_{\alpha\beta}$ and $b_{\alpha\beta}$ belong to $L^\infty(X)$ and
$|a_{\alpha\beta}(x)-b_{\alpha\beta}(x)|\rarrow0$ as
$|x|\rarrow\infty$. Of course, the assumptions on the lowest order
coefficients are much more general.

\medskip

\noindent{\bf Example:}
We show here that ``highly oscillating potentials'' do not modify
the essential spectrum. If $m=1$ then the terms of order one of
$D^*aD$ are of the form $S=\sum_{k=1}^n(P_kv'_k+v_k''P_k)$, where
$v'_k\in\Bc(\Hr^1,\Hr)$ and $v''_k\in\Bc(\Hr,\Hr^{-1})$.  Choose
$v_k\in\Bc(\Hr^1,\Hr)$ symmetric in $\Hr$ and let
$v_k'=iv_k,v_k''=-iv_k$. Then $S=[iP,v]\equiv{\rm div\ \!} v$, with
natural notations, can also be thought as a term of order zero. Now
assume that $v_k$ are bounded Borel functions and consider a similar
term $T=[iP,w]$ for $D^*bD$. Then the condition
$|v_k(x)-w_k(x)|\rarrow0$ as $|x|\rarrow\infty$ suffices to ensure
the stability of the essential spectrum. However, the difference
$S-T$ could be a function which does not tend to zero at infinity in
a simple sense, being only ``highly oscillating''.  An explicit
example in the case $n=1$ is the following: a perturbation of the
form $\exp(x)(1+|x|)^{-1}\cos(\exp(x))$ is allowed because it is the
derivative of $(1+|x|)^{-1}\sin(\exp(x))$ plus a function which
tends to zero at infinity.

\medskip

In order to apply Proposition \ref{p:applic1} we need that $D^*
aD-z:\Hr^m\rarrow\Hr^{-m}$ be bijective for some $z\in\C$, and
similarly for $b$. A standard way of checking this is to require the
following \emph{coercivity condition}:
$$
\leqno{\mbox{\bf($C$)}}\hspace{1mm}
\left\{
\begin{array}{ll}
\mbox{\rm there are } \mu, \nu>0 \mbox{ such that for all }u\in\Hr^m:
\\
\sum_{|\alpha|,|\beta|\leq m} \Re \langle P^\alpha u, a_{\alpha\beta}
P^\beta u\rangle\geq \mu\|u\|^2_{\Hr^m}-\nu\|u\|^2_{\Hr}
\end{array}
\right.
$$

\noindent{\bf Example:}
One often imposes a stronger ellipticity condition that we describe
below.  Observe that the coefficients of the highest order part of
$D^*aD$ defined by $A_0=\sum_{|\alpha|=|\beta|=m}P^\alpha
a_{\alpha\beta}P^\beta$ are operators $a_{\alpha\beta}\in\Bc(\Hr)$.
Then ellipticity means:
 $$
\leqno{\mbox{\bf($Ell$)}}\hspace{1mm}
\left\{
\begin{array}{ll}
\mbox{\rm there is } \mu>0 \mbox{ such that if } u_\alpha\in\Hr
\mbox{ for }|\alpha|=m \mbox{ then }
\\
\sum_{|\alpha|=|\beta|= m} \Re \langle u_\alpha, a_{\alpha\beta}
u_\beta\rangle\geq \mu\sum_{|\alpha|=m}\|u_\alpha\|^2_\Hr.
\end{array}
\right.
$$
But we emphasize that, our conditions on the lower order terms
being quite general, e.g. the $a_{\alpha\beta}$ could be differential
operators, so the terms of formally lower order could be of order
$2m$ in fact, we have to supplement the ellipticity condition
{\bf($Ell$)} with a condition saying that the rest of the terms
$A_1=\sum_{|\alpha|+|\beta|<2m}P^\alpha a_{\alpha\beta}P^\beta$ is
small with respect to $A_0$. For example, we may require the
existence of some $\delta<\mu$ and $\gamma>0$ such that
\begin{equation}\label{***}
|\sum_{|\alpha|+|\beta|<2m} \Re \langle P^\alpha u, a_{\alpha\beta}
P^\beta u\rangle| \leq \delta\|u\|^2_{\Hr^m}+\gamma\|u\|^2_{\Hr}.
\end{equation}
This is satisfied if $A_1\Hr^m\subset \Hr^{-m+\theta}$ for some
$\theta>0$, because for each $\varepsilon>0$ there is
$c(\varepsilon)<\infty$ such that $\|u\|_{\Hr^{m-\theta}}\leq
\varepsilon\|u\|_{\Hr^m}+c(\varepsilon) \|u\|_{\Hr}$.

\begin{remark}\label{r:00-0}{\rm
    If we use Theorem \ref{t:main-2} in the context of this section
    then we get the same conditions on the coefficients
    $a_{\alpha\beta}-b_{\alpha\beta}$ of the principal part (i.e.
    such that $|\alpha|=|\beta|=m$) of the operator $a-b$ but those
    on the lower order coefficients are less general.  Indeed, if
    $s+t>0$ the space $\bool(\Hr^s,\Hr^{-t})$ defined as the closure
    of $\bol(\Hr)$ in $\Bc(\Hr^s,\Hr^{-t})$ does not contain
    operators of order $s+t$, while $\bol(\Hr^s,\Hr^{-t})$ contains
    such operators.
}\end{remark}

\section{Weak decay preserving operators}\label{wvpr}
The purpose of the next two sections is to reconsider the examples
treated in Section \ref{div} and to prove some stability results for
perturbations which decay in a generalized sense, as described in
Examples \ref{ex:xmodf}--\ref{ex:meas}. This will be done in the
next section, this one contains some preparatory material concerning
weak decay preserving operators.

We first consider the setting of Example \ref{ex:meas}: $(X,\mu)$ is
a positive measure such that $\mu(X)=\infty$, $\Fc_\mu$ is the
filter of sets of co-finite measure\footnote{ Note that if $X$ is a
  locally compact space and $\mu$ a Radon measure then $\Fc_\mu$ is
  finer than the Fr\'echet filter.  Moreover, if $X$ is an abelian
  locally compact non-compact group then $\Fc_\mu$ is strictly
  included in the filter $\Fc_\w$ which will be defined below.}, and
$B_\mu(X)$ is the algebra of bounded measurable $\Fc_\mu$-vanishing
functions.  We recall that any direct integral of Hilbert spaces
over $X$ has a canonical Hilbert module structure with $B_\mu(X)$ as
multiplier algebra. To avoid ambiguities, we shall speak of
$\Fc_\mu$-decay preserving operators when we refer to this algebra.
Let $\{\Hr(x)\}_{x\in X}$ and $\{\Kr(x)\}_{x\in X}$ be measurable
families of Hilbert spaces with dimensions $\leq N$ for some
finite $N$. We shall use the notations introduced before Corollary
\ref{co:bmsft}.

\begin{theorem}\label{th:meas} 
  Let $S\in\Bc(\Hr,\Kr)\cap\Bc(\Hr_p,\Kr_p)$ for some $p\neq2$. If
  $p<2$ then $S$ is left $\Fc_\mu$-decay preserving and if $p>2$
  then $S$ is right $\Fc_\mu$-decay preserving
\end{theorem}
\proof We shall consider only the case $p<2$, the assertion for
$p>2$ follows by observing that
$S^*\in\Bc(\Kr,\Hr)\cap\Bc(\Kr_{p'},\Hr_{p'})$ and then using
Proposition \ref{p:aql}. We prove that for each measurable set $N$
of finite measure the operator $T=S\cchi_N(Q)$ has the property: if
$\varepsilon>0$ then there is a Borel set $F\in\Fc_\mu$ such that
$\|\cchi_F(Q)T\|\leq\varepsilon$ (then Proposition \ref{p:dql}
implies that $S$ is left $\Fc_\mu$-decay preserving).  Since $N$ is
of finite measure, $\cchi_N(Q)$ is a bounded operator $\Hr\rarrow
\Hr_p$, hence $T\in\Bc(\Hr,\Kr_p)$. The rest of the proof is a
straightforward application of Corollary \ref{co:bmsft}. Let $a>0$
real and let $F$ be the set of points $x$ such that $|g(x)|\leq a$.
Since $g\in L^q$ with $q<\infty$, we have $F\in\Fc_\mu$ and
$$
\|\cchi_F(Q)T\|_{\Bc(\Hr,\Kr)}=\|\cchi_F(Q)g(Q)R\|_{\Bc(\Hr,\Kr)}
\leq a\|R\|_{\Hr,\Kr)}.
$$
Thus it suffices to choose $a$ such that
$a\|R\|_{\Bc(\Hr,\Kr)}=\varepsilon$.  \qed

Let $X$ be a locally compact non-compact topological space and let
$\Hr$ be a Hilbert $X$-module. Then, due to Remark \ref{re:xmod2},
the operator $\vphi(Q)\in\Bc(\Hr)$ is well defined for all $\vphi\in
B(X)$. If $\Fc$ is a filter finer than the Fr\'echet filter on $X$
then
\begin{equation}
B_\Fc(X):=\{\vphi\in B(X)\mid \lim_\Fc\vphi=0\}
\end{equation}
is a $C^*$-algebra and we can consider on $\Hr$ the Hilbert
module structure defined by the multiplier algebra
$\Mc_\Fc:=\{\vphi(Q)\mid \vphi\in B_\Fc(X)\}$. We are
interested in the corresponding classes of decay improving or decay
preserving operators. To be precise, we shall speak in this
context of (left or right) \emph{$\Fc$-vanishing at infinity}
or of (left or right) \emph{$\Fc$-decay preserving} operators.
Below and later on we use the notation $N^\crm=X\setminus N$.
\begin{lemma}\label{l:filter}
Let $\Hr,\Kr$ be Hilbert $X$-modules. Then an operator
$S\in\Bc(\Hr,\Kr)$
is left $\Fc$-decay preserving if and only if for each Borel set $N$
with $N^\crm\in\Fc$ and for each $\veps>0$ there is a Borel set
$F\in\Fc$ such that $\|\cchi_F(Q)S\cchi_N(Q)\|\leq\veps$.
\end{lemma}
\proof We note first that the family of operators $\cchi_N(Q)$,
where $N$ runs over the family of Borel sets with complement in
$\Fc$, is an approximate unit for $B_\Fc(X)$. Indeed, if $\veps>0$
and $\vphi\in B_\Fc(X)$ then the set $N=\{x\mid |\vphi(x)|>\veps\}$
is Borel, its complement is in $\Fc$, and
$\sup_x|\vphi(x)(1-\cchi_N(x))|\leq\veps$.  Thus, according to
Proposition \ref{p:dql}, $S$ is left $\Fc$-decay preserving if and
only if $S\cchi_N(Q)$ is left $\Fc$-vanishing at infinity for each
$N$. Now the result follows from (\ref{e:lv}).  \qed

The main restriction we have to impose on $\Fc$ comes from the fact
that the Friedrichs couple $(\Gr,\Hr)$ which is involved in our
abstract compactness criteria must be such that
$\vphi(Q)\in\Kc(\Gr,\Hr)$ if $\vphi\in\ B_\Fc(X)$. Sometimes this
can be stated quite explicitly:
\begin{lemma}\label{l:best}
Let $X$ be an Euclidean space, $\Hr=L^2(X)$, and let $\Gr=\Hr^s$
be a Sobolev space of order $s>0$. If $\vphi\in\ B(X)$ then
$\vphi(Q)\in\Kc(\Gr,\Hr)$ if and only if
\begin{equation}\label{e:wv8}
\lim_{a\rarrow\infty}\int_{|x-a|\leq1}|\vphi(x)|dx=0.
\end{equation}
\end{lemma}
The importance of such a condition in questions of stability of the
essential spectrum has been noticed in \cite{He,LV,OS,We}. That it
is a natural condition follows also from the characterizations that
we shall give below in a more general context.

\medskip

Let $X$ be a locally compact non-compact abelian group. We shall
say that a function $\vphi\in B(X)$ is \emph{weakly vanishing
(at infinity)} if
\begin{equation}\label{e:wv8g}
\lim_{a\rarrow\infty}\int_{a+K}|\vphi(x)|dx=0
\mbox{ for each compact set }K.
\end{equation}
We shall denote by $B_\w(X)$ the set of functions $\vphi$ satisfying
(\ref{e:wv8g}). This is clearly a $C^*$-algebra. Note that it
suffices that the convergence condition in (\ref{e:wv8g}) be
satisfied for only one compact set $K$ with non-empty interior.

Let us now express the condition (\ref{e:wv8g}) in terms of
convergence to zero along a filter. We denote $|K|$ the exterior
(Haar) measure of a set $K\subset X$ and we set $K_a=a+K$ if $a\in
X$. A subset $N$ is called \emph{$\w$-small (at infinity)} if there
is a compact neighborhood $K$ of the origin such that
$\lim_{a\rarrow\infty}|N\cap K_a|=0$. The complement of a $\w$-small
set will be called \emph{$\w$-large (at infinity)}. The family
$\Fc_\w$ of all $\w$-large sets is clearly a filter on $X$ finer
than the Fr\'echet filter.

We give now a characterization of weakly vanishing functions in
terms of compactness properties. This characterization implies that
of Lemma \ref{l:best} if $X=\R^n$.  Observe that a Borel set is
$\w$-small if and only if its characteristic function weakly
vanishes at infinity. Denote $f*g$ the convolution of two functions
on $X$.
\begin{lemma}\label{lm:wvf}
For a function $\vphi\in B(X)$ the following conditions are
equivalent: {\rm(1)} $\vphi$ is weakly vanishing;
{\rm(2)} $\theta*|\vphi|\in\co(X)$ if $\theta\in\cc(X)$;
{\rm(3)} $\lim_{\Fc_\w}\vphi=0$;
{\rm(4)} $\vphi(Q)\psi(P)$ is a compact operator on $L^2(X)$ for all
$\psi\in\co(X)$. 
\end{lemma}
\proof The equivalence of (1) and (2) is clear because
$\int_{K_a}|\vphi|dx=(\cchi_K*|\vphi|)(a)$.  Then (3) means that for
each $\veps>0$ the Borel set $N$ where $|\vphi(x)|>\veps$ is
$\w$-small. Since $\cchi_N\leq\vphi/\veps$, the implication (2)
$\Rightarrow$ (3) is clear, while the reciprocal implication follows
from $\cchi_K*|\vphi|\leq\sup|\vphi|\cchi_K*\cchi_N+\veps|K|$.  If
(4) holds, let us choose $\psi$ such that its Fourier transform
$\what\psi$ be a positive function in $\cc(X)$ and let $f\in\cc(X)$
be positive and not zero. Since $\psi(P)f$ is essentially the
convolution of $\what\psi$ with $f$, there is a compact set $K$ with
non-empty interior such that $\psi(P)f\geq c\cchi_K$ with a number
$c>0$. Let $U_a$ be the unitary operator of translation by $a$ in
$L^2(X)$, then $U_af\rarrow0$ weakly when $a\rarrow\infty$, hence
$\|\vphi(Q)U_a\psi(P)f\|=\|\vphi(Q)\psi(P)U_af\|\rarrow0$.  Since
$U_a^*\vphi(Q)U_a=\vphi(Q-a)$ we get
$\|\vphi(Q-a)\cchi_K\|\rarrow0$, hence (1) holds.

Finally, let us prove that (1) $\Rightarrow$ (4). It suffices to
prove that $\vphi(Q)\psi(P)$ is compact if $\what\psi\in\cc(X)$
and for this it suffices that $\bar\psi(P)|\vphi|^2(Q)\psi(P)$
be compact. Since $\xi:=|\vphi|^2\in B_\w(X)$ and since $\psi(P)$
is the operator of convolution by a function $\theta\in\cc(X)$,
we are reduced to proving that the integral operator $S$ with
kernel $S(x,y)=\int\bar\theta(z-x)\xi(z)\theta(z-y)dz$ is compact.
If $K=\supp\theta$ and $\Lambda$ is the compact set $K-K$, then
clearly there is a number $C$ such that
$$|S(x,y)|\leq C\int_{K_x}\xi(z)dz\cchi_\Lambda(x-y)
\equiv\phi(x)\cchi_\Lambda(x-y)$$
where $\phi\in\co(X)$.  The last
term here is a kernel which defines a compact operator $T$.  Thus
$\eta(Q)S$ is a Hilbert-Schmidt operator for each $\eta\in\cc(X)$
and from the preceding estimate we get
$\|(S-\eta(Q)S)u\|\leq\|(1-\eta(Q))T|u|\|$ for each $u\in L^2(X)$.
Thus $\|S-\eta(Q)S\|\leq\|(1-\eta(Q))T\|$ and the right hand side
tends to zero if $\eta\equiv\eta_\alpha$ is an approximate unit for
$\co(X)$.  \qed

We shall consider now a general class of filters defined in terms of
the metric and measure space structure.  We consider only the case
of an Euclidean space $X$, the extension to the case of locally
compact groups or metric spaces being obvious. We set $B_a(r)=\{x\in
X\mid |x-a|<r\}$, $B_a=B_a(1)$ and $B(r)=B_0(r)$.  To each function
$\nu:X\rarrow]0,\infty[$ such that
$\liminf_{a\rarrow\infty}\nu(a)=0$ we associate a set of subsets of
$X$ as follows:
\begin{equation}\label{e:nu} 
\Nr_\nu=\{N\subset X\mid \limsup_{a\rarrow\infty}
\nu(a)^{-1}|N\cap B_a|<\infty\}.
\end{equation}
Clearly $\Fr_\nu=\{F\subset X\mid F^\crm\in\Nr_\nu\}$ is a filter on
$X$ finer than the Fr\'echet filter.

\begin{theorem}\label{t:fquasi} 
  Let $X=\R^n$ and let $\nu:X\rarrow]0,\infty[$ such that
  $\liminf_{a\rarrow\infty}\nu(a)=0$ and $\sup_{|b-a|\leq
    r}\nu(b)/\nu(a)<\infty$ for each real $r$.  If $S\in\Bc(L^2(X))$
  is of class $\cu(Q)$ and if $S\in\Bc(L^p(X))$ for some $p<2$, then
  $S$ is left $\Fr_\nu$-decay preserving.
\end{theorem} 
\proof We can approximate in norm in $\Bc(L^2(X))$ the operator $S$
by operators which are in $\Bc(L^2(X))\cap\Bc(L^p(X))$ and have
finite range. Indeed, the approximation procedure (\ref{eq21}) used
in the proof of Proposition \ref{p:frange} is such that it leaves
$\Bc(L^2(X))\cap\Bc(L^p(X))$ invariant (because $V_k$ are isometries
in $L^p$ too).  Since the set of left $\Fr_\nu$-decay preserving
operators is norm closed in $\Bc(L^2(X))$, we may assume in the rest
of the proof that $S$ is of finite range. According to Lemma
\ref{l:filter}, it suffices to show that, for a given Borel set
$N\in\Nr_\nu$ and for any number $\varepsilon>0$, there is a Borel
set $M\in\Nr_\nu$ such that
$\|\cchi_{M^\crm}(Q)S\cchi_N(Q)\|<\varepsilon$.

In the rest of the proof we shall freely use the notations
introduced in Section \ref{mft} (see also the
proof of Proposition \ref{p:frange}).  In particular, $q$ is defined
by $\frac{1}{p}=\frac{1}{2}+\frac{1}{q}$.  If $f\in L^2(X)$ we have
$$
\|\cchi_Nf\|_{L^p(K_a)}\leq \|\cchi_N\|_{L^q(K_a)}\|f\|_{L^2(K_a)}
\leq|N\cap K_a|^{1/q}\|f\|_{L^2(K_a)}.
$$
Since $N\in\Nr_\nu$ we can find a constant $c$ such that $|N\cap
K_a|\leq c\nu(a)$ (note that the definition (\ref{e:nu}) does not
involve the restriction of $\nu$ to bounded sets). Thus, if we take
$\lambda_a=\nu(a)^{-1/q}$ for $a\in Z\equiv\Z^n$, we get
$\cchi_Nf\in\Lr$ with the notations of Section \ref{mft}. In other
terms, we see that we have $\cchi_N(Q)\in\Bc(L^2(X),\Lr)$.  Let
$T=S\cchi_N(Q)$ and let us assume that we also have $S\in\Bc(\Lr)$.
Then $T\in\Bc(L^2(X),\Lr)$ and we can apply the Maurey type
factorization theorem Theorem \ref{t:maurey}, where $\Hr=L^2(X)$.
Thus we can write $T=g(Q)R$ for some $R\in\Bc(L^2(X))$ and some
function $g\in\Mr$, which means that $G:=\sup_{a\in
  Z}\nu(a)^{-1/q}\|g\|_{L^q(K_a)}$ is a finite number. If $t>0$ and
$M=\{x\mid g(x)>t\}$ then we get for all $a\in Z$:
$$
|M\cap K_a|=\|\cchi_M\|^q_{L^q(K_a)}\leq\|g/t\|^q_{L^q(K_a)}
\leq (G/t)^q \nu(a).
$$

Note that the second condition imposed on $\nu$ in Theorem
\ref{t:fquasi} can be stated as follows: there is an increasing
strictly positive function $\delta$ on $[0,\infty[$ such that
$\nu(b)\leq \delta(|b-a|)\nu(a)$ for all $a, b$.  Indeed, we may
take $\delta(r)=\sup_{|b-a|\leq r}\nu(b)/\nu(a)$.  Now let $a\in X$
and let $D(a)$ be the set of $b\in Z$ such that $K_b$ intersects
$B_a$. Clearly $D(a)$ contains at most $2^n$ points $b$ all of them
satisfying $|b-a|\leq\sqrt n+1$. Hence:
$$
|M\cap K_a|\leq\sum_{b\in D(a)}|M\cap K_b|\leq
2^n\sup_{b\in D(a)}(G/t)^q \nu(b)
\leq 2^n (G/t)^q \delta(\sqrt n+1)\nu(a),
$$
which proves that $M$ belongs to $\Nr_\nu$.
On the other hand, we have:
$$
\|\cchi_{M^\crm}(Q)T\|=\|\cchi_{M^\crm}(Q)g(Q)R\|\leq
\|\cchi_{M^\crm}g\|_{L^\infty}\|R\|\leq t\|R\|.
$$
To finish the proof of the theorem it suffices to take
$t=\varepsilon/\|R\|$.

We still have to prove that $S\in\Bc(\Lr)$. Since $S$ is of finite
range, there is a number $r$ such that $\cchi_a(Q)\cchi_b(Q)=0$ if
$|a-b|\geq r$.  Then for any $f\in\Lr$:
$$
\sum_a\lambda_a^2\|\cchi_aSf\|^2_{L^p}=
\sum_a\lambda_a^2\|\sum_{|b-a|<r}\cchi_aS\cchi_bf\|^2_{L^p}\leq
C\sum_{|b-a|<r}\lambda_a^2\|\cchi_aS\cchi_bf\|^2_{L^p}
$$
where $C$ is a number depending only on $r$ and $n$. Since $S$ is
bounded in $L^p$ the last term is less than
$CC'\sum_{|b-a|<r}\lambda_a^2\|\cchi_bf\|^2_{L^p}$ for some constant
$C'$. Finally, from
$\nu(b)\leq\delta(|b-a|)\nu(a)\leq\delta(r)\nu(a)$ we get
$$
\sum_{|a-b|<r}\lambda_a^2=\sum_{|a-b|<r}\nu(a)^{-2/q}\leq
L(r)\delta(r)^{2/q}\lambda_b^2
$$
where $L(r)$ is the maximum number of points from $Z$ inside a
ball of radius $r$. Thus we have $\|S\|^2_{\Bc(\Lr)}\leq
CC'L(r)\delta(r)^{2/q}$.  \qed

\begin{theorem}\label{th:fquasi} 
  Let $X=\R^n$ and let $S$ be a pseudo-differential operator of
  class $S^0$.  Then $S$ is $\Fc_\w$-decay preserving in $L^2(X)$,
  i.e.\ if $\vphi\in B_\w(X)$ then $\vphi(Q)S=T_1\psi_1(Q)$ and
  $S\vphi(Q)=\psi_2(Q)T_2$ for some $\psi_1,\psi_2\in B_\w(X)$ and
  $T_1,T_2\in\Bc(L^2(X))$.
\end{theorem} 
\proof Since the adjoint of $S$ is also a pseudo-differential
operator of class $S^0$, it suffices to show that $S$ is left
$\Fc_\w$-decay preserving. We have $S\in\Bc(L^p(X))$ for all
$1<p<\infty$ and $S$ is of class $\cu(Q)$ because the commutators
$[Q_j,S]$ are bounded operators for all $1\leq j\leq n$.  Thus we
can apply Theorem \ref{t:fquasi} and deduce that for any function
$\nu$ as in the statement of the theorem, for any $\varepsilon >0$,
and for any $N\in\Nr_\nu$ there is $M\in\Nr_\nu$ such that
$\|\cchi_{M^\crm}(Q)S\cchi_N(Q)\|\leq\varepsilon$.  Now let $N$ be a
Borel $\w$-small set, i.e.\ such that $|N\cap B_a|\rarrow 0$ if
$a\rarrow\infty$. We shall prove that there is a function $\nu$ with
the properties required in Theorem \ref{t:fquasi} and with
$\lim_{a\rarrow\infty}\nu(a)=0$ such that $N\in\Nr_\nu$. This
finishes the proof of the corollary because the relation
$M\in\Nr_\nu$ implies now that $M$ is $\w$-small.

We construct $\nu$ as follows. The relation $\theta(r)=\sup_{|a|\geq
  r}|N\cap B_a|$ defines a positive decreasing function on
$[0,\infty[$ which tends to zero at infinity and such that $|N\cap
B_a|\leq\theta(|a|)$ for all $a\in X$.  We set $\xi(t)=\theta(0)$ if
$0\leq t < 1$ and for $k\geq0$ integer and $2^k\leq t < 2^{k+1}$ we
define $\xi(t)=\max\{\xi(2^{k-1})/2,\theta(2^k)\}$. So $\xi$ is a
strictly positive decreasing function on $[0,\infty[$ which tends to
zero at infinity and such that $\theta\leq\xi$. Moreover, if
$2^k\leq s < 2^{k+1}$ and $2^{k+p}\leq t < 2^{k+p+1}$ then
$$
\xi(t)=\xi(2^{k+p})\geq \xi(2^{k+p-1})/2\geq\ldots\geq
2^{-p}\xi(2^{k})=2^{-p}\xi(s)
$$
hence $\xi(s)\geq\xi(t)\geq \frac{s}{2t}\xi(s)$ if $1\leq s\leq
t$.  We take $\nu(a)=\xi(|a|)$, so $\nu$ is a bounded strictly
positive function on $X$ with $\lim_{a\rarrow\infty}\nu(a)=0$ and
$|N\cap B_a|\leq\nu(a)$ for all $a$. If $a,b$ are points with
$|a|,|b|\geq1$ and $|a-b|\leq r$ then $\nu(b)/\nu(a)\leq1$ if
$|a|\leq|b|$ and if $|a|>|b|$ then
$$
\frac{\nu(b)}{\nu(a)}=
\frac{\xi(|b|)}{\xi(|a|)}\leq\frac{2|a|}{|b|}\leq2(1+r).
$$
Thus the second condition imposed on $\nu$ in Theorem
\ref{t:fquasi} is also satisfied.  \qed

As a final example, we introduce now classes of vanishing at
infinity functions of a more topological nature.  Let us fix a
\emph{uniformly discrete set} $L\subset X$, i.e.\ a set such that
$\inf|a-b|>0$ where the infimum is taken over couples of distinct
points $a,b\in L$.  Let $L_\varepsilon=L+B(\varepsilon)$ be the set
of points at distance $<\varepsilon$ from $L$.  We say that a subset
$N\subset X$ is \emph{$L$-thin} if for each $\varepsilon>0$ there is
$r<\infty$ such that $N\setminus B(r)\subset L_\veps$.  In other
terms, $N$ is $L$-thin if there is a family $\{\delta_a\}_{a\in L}$
of positive real numbers with $\delta_a\rarrow0$ as $a\rarrow\infty$
such that $N\subset\bigcup B_a(\delta_a)$.  The complement of such a
set will be called \emph{$L$-fat}. We denote $\Fc_L$ the family of
$L$-fat sets, we note that $\Fc_L$ is a filter on $X$ contained in
$\Fc_\w$ and finer than the Fr\'echet filter, and we denote $B_L(X)$
the set of bounded Borel functions such that $\lim_{\Fc_L}\vphi=0$.
So $\vphi\in B(X)$ belongs to $B_L(X)$ if and only if the set
$\{|\vphi|\geq\lambda\}$ is $L$-thin for each $\lambda>0$.

\begin{proposition}\label{p:Lql}
  Let $X=\R^n$ and let $S$ be a bounded operator on $L^2(X)$ such
  that on the region $x\neq y$ its distribution kernel is a function
  satisfying the estimate $|S(x,y)|\leq c|x-y|^{-m}$ for some $m>n$.
  Then $S$ is $\Fc_L$-decay preserving.
\end{proposition}
\proof Let $\theta\in\cb(X)$ such that $\theta(x)=0$ on a
neighborhood of the origin and $S_\theta(x,y)=\theta(x-y)S(x,y)$. If
$\xi(x)=\theta(x)|x|^{-m}$ then for the operator $S_\theta$ of
kernel $S_\theta(x,y)$ we have $\|S_\theta u\|\leq c\|\xi*|u|\|$
hence $\|S_\theta\|\leq c\|\xi\|_{L^1}$ By choosing a convenient
sequence of functions $\theta$ we see that $S$ is the norm limit of
a sequence of operators which besides the properties from the
statement of the proposition are such that $S(x,y)=0$ if
$|x-y|>R(S)$. Since the set of $\Fc_L$-decay preserving operators is
closed in norm (see Subsection \ref{ql}), we may assume in the rest
of the proof that the kernel of $S$ has this property. In fact, in
order to simplify the notations and without loss of generality, we
shall assume $S(x,y)=0$ if $|x-y|>1$.

Let $N$ be an $L$-thin Borel set and let $\veps>0$. We shall
construct an $L$-fat Borel set with $F\subset N^\crm$
such that $\|\cchi_N(Q)S\cchi_F(Q)\|\leq\veps$.
Since the adjoint operator $S^*$ has the same properties as $S$, this
suffices to prove that it is decay preserving.

We shall only need two simple estimates. First, if $\rho_x(G)$
is the distance from a Borel set $G$ to a point $x$, then
\begin{equation}\label{e:1}
\int_G\frac{dy}{|x-y|^{2m}}\leq C(m,n)\rho_x(G)^{n-2m}.
\end{equation} 
Then, if $B_0, B$ are two balls with the same center and radiuses
$\delta$ and $\delta+\veps$, then
\begin{equation}\label{e:2}
\int_{B_0}\rho_x(B^\crm)^{n-2m}dx\leq
C(m,n)\varepsilon^{n-2m}\delta^n. 
\end{equation} 
We shall choose $\varepsilon=\delta^{n/2m}$.  Then
$\cchi_{B_0}(Q)S\cchi_{B^\crm}(Q)$ is an operator with integral
kernel and we can estimate its Hilbert-Schmidt norm as follows:
\begin{eqnarray}\label{e:str}
\|\cchi_{B_0}(Q)  S \cchi_{B^\crm}(Q)\|_{HS}^2 &=& \int_{X\times X}
\cchi_{B_0}(x)|S(x,y)|^2\cchi_{B^\crm}(y)dxdy  \nonumber\\
& \leq &
c\int_{B_0}dx\int_{B^\crm}\frac{dy}{|x-y|^{2m}}
\leq C\int_{B_0}\rho_x(B^\crm)^{n-2m}dx \nonumber \\
& \leq & C'\varepsilon^{n-2m}\delta^n=C'\delta^{\lambda}
\end{eqnarray}
where $\lambda=n^2/2m>0$.

We can assume that $N=\bigcup_a B_a(\delta_a)$, where the sequence
of numbers $\delta_a$ satisfies $\delta_a\rarrow0$ as
$a\rarrow\infty$. Denote $N_a=B_a(\delta_a)$ and
$M_a=B_a(\delta_a+\varepsilon_a)$, where we choose
$\varepsilon_a=\delta_a^{n/2m}$ as above.  Choose $r$ such that the
balls $N_a$ are pairwise disjoint and $\delta_a+\varepsilon_a<1$ if
$|a|>r$ and let $R$ such that
$\cchi_{N_a}(Q)S\cchi_{B(R)^\crm}(Q)=0$ if $|a|\leq r$.  Let
$M=\bigcup M_a$ and $F=M^\crm\setminus B(R)$, so that $F$ is a
closed $L$-fat set.  Then for any $u\in L^2(X)$ we have:
\begin{equation*}
\|\cchi_{N}(Q)S\cchi_F(Q)u\|^2=
\sum_{|a|>r}\|\cchi_{N_a}(Q)S\cchi_F(Q)u\|^2.
\end{equation*} 
Since $S$ is of range $1$ we have
$\cchi_{N_a}(Q)S\cchi_{B_a(2)^\crm}(Q)=0$ if $\delta_a<1$.  Thus
$$
\|\cchi_{N}(Q)S\cchi_F(Q)u\|^2\leq
\sum_{|a|>r}\|\cchi_{N_a}(Q)S\cchi_{F\cap B_a(2)}(Q)\|^2\,
\|\cchi_{B_a(2)}(Q)u\|^2
$$
The number of $b\in L$ such that $B_b(2)$ meets $B_a(2)$ is a
bounded function of $a$, hence there is a constant $C$ depending
only on $L$ such that
\begin{equation*}
\|\cchi_{N}(Q)S\cchi_F(Q)u\|\leq 
C\sup_{|a|>r}\|\cchi_{N_a}(Q)S\cchi_{F\cap B_a(2)}(Q)\|\,\|u\|.
\end{equation*} 
We have $F\subset M^\crm\subset M_a^\crm$ hence
$$
\|\cchi_{N_a}(Q)S\cchi_{F\cap B_a(2)}(Q)\|\leq
\|\cchi_{N_a}(Q)S\cchi_{M_a^\crm}(Q)\|_{HS}\leq
C'\delta_a^{\lambda/2} 
$$
because of (\ref{e:str}). So the norm $\|\cchi_{N}(Q)S\cchi_F(Q)\|$
can be made as small as we wish by choosing $r$ large enough.
\qed

\begin{corollary}\label{c:Lql} 
  Let $X=\R^n$, $\mu$ the Lebesgue measure, and $L$ a uniformly
  discrete subset of $\R^n$.  Then a pseudo-differential operator of
  class $S^0$ on $L^2(X)$ is decay preserving with respect to 
  $\Fc_\mu$, $\Fc_\w$ and $\Fc_L$.
\end{corollary}
\proof In the first case we use Theorem \ref{th:meas} by taking into
account that a pseudo-differential operator of class $S^0$ belongs
to $\Bc(L^p(X))$ for all $1<p<\infty$ and that the adjoint of such
an operator is also pseudo-differential of class $S^0$.  The second
case has already been considered in Theorem \ref{th:fquasi}. For the
third case, note that the distribution kernel of such an operator
verifies the estimates $|S(x,y)|\leq C_k|x-y|^{-n}(1+|x-y|)^{-k}$
for any $k>0$, see \cite{Ho}.  \qed

\section{Weakly vanishing perturbations}\label{wvp}

In this subsection we reconsider the framework of Subsection
\ref{div} and improve, but with a stronger assumption
$a\in\Bc(\Kr)$, the decay condition (\ref{e:sd8}).  We shall
consider on $\Hr$ the class of ``vanishing at infinity'' functions
corresponding to the algebra $B_\w(X)$, in other terms we equip
$\Hr$ with the Hilbert module structure associated to the multiplier
algebra $\{\varphi(Q)|\varphi\in B_\w(X)\}$. By Lemma \ref{l:best},
$(\Gr,\Hr)$ remains a compact Friedrichs module. The space $\Kr$
inherits a natural direct sum Hilbert module structure.

We keep the notations and terminology of Sections \ref{div} and
\ref{wvpr}. We recall that an operator
$D^*aD:\Hr^m\rarrow\Hr^{-m}$ is coercive if there are numbers
$\mu,\nu>0$ such that
\begin{equation}\label{e:coercivity}
\Re\la Du,a Du\rangle\geq \mu\|u\|^2_{\Hr^m}-\nu\|u\|^2_{\Hr}
\hspace{3mm}\forall  u\in\Hr^m.
\end{equation}   
Clearly the next lemma remains true if the filter $\Fc_\w$ is
replaced by $\Fc_\mu$ or $\Fc_L$.

\begin{lemma}\label{l:analytic} 
  Assume that $a\in\Bc(\Kr)$ is $\Fc_\w$-decay preserving and that
  the operator $D^*aD:\Hr^m\rarrow\Hr^{-m}$ is coercive.  Then
  $D(\Delta_a-z)^{-1}$ is $\Fc_\w$-decay preserving if $\Re
  z\leq-\nu$, where $\nu$ is as in (\ref{e:coercivity}).
\end{lemma} 
\proof We shall use Proposition \ref{pr:pco} with $c$ the the
identity operator in $\Kr$, so $\Delta\equiv\Delta_c$ is the
operator in $\Hr$ associated to $D^*D=\sum_{|\alpha|\leq
m}P^{2\alpha}$, which is the canonical (Riesz) positive isomorphism
of $\Gr$ onto $\Gr^*$ and (\ref{e:coercivity}) means $\Re
D^*aD\geq\mu D^*D -\nu$. We have $D(\Delta-z)^{-1}\in\bq(\Hr,\Kr)$
and $D(D^*D-z)^{-1}D^*\in\bq(\Kr)$ if $\Re z<0$ because these
operators consist of matrices of pseudo-differential operators with
constant coefficients of class $S^0$, so we can use Theorem
\ref{th:fquasi}.  \qed

We now consider two operators $\Hr^m\rarrow\Hr^{-m}$ of the form
$$
D^*aD=\sum_{|\alpha|, |\beta|\leq m} 
P^\alpha a_{\alpha\beta} P^\beta \hspace{3mm} \mbox{and}
\hspace{3mm} 
D^*bD=\sum_{|\alpha|, |\beta|\leq m} 
P^\alpha b_{\alpha\beta} P^\beta
$$
where the coefficients are continuous operators
$
a_{\alpha\beta},b_{\alpha\beta}:
\Hr^{m-|\beta|}\rarrow\Hr^{|\alpha|-m}
$
satisfying some other conditions stated below and denote as usual
$\Delta_a$ and $\Delta_b$ the operators in $\Hr$ associated to
them. 

\begin{theorem}\label{t:w} 
Assume that the operators $D^*aD$ and $D^*bD$ are coercive and that
their coefficients satisfy the following conditions: (1)
$a_{\alpha\beta}\in\Bc(\Hr)$ and are $\Fc_\w$-decay preserving
operators; (2) if $|\alpha|+|\beta|=2m$ then
$a_{\alpha\beta}-b_{\alpha\beta}$ is left $\Fc_\w$-vanishing at
infinity; (3) if $|\alpha|+|\beta|<2m$ then
$a_{\alpha\beta}-b_{\alpha\beta}\in
\Kc(\Hr^{m-|\beta|},\Hr^{|\alpha|-m})$.  Then the operator
$\Delta_b$ is a compact perturbation of $\Delta_a$, in particular
$\se(\Delta_a)=\se(\Delta_b)$.
\end{theorem}
\proof
We check the conditions of Theorem \ref{t:main-2}. Because of
the coercivity assumptions, condition (1) is fulfilled, and (3) is
satisfied by Lemma \ref{l:analytic}. The part of condition (2)
involving the coefficients such that $|\alpha|+|\beta|=2m$ is
satisfied by definition, for the lower order coefficients it
suffices to use (\ref{e:coo}).  \qed

\begin{remark}\label{r:..}{\rm
    If $a_{\alpha\beta}$ and $b_{\alpha\beta}$ are bounded Borel
    functions and $a_{\alpha\beta}-b_{\alpha\beta}\in B_\w(X)$ for
    all $\alpha,\beta$, then the conditions (1)-(3) of the theorem
    are satisfied. Indeed, in order to check the compactness
    conditions on the lower order coefficients note that, by Lemma
    \ref{l:best}, if $\varphi\in B_\w(X)$ then the operator
    $\varphi(Q):\Hr^s\rarrow\Hr^{-t}$ is compact if $s,t\geq0$ and
    one of them is not zero.  }\end{remark}

The next result is a more general but less explicit version of
Theorem \ref{t:w}. This is an improvement of \cite[Theorem 2.1]{OS},
thus it covers some subelliptic operators.

\begin{theorem}\label{t:wos} 
  Assume that $D^*aD$ satisfies (\ref{e:coercivity}) and that
  $\Delta_b$ is a closed densely defined operator such that there is
  $z\in\rho(\Delta_b)$ with $\Re z\leq-\nu$. Moreover, assume that
  $a,b$ satisfy the conditions (1)-(3) of Theorem \ref{t:w}. Then
  the operator $\Delta_b$ is a compact perturbation of $\Delta_a$.
\end{theorem}
\proof We shall apply Theorem \ref{t:main} with $A=\Delta_a$ and
$B=\Delta_b$. The assumption (AB) is clearly satisfied and we take
$\wtilde A=D^*aD$ and $\wtilde B=D^*bD$, hence $\wtilde B-\wtilde
A=D^*(b-a)D$. Then let $S=D^*$ and $T=(b-a)D$.
\qed

Finally, let us note that one should be able to use Theorem
\ref{t:main} to treat situations when the coefficients
$a_{\alpha\beta}$ and $b_{\alpha\beta}$ are unbounded operators even
if $|\alpha|=\beta|=m$ (as in \cite[Theorem 3.1]{OS} and
\cite{Ba1,Ba2}), see the framework of Example \ref{e:ab} and
Corollary \ref{c:main-aa}, but we shall not pursue this idea here.

\section{Riemannian manifolds}\label{lrm}

Let $\Hr, \Kr$ be two Hilbert spaces identified with their adjoints
and $\d$ a closed densely defined operator mapping $\Hr$ into $\Kr$.
Let $\Gr=\Dc(\d)$ equipped with the graph norm, so $\Gr\subset\Hr$
continuously and densely and $\d\in \Bc(\Gr, \Kr)$.

Then the quadratic form $\|\d u\|^2_\Kr$ on $\Hr$ with domain $\Gr$
is positive densely defined and closed. Let $\Delta$ be the positive
self-adjoint operator on $\Hr$ associated to it. In fact
$\Delta=\d^*\d$, where the adjoint $\d^*$ of $\d$ is a closed
densely defined operator mapping $\Kr$ into $\Hr$.

Now let $\lambda\in\Bc(\Hr)$ and $\Lambda\in\Bc(\Kr)$ be
self-adjoint and such that $\lambda\geq c$ and $\Lambda\geq c$ for
some real $c>0$.  Then we can define new Hilbert spaces
$\widetilde{\Hr}$ and $\widetilde{\Kr}$ as follows:
$$
\leqno{\mbox{\bf($*$)}}\hspace{1mm}
\left\{
\begin{array}{ll}
\widetilde{\Hr}=\Hr \mbox{ as vector space and } \langle u\mid v
\rangle_{\widetilde\Hr}=\langle u\mid \lambda v\rangle_\Hr,
\\
\widetilde{\Kr}=\Kr \mbox{ as vector space and } \langle u\mid v
\rangle_{\widetilde\Kr}=\langle u\mid \Lambda v\rangle_\Kr.
\end{array}
\right.
$$

Since $\Hr=\widetilde\Hr$ and $\Kr=\widetilde\Kr$ as topological
vector spaces, the operator $\d:\Gr\subset\widetilde\Hr\rightarrow
\widetilde\Kr$ is still a closed densely defined operator, hence the
quadratic form $\|\d u\|_{\widetilde{\Kr}}^2$ on $\widetilde{\Hr}$
with domain $\Gr$ is positive, densely defined and closed. We shall
denote by $\widetilde\Delta$ the positive self-adjoint operator on
$\widetilde{\Hr}$ associated to it.

We can express $\widetilde \Delta$ in more explicit terms as
follows.  Denote by $\widetilde \d$ the operator $\d$ when viewed as
acting from $\widetilde\Hr$ to $\widetilde\Kr$. Then
$\widetilde\Delta=\widetilde{\d}^*\widetilde \d$, where
$\widetilde{\d}^*:\Dc(\widetilde{\d}^*)\subset
\widetilde{\Kr}\rightarrow\widetilde{\Hr}$ 
is the adjoint of $\widetilde \d=\d$ with respect
to the new Hilbert space structures (the spaces $\widetilde{\Hr},
\widetilde{\Kr}$ being also identified with their adjoints). It is
easy to check that $\widetilde{\d}^*=\lambda^{-1}\d^*\Lambda$. Thus
$\widetilde{\Delta}=\lambda^{-1}\d^*\Lambda \d$.

\medskip

Now let $(X,\rho)$ be a proper locally compact metric space
(see the definition before Corollary \ref{ex27}) and let us
assume that $\Hr$ and $\Kr$ are Hilbert $X$-modules.

\begin{definition}
A closed densely defined map $\d:\Dc(\d)\subset\Hr\rarrow\Kr$ is a
\emph{first order operator} if there is $C\in\R$ such that for
each bounded Lipschitz function $\varphi$ on $X$
the form $[\d, \varphi(Q)]$ is a bounded operator and
$\|[\d, \varphi(Q)]\|_{\Bc(\Hr, \Kr)}\leq
C\,\mbox{\rm  Lip } \varphi$.
\end{definition}
Here
\begin{equation*}
\mbox{ Lip } \varphi=
\inf_{x\neq y}|\varphi(x)-\varphi(y)|\rho(x,y)^{-1}.
\end{equation*}
In more explicit terms, we require
\begin{equation*}
|\langle \d^* u, \varphi(Q)v\rangle_\Hr- \langle u,
\varphi(Q)\d v\rangle_{\Kr}|
\leq C\, \mbox{\rm  Lip }\varphi\, \|u\|_\Kr\|v\|_\Hr
\end{equation*}
for all $u\in\Dc(\d^*)$ and $v\in\Dc(d)$. Thus $\langle \d^* u,
\varphi(Q) v\rangle- \langle u, \varphi(Q) \d v\rangle$ is a
sesquilinear form on the dense subspace $\Dc(\d^*)\times\Dc(\d)$ of
$\Kr\times\Hr$ which is continuous for the topology induced by
$\Hr\times \Kr$. Hence there is a unique continuous operator $[\d,
\varphi(Q)] :\Hr\rightarrow\Kr$ such that
\begin{equation*}
\langle \d^* u, \varphi(Q)v\rangle_\Hr
- \langle u, \varphi(Q)\d v\rangle_\Kr =\langle u,[\d, \varphi(Q)] v
\rangle_\Kr
\end{equation*}
for all $u\in \Dc(\d^*)$, $v\in\Dc(\d)$ and $\|[\d,
\varphi(Q)]\|_{\Bc(\Hr, \Kr)}\leq C\,\mbox{\rm Lip }\varphi$.
\begin{lemma}\label{l:40}
The operator $\d(\Delta+1)^{-1}$ is decay preserving.
\end{lemma}
\proof We shall prove that $S:=\d(\Delta+1)^{-1}$ is a decay preserving
operator with the help of Corollary \ref{ex27}, more precisely we
show that $[S, \varphi(Q)]$ is a bounded operator if $\varphi$ is a
positive Lipschitz function. Let $\varepsilon >0$ and
$\varphi_\epsilon=\varphi(1+\varepsilon\varphi)^{-1}$.  Then
$\varphi_\varepsilon$ is a bounded function with
$|\varphi_\varepsilon|\leq\varepsilon^{-1}$ and
\begin{equation*}
|\varphi_\varepsilon(x)-\varphi_\varepsilon(y)|
=\frac{|\varphi(x)-\varphi(y)|}{(1+\varepsilon\varphi(x))
(1+\varepsilon\varphi(y))}\leq|\varphi(x)-\varphi(y)|
\end{equation*}
hence Lip $\varphi_\varepsilon\leq$ Lip $\varphi$. Let $v\in\Dc(\d)$
we have for all $u\in \Dc(\d^*)$:
\begin{eqnarray*}
|\langle \d^* u, \varphi_\epsilon(Q)v\rangle_\Hr|
&=&|\langle u, \varphi_\epsilon(Q)\d v\rangle_\Kr
+\langle u,[\d, \varphi_\varepsilon(Q)]  v\rangle_\Kr|\\
&\leq& \|u\|_\Kr(\varepsilon^{-1} \|\d v\|_\Kr+
C\,\mbox{\rm  Lip }\varphi_{\varepsilon}\,\|u\|_\Hr).
\end{eqnarray*}
Hence $\varphi_\varepsilon(Q)v\in\Dc(\d^{**})=\Dc(\d)$ because $\d$
is closed.  Thus $\varphi_\varepsilon(Q)\Dc(\d)\subset\Dc(\d)$ and
by the closed graph theorem we get
$\varphi_\varepsilon(Q)\in\Bc(\Gr)$, where $\Gr$ is the domain of
$\d$ equipped with the graph topology.  This also implies that
$\varphi_\varepsilon(Q)$ extends to an operator in $\Bc(\Gr^*)$
(note that $\varphi_\varepsilon(Q)$ is symmetric in $\Hr$).

Now, if we think of $\d$ as a continuous operator
$\Gr\rightarrow\Kr$, then it has an adjoint
$\d^*:\Kr\rightarrow\Gr^*$ which is the unique continuous extension
of the operator $\d^*:\Dc(\d^*)\subset\Kr\rightarrow\Hr\subset
\Gr^*$. Thus the canonical extension of $\Delta$ to an element of
$\Bc(\Gr, \Gr^*)$ is the product of $\d:\Gr\rightarrow \Kr$ with
$\d^*:\Kr\rightarrow\Gr^*$ (note $\Dc(\d)$ is the form domain of
$\Delta$). Then it is trivial to justify that we have in $\Bc(\Gr,
\Gr^*)$:
\begin{eqnarray*}
[\Delta, \varphi_\varepsilon(Q)]=[\d^*, \varphi_\varepsilon(Q)]\d
+\d^*[\d, \varphi_\varepsilon(Q)].
\end{eqnarray*}
Here $[\d^*, \varphi_\varepsilon(Q)]=[\varphi_\varepsilon(Q),\d]^*\in
\Bc(\Kr, \Hr)$. Since $\Delta+1:\Gr\rightarrow\Gr^*$ is a linear
homeomorphism, we then have in $\Bc(\Gr^*, \Gr)$:
\begin{eqnarray*}
[\varphi_\varepsilon(Q), (\Delta+1)^{-1}]&=&
(\Delta+1)^{-1}[\Delta, \varphi_\varepsilon(Q)](\Delta+1)^{-1}\\
&=&(\Delta+1)^{-1}[\varphi_\varepsilon(Q),\d]^* \d(\Delta+1)^{-1}\\
&+&(\Delta+1)^{-1}\d^*[\d, \varphi_\varepsilon(Q)](\Delta+1)^{-1}.
\end{eqnarray*}
Finally, taking once again into account the fact that
$\varphi_\varepsilon(Q)$ leaves $\Gr$ invariant, we have:
\begin{eqnarray*}
[\varphi_\varepsilon(Q), \d(\Delta+1)^{-1}]&=&
[\varphi_\varepsilon(Q),\d](\Delta+1)^{-1}\\
&&+\, \d(\Delta+1)^{-1}
[\varphi_\varepsilon(Q),\d]^* \d(\Delta+1)^{-1}
\\
&&+ \d(\Delta+1)^{-1}\d^*[\d, \varphi_\varepsilon(Q)](\Delta+1)^{-1}.
\end{eqnarray*}
Hence:
\begin{eqnarray*}
\|[\varphi_\varepsilon(Q), \d(\Delta+1)^{-1}]\|_{\Bc(\Hr, \Kr)}&\leq&
\|[\varphi_\varepsilon(Q),\d]\|_{\Bc(\Hr, \Kr)}
\|(\Delta+1)^{-1}\|_{\Bc(\Hr)}\\
&&\hspace{-5cm}+\, \|\d(\Delta+1)^{-1}\|_{\Bc(\Hr, \Kr)}
\|[\varphi_\varepsilon(Q),\d]^*\|_{\Bc(\Kr, \Hr)}
\|\d(\Delta+1)^{-1}\|_{\Bc(\Hr, \Kr)}\\
&&\hspace{-5cm}+\, \|\d(\Delta+1)^{-1}\d^*\|_{\Bc(\Kr, \Kr)}
\|[\d, \varphi_\varepsilon(Q)]\|_{\Bc(\Hr, \Kr)}
\|(\Delta+1)^{-1}\|_{\Bc(\Hr)}.
\end{eqnarray*}
The most singular factor here is
\begin{eqnarray*}
\|\d(\Delta+1)^{-1}\d^*\|_{\Bc(\Kr, \Kr)}
\leq \|\d\|_{\Bc(\Gr, \Kr)}\|(\Delta+1)^{-1}\|_{\Bc(\Gr^*, \Gr)}
\|\d^*\|_{\Bc(\Kr, \Gr^*)}
\end{eqnarray*}
and this is finite. Thus we get for a finite constant $C_1$:
\begin{eqnarray*}
\|[\varphi_\varepsilon(Q), \d(\Delta+1)^{-1}]\|_{\Bc(\Hr, \Kr)}
&\leq& C_1\|[\d, \varphi_\varepsilon(Q)]\|_{\Bc(\Hr, \Kr)}\\
&\leq&C_1\, C\, \mbox{\rm  Lip }
 \varphi_\varepsilon\leq C_1\, C\, \mbox{\rm  Lip }\varphi
\end{eqnarray*}
Now let $u\in\Kr_{\rm c}$ and $v\in\Hr_{\rm c}$. We get:
\begin{eqnarray*}
|\langle\varphi(Q)u, \d(\Delta+1)^{-1}v \rangle-
\langle u, \d(\Delta+1)^{-1}\varphi(Q)v\rangle|&=&\\
&&\hspace{-7cm}=\,\lim_{\varepsilon\rightarrow 0}
|\langle\varphi_\varepsilon(Q)u, \d(\Delta+1)^{-1}v \rangle-
\langle u, \d(\Delta+1)^{-1}\varphi_\varepsilon(Q)v\rangle|\\
&&\hspace{-7cm} \leq\, C_1\, C\, \mbox{\rm  Lip }\varphi
\end{eqnarray*}
Thus $[\varphi(Q), \d(\Delta+1)^{-1}]$ is a bounded operator.\qed
\begin{theorem}\label{metric}
Let $(X, \rho)$ be a proper locally compact metric space. Assume that
$(\Gr, \Hr)$ is a compact Friedrichs $X$-module and that $\Kr$
is a Hilbert $X$-module. Let $\d, \lambda, \Lambda$ be operators
satisfying the following conditions:\\
(i) $\d$ is a closed first order operator from $\Hr$ to $\Kr$ with
$\Dc(\d)=\Gr$; \\
(ii) $\lambda$ is a bounded self-adjoint operator on $\Hr$ with $\inf
\lambda>0$ and such that $\lambda -1\in \Kc(\Gr, \Hr)$ (e.g.
$\lambda-1\in\Bc_0(\Hr)$);\\
(iii) $\Lambda$ is a bounded self-adjoint operator on $\Kr$ with
$\inf\Lambda>0$ and such that $\Lambda-1\in \Bc_0(\Kr)$.\\
Then the self-adjoint operators $\Delta$ and $\widetilde{\Delta}$
have the same essential spectrum.
\end{theorem}
\proof In this proof, we shall consider $\widetilde{\Delta}$ as an
operator acting on $\Hr$. Since $\widetilde{\Hr}=\Hr$ as topological
vector spaces and the notion of spectrum is purely topological,
$\widetilde{\Delta}$ is a closed densely defined operator on $\Hr$
and it has the same spectrum as the self-adjoint
$\widetilde{\Delta}$ on $\widetilde{\Hr}$.  Moreover, if we define
the essential spectrum $\sigma_{\rm ess} (A)$ as the set of $z\in\C$
such that either $\ker(A-z)$ is infinite dimensional or the range of
$A-z$ is not closed, we see that the essential spectrum is a
topological notion, so $\sigma_{\rm ess} (\widetilde{\Delta})$ is
the same, whether we think of $\widetilde{\Delta}$ as operator on
$\Hr$ or on $\widetilde\Hr$. Finally, with this definition of
$\sigma_{\rm ess}$ we have $\sigma_{\rm ess}(A)=\sigma_{\rm ess}(B)$
if $(A-z)^{-1}-(B-z)^{-1}$ is a compact operator for some
$z\in\rho(A)\cap \rho(B)$.

Thus it suffices to prove that
$(\Delta+1)^{-1}-(\widetilde{\Delta}+1)^{-1} \in\Kc(\Hr)$. Now we
observe that
\begin{equation*}
\widetilde{\Delta}+1=\lambda^{-1}d^*\Lambda d+1=\lambda^{-1}
(\d^*\Lambda \d+ \lambda)
\end{equation*}
and $\Delta_\Lambda=\d^*\Lambda \d$ is the positive self-adjoint
operator on $\Hr$ associated to the closed quadratic form $\|\d
u\|^2_{\widetilde{\Kr}}$ on $\Hr$ with domain $\Gr$. Thus
$(\widetilde{\Delta}+1)^{-1}= (\Delta_\Lambda+\lambda)^{-1}\lambda$
and
\begin{equation*}
(\widetilde{\Delta}+1)^{-1}-(\Delta_\Lambda+\lambda)^{-1} =
(\Delta_\Lambda+\lambda)^{-1}(\lambda-1)=
[(\lambda-1)(\Delta_\Lambda+\lambda)^{-1}]^*
\end{equation*}
The range of $(\Delta_\Lambda+\lambda)^{-1}$ is included in the form
domain of $\Delta_\Lambda+\lambda$, which is $\Gr$. The map
$(\Delta_\Lambda+\lambda)^{-1}:\Hr\rightarrow\Gr$ is continuous, by
the closed graph theorem, and $\lambda-1:\Gr\rightarrow\Hr$ is
compact.  Hence
$(\widetilde{\Delta}+1)^{-1}-(\Delta_\Lambda+\lambda)^{-1}$ is
compact.  Similarly:
\begin{equation*}
(\Delta+1)^{-1}-(\Delta_\Lambda+\lambda)^{-1} =
(\d^*\d+1)^{-1}-(\d^*\Lambda \d+1)^{-1}\in\Kc(\Hr)
\end{equation*}
For this we use Theorem \ref{t:main-2} with: $\Er=\Kr$, $D=\d$,
$a=1$, $b=\Lambda$ and $z=-1$. Since $\d^*\d$ and $\d^*\Lambda d$
are positive self-adjoint operators on $\Hr$ with the same form
domain $\Gr$, the first condition of Theorem \ref{t:main-2} is
satisfied. Then the second condition holds because
$\Lambda-1\in\bol(\Kr)$. Thus it remains to observe that the
operator $\d(\Delta+1)^{-1}$ is decay preserving by Lemma \ref{l:40}.
\qed

We shall consider now an application of Theorem \ref{metric} to
concrete Riemannian manifolds. It will be clear from what follows
that we could treat Lipschitz manifolds with measurable metrics (see
\cite{DP,Hi,Te,We} for example), but the case of $C^1$ manifolds
with locally bounded metrics suffices as an example. Note, however,
that the arguments of the proof of Theorem \ref{t:metric} cover
without any modification the case when $X$ is not $C^1$ but is a
Lipschitz manifold and a countable atlas has been specified, because
then the tangent space are well defined almost everywhere and the
absolute continuity notions that we use make sense.

\medskip

From now on in this section
\emph{ $X$ is a non-compact differentiable manifold of class $C^1$}.
Then its cotangent manifold $T^*X$ is a topological vector fiber
bundle over $X$ whose fiber over $x$ will be denoted $T^*_xX$.  If
$u:X\rarrow\R$ is differentiable then $\d u(x)\in T^*_xX$ is its
differential at the point $x$ and its differential $\d u$ is a
section of $T^*X$. Thus for the moment $\d$ is a linear map defined
on the space of real $C^1(X)$ functions to the space of sections of
$T^*X$.

A measurable locally bounded Riemannian structure on $X$ will be
called an \emph{R-structure} on $X$. To be precise, an R-structure
is given on $X$ if each $T^*_xX$ is equipped with a quadratic (i.e.\
generated by scalar product) norm $\|\cdot\|_x$ such that:
$$
\leqno{\mbox{\bf($R$)}}\hspace{1mm}
\left\{
\begin{array}{ll}
\mbox{if } v \mbox{ is a continuous section of  }T^*X
\mbox{ over a compact set }K \mbox{ such that }\\v(x)\neq0
\mbox{ for }x\in K,\mbox{ then }x\mapsto\|v(x)\|_x
\mbox{ is a bounded Borel map on }\\
K\mbox{ and }\|v(x)\|_x\geq c
\mbox{ for some number }c>0 \mbox{ and all }x\in K.
\end{array}
\right.
$$
Such a structure allows one to construct a metric compatible with
the topology on $X$, the distance between two points being the
infimum of the length of the Lipschitz curves connecting the points
(see the references above). Since $X$ was assumed to be non-compact,
the metric space $X$ is proper in the sense defined in Subsection
\ref{ql} if and only if it is complete. If this is the case, we say
that \emph{the R-structure is complete}.

It will also be convenient to complexify these structures (i.e.\ 
replace $T^*_xX$ by $T^*_xX\otimes\C$ and extend the scalar product
as usual) and to keep the same notations for the complexified
objects.

We shall consider positive measures  $\mu$ on $X$ such that:
$$
\leqno{\mbox{\bf($M$)}}\hspace{1mm}
\left\{
\begin{array}{ll}
\mu \mbox{ is absolutely continuous and its
density is locally bounded}\\
\mbox{and locally bounded from below by
strictly positive constants}.
\end{array}
\right.
$$
A couple consisting of an R-structure and a measure verifying (M)
on $X$ will be called an \emph{RM-structure on $X$}. The definition
of a complete RM-structure is obvious. To an R-structure we may
canonically associate an RM-structure by taking $\mu$ equal to the
Riemannian volume element.

If an RM-structure is given on $X$ then we may consider the two
Hilbert spaces $\Hr=L^2(X,\mu)$ and $\Kr$ defined as the completion
of the space of continuous sections with compact support of $T^*X$
equipped with the norm
$$\|v\|_\Kr^2=\int_X\|v(x)\|_x^2 d\mu(x).$$
In fact, $\Kr$ is the space of (suitably defined) square integrable
sections of $T^*X$.

The operator of exterior differentiation $\d$ induces a linear map
$\cc^1(X)\rarrow\Kr$ which is easily seen to be closable as operator
from $\Hr$ to $\Kr$ (this is a purely local problem and the
hypotheses we put on the metric and the measure allow us to reduce
ourselves to the Euclidean case). We shall keep the notation $\d$
for its closure and we note that its domain $\Gr$ is the first order
Sobolev space $\Hr^1$ defined in this context as the closure of
$\cc^1(X)$ under the norm
$$
\|u\|^2_{\Hr^1}=\int_X\Big(|u(x)|^2+\|\d u(x)\|^2_x\Big)d\mu(x).
$$
The self-adjoint operator $\Delta=\d^*\d$ in $\Hr$ associated to
the quadratic form $\|\cdot\|^2_{\Hr^1}$ is the \emph{Laplace
  operator} associated to the given RM-structure. This is a
generalized form of the Laplace operator associated to the
Riemannian structure of $X$ because $\mu$ is not necessarily the
Riemannian volume element.

Two RM-structures $(\{\|\cdot\|_x\}_{x\in X},\mu)$ and
$(\{\|\cdot\|'_x\}_{x\in X},\mu')$ on $X$ are called
\emph{equivalent} if there are bounded Borel functions
$\alpha,\beta,\lambda$ on $X$ with $\alpha\geq c$ and $\lambda\geq
c$ for some number $c>0$ such that
$\alpha(x)\|\cdot\|_x\leq\|\cdot\|_x'\leq\beta(x)\|\cdot\|_x$ for
all $x$ and $\mu'=\lambda\mu$. The distances $\rho,\rho'$ on $X$
associated to these structures satisfy $a\rho\leq\rho'\leq b\rho$ for
some numbers $b\geq a>0$, hence if one of the RM-structures is
complete, the second one is also complete.
Notice that the spaces $\Hr,\Kr$ associated to equivalent
RM-structures are identical as topological vector spaces.

Two RM-structures are \emph{strongly equivalent} if they are
equivalent and if the functions $\alpha,\beta,\lambda$ can be chosen
such that $\lambda(x)\rarrow 1$, $\alpha(x)\rarrow 1$ and
$\beta(x)\rarrow 1$ as $x\rarrow\infty$.

\begin{theorem}\label{t:metric}
  The Laplace operators associated to strongly equivalent complete
  RM-structures on $X$ have the same essential spectrum.
\end{theorem}
\proof We check that the assumptions of Theorem \ref{metric} are
satisfied.  We noted above that $X$ is a proper metric space for the
metric associated to the initial Riemann structure.  The spaces
$\Hr,\Kr$ have obvious $X$-module structures and for each
$\vphi\in\cc(X)$ the operator $\vphi(Q):\Hr^1\rarrow\Hr$ is compact.
Indeed, by using partitions of unity, we may assume that the support
of $\vphi$ is contained in the domain of a local chart and then we
are reduced to a known fact in the Euclidean case.  Thus $(\Gr,\Hr)$
is a compact Friedrichs $X$-module.  To see that $\d$ is a first
order operator we observe that if $\vphi$ is Lipschitz then
$[\d,\vphi]$ is the operator of multiplication by the differential
$\d\vphi$ of $\vphi$ and the estimate
$\mbox{ess-sup\,}\|\d\vphi(x)\|_x\leq \mbox{Lip }\vphi$ is easy to
obtain. The conditions on $\lambda$ in Theorem \ref{metric} are
trivially verified. So it remains to consider the operator
$\Lambda$. For each $x\in X$ there is a unique operator
$\Lambda_0(x)$ on $T^*_xX$ such that $\la u|v\ra_x '=\la
u|\Lambda_0(x)v\ra_x $ for all $u,v\in T^*_xX$ and we have
$\alpha(x)^2\leq\Lambda_0(x)\leq\beta(x)^2$ by hypothesis. Here the
inequalities must be interpreted with respect to the initial scalar
product on $T^*_xX$.  Thus the operator $\Lambda$ on $\Kr$ is just
the operator of multiplication by the function
$\Lambda(x)=\lambda(x)\Lambda_0(x)$ and the condition (iii) of
Theorem \ref{metric} is clearly satisfied.  \qed

The (strong) equivalence of two R-structures is defined in an
obvious way. Note that if $\mu,\mu'$ are the Riemannian measures
associated to two strongly equivalent R-structures then the unique
function $\lambda$ such that $\mu'=\lambda\mu$ satisfies
$\lambda(x)\rarrow 1$ as $x\rarrow\infty$.

\begin{corollary}\label{co:metric}
  The Laplace operators associated to strongly equivalent
  complete R-structures on $X$ have the same essential spectrum.
\end{corollary}

We stress that if one of the Riemannian structures is locally
Lipschitz then this result is easy to prove by using local
regularity estimates for elliptic equations.

An assumption of the form $\alpha(x)\rarrow1$ as $x\rarrow\infty$
imposed in the definition of strong equivalence means that the set
where $|\alpha(x)-1|>\varepsilon$ is relatively compact for any
$\varepsilon>0$. We shall consider now a weaker notion of
equivalence associated to the filter $\Fc_\mu$ introduced in Example
\ref{ex:meas}.

We first introduce two notions which clearly depend only on the
equivalence class of an RM-structure.  We say that an RM-structure
is of \emph{infinite volume} if $\mu(X)=\infty$. We say that it has
the \emph{F-embedding property} if for each Borel set $F\subset X$
of finite measure the operator $\cchi_F(Q):\Hr^1\rarrow\Hr$ is
compact.

\begin{remark}\label{re:fprop}{\rm
    The F-embedding property is satisfied under quite general
    conditions. Indeed, the compactness of
    $\cchi_F(Q):\Hr^1\rarrow\Hr$ is equivalent to the compactness of
    the operator $\cchi_F(Q)(\Delta+1)^{-1/2}$ in $\Hr$. Or the set
    of functions $\varphi\in C([0,\infty[)$ such that
    $\cchi_F(Q)\vphi(\Delta)$ is compact is a closed
    $C^*$-subalgebra of $C([0,\infty[)$ so it suffices to find one
    function $\vphi$ which generates this algebra such that
    $\cchi_F(Q)\vphi(\Delta)$ be compact. But
    $\cchi_F(Q)\vphi(\Delta)$ is compact if and only if
    $\cchi_F(Q)|\vphi(\Delta)|^2\cchi_F(Q)$ is compact, so we see
    that it suffices to show that for each Borel set $F$ of finite
    measure there is $t>0$ such that
    $\cchi_F(Q){\rm e}^{-t\Delta}\cchi_F(Q)$ be compact. For
    example, it suffices that this operator be Hilbert-Schmidt,
    i.e.\ that the integral kernel $P_t$ of ${\rm e}^{-t\Delta}$ be
    such that
    $\int_{F\times F} |P_t(x,y)|^2d\mu(x)d\mu(y)<\infty$, which is
    true if $P_t$ satisfies a Gaussian upper estimate and the
    measure of a ball of radius $t^{1/2}$ is bounded below by a
    strictly positive constant (see \cite{AC,ACDH} and references
    there). 
}\end{remark}

Two infinite volume RM-structures will be called
\emph{$\mu$-strongly equivalent} if they are equivalent and if the
functions $\alpha,\beta,\lambda$ can be chosen such that for each
$\varepsilon >0$ the set where one of the inequalities
$|\alpha(x)-1|>\varepsilon$, $|\beta(x)-1|>\varepsilon$ or
$|\alpha(x)-1|>\varepsilon$ holds is of finite measure.

We say that an RM-structure is \emph{regular} if there is $p>2$
such that $\d(\Delta+1)^{-1}$ induces a bounded operator in $L^p$.
More precisely, this means that there is a constant $C$ such that if
$u\in\ L^2(X)\cap L^p(X)$ then $\d(\Delta+1)^{-1} u$, which is a
section of $T^*X$ of finite $L^2$ norm, has an $L^p$ norm bounded by
$C\|u\|_{L^p}$. If the operator $\d(\Delta+1)^{-1}d^*$ also induces
a bounded operators in $L^p$ (in an obvious sense), we say that
the RM-structure is \emph{strongly regular}. From the relation
$\d(\Delta+1)^{-1}d^*=[\d(\Delta+1)^{-1/2}][\d(\Delta+1)^{-1/2}]^*$
we see that strong regularity follows from: there is $\varepsilon>0$
such that $\d(\Delta+1)^{-1/2}$ induces a bounded operator in $L^p$
for $2-\varepsilon<p<2+\varepsilon$.

\begin{theorem}\label{th:wmetric}
  Let $\Delta$ be the Laplace operator associated to an infinite
  volume complete RM-structure on $X$ which has the F-embedding
  property and is regular. Then the Laplace operator associated to
  an RM-structure $\mu$-strongly equivalent to the given structure
  has the same essential spectrum as $\Delta$.
\end{theorem}
\proof Let $\Lambda(x)$ be as in the proof of Theorem
\ref{t:metric}.  Clearly there is a number $C>0$ such that
$C^{-1}\leq \Lambda(x)\leq C$ for all $x$ and such that for each
$\varepsilon >0$ the set where $\|\Lambda(x)-1\|>\varepsilon$ is of
finite measure (the inequalities and the norm are computed on
$T^*_xX$, which is equipped with the initial scalar product).

Now we proceed as in the proof of Theorem \ref{metric} but this time
we equip $\Hr$ and $\Kr$ with the Hilbert module structures
described in Example \ref{ex:meas}. To avoid confusions, we denote
$\Bc_\mu(\Hr)$ and $\Bc_\mu(\Kr)$ the space of decay improving
operators relatively to these new module structures. The F-embedding
property implies that $(\Hr^1,\Hr)$ is a compact Friedrichs module.
Moreover, the operator $\lambda(Q)-1:\Hr^1\rarrow\Hr$ is compact.
Then, as in the proof of Theorem \ref{metric}, we see that it
suffices to prove that 
$$
(\d^*\d+1)^{-1}-(\d^*\Lambda(Q) \d+1)^{-1}\in\Kc(\Hr).
$$
Clearly $\Lambda(Q)-1\in\Bc_\mu(\Kr)$. Now we use Theorem
\ref{t:main-2} exactly as in the proof of Theorem \ref{metric} and
we see that the only condition which remains to be checked is (3) of
Theorem \ref{t:main-2}, i.e.\ in our case
$\d(\Delta+1)^{-1}\in\bqr(\Hr,\Kr)$, where the decay preserving
property is relatively to the algebra $B_\mu(X)$. But this follows
from Theorem \ref{th:meas}. 
\qed

One may check the regularity property needed in Theorem
\ref{th:wmetric} by using the results from \cite{AC,ACDH} concerning
the boundedness in $L^p$ of the operator $d\Delta^{-1/2}$.  For
example, it suffices that $X$ be complete, with the doubling volume
property, and such that the Poincar\'e inequality holds in $L^2$
sense. Note, however, that these results are much stronger than
necessary in our context and it seems reasonable to think that the
boundedness of $d(\Delta+1)^{-1/2}$ holds under less restrictive
assumptions.

The next result does not require regularity assumptions on any of
the RM-structures that we want to compare but only on a third one in
their equivalence class. Observe that each equivalence class of
RM-structures contains one of the same degree of local smoothness as
the manifold $X$ (make local regularizations and use a partition of
unity).

\begin{theorem}\label{th:wmetric*}
  Let $\Delta_a,\Delta_b$ be the Laplace operators associated to
  $\mu$-strongly equivalent complete RM-structures of infinite
  volume and having the F-embedding property. If these structures
  are equivalent to a strongly regular RM-structure, then
  $\Delta_a$ and $\Delta_b$ have the same essential spectrum.
\end{theorem}
\proof Let $\Delta_c$ be the Laplace operator associated to the
third structure. From Theorem \ref{th:meas} it follows that
$\d(\Delta_c+1)^{-1}$ and $\d(\Delta_c+1)^{-1}D^*$ are right
$\Fc_\mu$-decay preserving. Then from Proposition \ref{pr:pco} we
see that $\d(\Delta_a+1)^{-1}$ is right $\Fc_\mu$-decay preserving
and we may conclude as in the proof of Theorem \ref{th:wmetric}.
\qed

\begin{remark}\label{re:conj}{\rm
    It is natural to consider an analog of the filter $\Fc_\w$
    introduced in Section \ref{wvpr} to get an optimal weak decay
    condition for the stability of the essential spectrum in the
    present context. The techniques of Section \ref{wvpr} should be
    relevant for this question.  }\end{remark}

\noindent{\bf Remarks on Laplace operators acting on forms:}
We shall describe here, without going into details, an abstract
framework for the study of the Laplace operator acting on forms. Let
$\Hr$ be a Hilbert space and $\d$ a closed densely defined operator
in $\Hr$ such that $\d^2=0$. For example, $\Hr$ could be the space
of square integrable differential forms over a Lipschitz manifold
and $\d$ the operator of exterior differentiation. We denote
$\delta=\d^*$ and we assume that $\Gr:=\Dc(\d)\cap\Dc(\delta)$ is
dense in $\Hr$ (which is a rather strong condition in the context of
this paper, e.g.\ in the preceding example it is a differentiability
condition on the metric). Then let $D=\d+\delta$ with domain $\Gr$,
observe that $\|Du\|^2=\|\d u\|^2+\|\delta u\|^2$ so $D$ is a closed
symmetric operator, assume that $D$ is self-adjoint, and define
$\Delta=D^2=\d\delta+\delta\d$ (form sum). Then
\begin{equation}\label{eq:forms}
(\Delta+1)^{-1}=(D+i)^{-1}(D-i)^{-1}.
\end{equation}
Now let $a\in\Bc(\Hr)$ with $a\geq\varepsilon >0$ and such that
$a^{\pm1}\Gr\subset\Gr$ and let $\Hr_a$ be the Hilbert space which
is equal to $\Hr$ as vector space but is equipped with the new the
scalar product $\la u,v\ra_a=\la u,av\ra$. Denote $\d_a$ the
operator $\d$ viewed as operator acting in $\Hr_a$ with adjoint
$\delta_a=a^{-1}\delta a$.  We can define as above operators $D_a$
(with domain $\Gr_a=\Gr$) and $\Delta_a=D_a^2$ which are
self-adjoint in $\Hr_a$ and satisfy a relation similar to
(\ref{eq:forms}). Then $\Delta_a$ is a compact perturbation of
$\Delta$ if the operators $(D_a\pm i)^{-1}-(D\pm i)^{-1}$ are
compact and this last condition is equivalent to the compactness of
the operator $D_a-D:\Gr\rarrow\Gr^*$. And this holds if $(\Gr,\Hr)$
is a compact Hilbert $X$-module over a metric space $X$ and
$a-1\in\Bc_0(\Hr)$.

\section{On Maurey's factorization theorem}\label{mft}
The subject of this section is quite different from that of the rest
of the paper: we shall prove a version of a factorization theorem
due to Bernard Maurey which plays an important role in several
arguments from the main part of this article. We first recall
Maurey's result, cf.\  Theorems 2 and 8 in \cite{Ma}. 
\begin{theorem}\label{th:bmsft}
  Let $1<p<2$ and let $T$ be an arbitrary continuous linear map from
  a Hilbert space $\Hr$ into $L^p$. Then there is $R\in\Bc(\Hr,L^2)$
  and there is a function $g\in L^q$, where
  $\frac{1}{p}=\frac{1}{2}+\frac{1}{q}$, such that $T=g(Q)R$.
\end{theorem}
We have stated only the particular case we need of the theorem (the
result extends easily to larger classes of Banach spaces $\Hr$).
The $L^p$ spaces refer to an arbitrary positive measure space
$(X,\mu)$.

Before going on to our main purpose, we shall state an easy
consequence of this theorem which is needed in Sections \ref{wvp}
and \ref{lrm}. Let $\{\Kr(x)\}_{x\in X}$ be a measurable family of
Hilbert spaces (see \cite[Ch.\ II]{Dix}) such that the dimension of
$\Kr(x)$ is $\leq N$ for some finite $N$.  Let
$\Kr=\int_X^\oplus\Kr(x) d\mu(x)$ be the corresponding direct
integral and for each $p\geq1$ let $\Kr_p$ be the space of
($\mu$-equivalence classes) of measurable vector fields $v$ such
that $\int_X\|v(x)\|^p_{\Kr(x)}d\mu(x)<\infty$. Thus $\Kr_p$ is
naturally a Banach space and $\Kr_2=\Kr$.

\begin{corollary}\label{co:bmsft}
  Let $\Hr$ be a Hilbert space and let $T\in\Bc(\Hr,\Kr_p)$ with
  $1<p<2$. Then there is $R\in\Bc(\Hr,\Kr)$ and there is a function
  $g\in L^q$, where $q=2p/(2-p)$, such that $T=g(Q)R$.
\end{corollary}
\proof For each $n=1,\dots,N$ let $X_n$ be the set of $x$ such that
the dimension of $\Kr(x)$ is equal to $n$. Then $X$ is the disjoint
union of the measurable sets $X_n$. For each $x$ there is $n$ such
that $x\in X_n$ and we can choose a unitary map
$j(x):\Kr(x)\rarrow\C^n$ such that $\{j_x\}$ be a measurable family
of operators. Let $J$ be the operator acting on vector fields
according to the rule $(Jv)(x)=j(x)v(x)$, let $\Pi_n$ be the
operator of multiplication by $\cchi_{X_n}$, and let $T_n\equiv
\Pi_nJT\in\Bc(\Hr,L^p(X_n;\C^n))$. We can write $T_n=(T_n^k)_{1\leq
  k\leq n}$ with $T_n^k\in\Bc(\Hr,L^p(X_n))$ and Maurey's theorem
gives us a factorization $T_n^k=g_n^k(Q)S_n^k$ with
$S_n^k\in\Bc(\Hr,L^2(X_n))$ and $g_n^k\in L^q(X_n)$, and clearly we
may assume $g_n^k\geq0$. Let $g_n=\sup_k g_n^k\in L^q(X_n)$ and
$S_n\in\Bc(\Hr,L^2(X_n;\C^n))$ be the operator with components
$(g_n^kg_n^{-1})(Q)R^k_n$. Then $T_n=g_n(Q)S_n$ and if we define
$R_n=J^{-1}S_n$ we get
$$
g_n(Q)R_n=J^{-1}g_n(Q)S_n=J^{-1}T_n=\Pi_nT.
$$
Thus, if we define $g=\sum_n\cchi_{X_n}g_n$ and $R=\sum\Pi_nR_n$, we
get $T=g(Q)R$.
\qed

\medskip

Our purpose in the rest of this section is to extend Theorem
\ref{th:bmsft} (in the case $X=\R^n$) to more general classes of
spaces of measurable functions, which do not seem to be covered by
the results existing in the literature, cf.\ \cite{Kr}.  Our proof
follows closely that of Maurey.  We first recall Ky Fan's Lemma,
see \cite[9.10]{DJT}.
\begin{proposition}\label{p:kfl} 
  Let $\Kc$ be a compact convex subset of a Hausdorff topological
  vector space and let $\Fr$ be a convex set of functions
  $F:\Kc\rarrow]-\infty,+\infty]$ such that each $F\in\Fr$ is convex
  and lower semicontinuous. If for each $F\in\Fr$ there is $g\in\Kc$
  such that $F(g)\leq0$, then there is $g\in\Kc$ such that
  $F(g)\leq0$ for all $F\in\Fr$.
\end{proposition} 

We need a second general fact that we state below. Let $(X,\mu)$ be
a $\sigma$-finite positive measure space and let $L^0(X)$ be the
space of $\mu$-equivalence classes of complex valued measurable
functions on $X$ with the topology of convergence in measure.  Let
$\Lr$ be a Banach space with $\Lr\subset L^0(X)$ linearly and
continuously and such that if $f\in L^0(X)$, $g\in\Lr$ and $|f|\leq
|g|$ ($\mu$-a.e.)\ then $f\in\Lr$ and $\|f\|_\Lr\leq\|g\|_\Lr$.  The
next result is a rather straightforward consequence of Khinchin's
inequality \cite[1.10]{DJT} (see also \cite[Section 8]{Pi}).
\begin{proposition}\label{p:type2}
  There is a number $C$, independent of $\Lr$, such that for any
  Hilbert space $\Hr$ and any $T\in\Bc(\Hr,\Lr)$ the following
  inequality holds
\begin{equation}\label{e:type2}
\|(\mbox{$\sum_j$}|Tu_j|^2)^{1/2}\|_\Lr\leq
C\|T\|_{\Bc(\Hr,\Lr)}(\mbox{$\sum_j$}\|u_j\|^2)^{1/2}
\end{equation}
for all finite families $\{u_j\}$ of vectors in $\Hr$.
\end{proposition} 

From now on we work in a setting adapted to our needs in Section
\ref{wvp}, although it is clear that we could treat by the same
methods a general abstract situation.  Let $X=\R^n$ equipped with
the Lebesgue measure, denote $Z=\Z^n$, and for each $a\in Z$ let
$K_a=a+K$, where $K=]-1/2,1/2]^n$, so that $K_a$ is a unit cube
centered at $a$ and we have $X=\bigcup_{a\in Z}K_a$ disjoint union.
Let $\cchi_a$ be the characteristic function of $K_a$ and if
$f:X\rarrow\C$ let $f_a=f|K_a$. We fix a number $1<p<2$ and a family
$\{\lambda_a\}_{a\in Z}$ of strictly positive numbers $\lambda_a>0$
and we define $\Lr\equiv\ell^2_\lambda(L^p)$ as the Banach space of
all (equivalence classes) of complex functions $f$ on $X$ such that
\begin{equation}\label{e:lr} 
\|f\|_\Lr:=
\Big(\sum_{a\in Z}\|\lambda_a\cchi_af\|^2_{L^p}\Big)^{1/2}<\infty.
\end{equation} 
Here $L^p=L^p(X)$ but note that, by identifying $\cchi_af\equiv
f_a$, we can also interpret $\Lr$ as a conveniently normed direct
sum of the spaces $L^p(K_a)$, see \cite[page XIV]{DJT}.  If
$\lambda_a=1$ for all $a$ we set $\ell^2_\lambda(L^p)=\ell^2(L^p)$.
Observe that $\ell^2(L^2)=L^2(X)$.

Let $q$ be given by $\frac{1}{p}=\frac{1}{2}+\frac{1}{q}$, so that
$1<p<2<q<\infty$.  We also need the space
$\Mr\equiv\ell^\infty_\lambda(L^q)$ defined by the condition
\begin{equation}\label{e:mr} 
\|g\|_\Mr:=\sup_{a\in Z}\|\lambda_a\cchi_ag\|_{L^q}<\infty.
\end{equation} 
The definitions are chosen such that $
\|gu\|_\Lr\leq\|g\|_\Mr\|u\|_{L^2} $ where $L^2=L^2(X)$. As
explained in \cite[page XV]{DJT}, the space $\Mr$ is naturally
identified with the dual space of the Banach space $\Mr_*\equiv
\ell^1_{\lambda^{-1}}(L^{q'})$, where $\frac{1}{q}+\frac{1}{q'}=1$,
defined by the norm
$$
\|h\|_{\Mr_*}:=\sum_{a\in Z}\|\lambda^{-1}_a\cchi_ah\|_{L^{q'}}.
$$
Below, when we speak about $w^*$-topology on $\Mr$ we mean the
$\sigma(\Mr,\Mr_*)$-topology. Clearly
$$
\Mr_1^+=\{g\in\Mr\mid g\geq0,\|g\|_\Mr\leq1\}
$$
is a convex compact subset of $\Mr$ for the $w^*$-topology.
\begin{lemma}\label{l:max} 
For each $f\in\Lr$ there is $g\in\Mr_1^+$ such that
$\|f\|_\Lr=\|g^{-1}f\|_{L^2}$.
\end{lemma} 
\proof
We can assume $f\geq0$. Since $1=\frac{p}{2}+\frac{p}{q}$, we have:
$$
\|f_a\|_{L^p}=\|f_a\|^{p/2}_{L^p}\|f_a\|^{p/q}_{L^p}=
\|f_a^{p/2}\|_{L^2}\|f_a^{p/q}\|_{L^q}=
\|f_a^{-p/q}f\|_{L^2}\|f_a^{p/q}\|_{L^q}
$$
with the usual convention $0/0=0$. Now we define $g_a$ on $K_a$
as follows.  If $f_a=0$ then we take any $g_a\geq0$ satisfying
$\lambda_a\|g_a\|_{L^q}=1$.  If $f_a\neq0$ let
$$
g_a=\lambda_a^{-1}\big(f_a/\|f_a\|_{L^p}\big)^{p/q}=
\lambda_a^{-1}\|f_a^{p/q}\|_{L^q}^{-1}f_a^{p/q}.
$$
Thus we have  $\lambda_a\|g_a\|_{L^q}=1$ for all $a$, in particular
$\|g\|_\Mr=1$. By the preceding computations we also have 
$\|f_a\|_{L^p}=\|g_a^{-1}f_a\|_{L^2}\|g_a\|_{L^q}$ and so
$$
\|f\|_\Lr^2=\sum\lambda_a^2\|f_a\|^2_{L^p}=
\sum\lambda_a^2\|g_a\|^2_{L^q} \|g_a^{-1}f_a\|^2_{L^2}=
\sum\|g_a^{-1}f_a\|^2_{L^2}
$$
which is just $\|g^{-1}f\|^2_{L^2}$.
\qed

The main technical result follows.
\begin{proposition}\label{p:maurey1} 
  Let $(f^u)_{u\in U}$ be a family of functions in $\Lr$ such that,
  for each $\alpha=(\alpha_u)_{u\in U}$ with $\alpha_u\in\R$,
  $\alpha_u\geq0$ and $\alpha_u\neq0$ for at most a finite number of
  $u$, the function $f^\alpha:=(\sum_u|\alpha_uf^u|^2)^{1/2}$
  satisfies $\|f^\alpha\|_\Lr\leq\|\alpha\|_{\ell^2(U)}$. Then there
  is $g\in\Mr^+_1$ such that $\|g^{-1}f^u\|_{L^2}\leq1$ for all
  $u\in U$.
\end{proposition} 
\proof
For each $\alpha$ as in the statement of the proposition we define
a function $F_\alpha:\Mr_1^+\rarrow ]-\infty,+\infty]$ as follows:
$$
F_\alpha(g) =\|g^{-1}f^\alpha\|^2_{L^2}-\|\alpha\|_{\ell^2(U)}^2=
\sum_u\alpha_u^2\big(\|g^{-1}f^u\|^2_{L^2}-1\big).
$$
Our purpose is to apply Proposition \ref{p:kfl} with
$\Kr=\Mr_1^+$ equipped with the $w^*$-topology and $\Fr$ equal to
the set of all functions $F_\alpha$ defined above. We saw before
that $\Kr$ is a convex compact set.  From the second representation
of $F_\alpha$ given above it follows that $\Fr$ is a convex set.
Each $F_\alpha$ is a convex function because
$\|g^{-1}f^\alpha\|^2_{L^2}=\int g^{-2}(f^{\alpha})^2\d x$ and the
map $t\mapsto t^{-2}$ is convex on $[0,\infty[$. We shall prove in a
moment that $F_\alpha$ is lower semicontinuous. From Lemma
\ref{l:max} it follows that there is $g_\alpha\in\Kr$ such that
$\|f^\alpha\|_\Lr=\|g_\alpha^{-1}f^\alpha\|_{L^2}$. Then by our
assumptions we have
$$
F_\alpha(g_\alpha)=\|f^\alpha\|_\Lr^2-\|\alpha\|_{\ell^2(U)}^2\leq0.
$$
From Ky Fan's Lemma it follows that one can choose $g\in\Kr$ such
that $F_\alpha(g)\leq0$ for all $\alpha$, which finishes the proof
of the proposition.

It remains to show the lower semicontinuity of $F_\alpha$. For this
it suffices to prove that $g\mapsto\|g^{-1}f\|^2_{L^2}\in[0,\infty]$
is lower semicontinuous on $\Kr$ if $f\in\Lr, f\geq0$. But
$$
\|g^{-1}f\|^2_{L^2}=\sum_a\int_{K_a}g_a^{-2}f_a^{2}\d x
$$
and the set of lower semicontinuous functions
$\Kr\rarrow[0,\infty]$ is stable under sums and upper bounds of
arbitrary families. Hence it suffices to prove that each map
$g\mapsto\int_{K_a}g_a^{-2}f_a^{2}\d x$ is lower semicontinuous.
This map can be written as a composition $\phi\circ J_a$ where
$J_a:\Mr\rarrow L^q(K_a)$ is the restriction map $J_ag=g_a$ and
$\phi:L^q(K_a)\rarrow[0,\infty]$ is defined by
$\phi(\theta)=\int_{K_a}\theta^{-2}f_a^{2}\d x$.  The map $J_a$ is
continuous if we equip $L^q(K_a)$ with the weak topology and $\Mr$
with the $w^*$-topology because it is the adjoint of the norm
continuous map $L^{q'}(K_a)\rarrow\Mr_*$ which sends $u$ into the
function equal to $u$ on $K_a$ and $0$ elsewhere. Thus it suffices
to show that $\phi$ is lower semicontinuous on the positive part of
$L^q(K_a)$ equipped with the weak topology and for this we can use
exactly the same argument as Maurey. We must prove that the set
$\{\theta\in L^q(K_a)\mid \theta\geq0,\phi(\theta)\leq r \}$ is
weakly closed for each real $r$. Since $\phi$ is convex, this set is
convex, so it suffices to show that it is norm closed.  But this is
clear by the Fatou Lemma.  \qed

\begin{theorem}\label{t:maurey} 
  Let $\Hr$ be a Hilbert space and $T:\Hr\rarrow\Lr$ a linear
  continuous map.  Then there exist a linear continuous map
  $R:\Hr\rarrow L^2(X)$ and a positive function $g\in\Mr$ such that
  $T=g(Q)R$.
\end{theorem} 
\proof Let $U$ be the unit ball of $\Hr$ and for each $u\in U$ let
$f^u=Tu$. From Proposition \ref{p:type2} we get
$$
\|f^\alpha\|_\Lr=\|(\mbox{$\sum_u$}|T(\alpha_uu)|^2)^{1/2}\|_\Lr\leq
A(\mbox{$\sum_u$}\|\alpha_uu\|^2)^{1/2}\leq
A(\mbox{$\sum_u$}|\alpha_u|^2)^{1/2}
$$
where $A=C\|T\|_{\Bc(\Hr,\Lr)}$. Since there is no loss of
generality in assuming $A\leq1$, we see that the assumptions of
Proposition \ref{p:maurey1} are satisfied. So there is $g\in\Mr_1^+$
such that $\|g^{-1}Tu\|_{L^2(X)}\leq 1$ for all $u\in U$. Thus it
suffices to define $R$ by the rule $Ru=g^{-1}Tu$ for all $u\in\Hr$.
\qed

\appendix
\section{Appendix}\label{pc}

Let $(\Gr,\Hr)$ be a Friedrichs couple and
$\Gr\subset\Hr\subset\Gr^*$ the Gelfand triplet associated to it.
To an operator $S\in \Bc(\Gr, \Gr^*)$ (which is the same as a
continuous sesquilinear form on $\Gr$) we associate an operator
$\widehat{S}$ acting in $\Hr$ according to the rules:
$\Dc(\widehat{S})= S^{-1}(\Hr)$, $\widehat{S}=S|{\Dc(\widehat{S})}$.
Due to the identification $\Gr^{**}=\Gr$, the operator $S^*$ is an
element of $\Bc(\Gr, \Gr^*)$, so $\widehat{S^*}$ makes sense. On the
other hand, if $\widehat{S}$ is densely defined in $\Hr$ then the
adjoint $\widehat{S}^*$ of $\widehat{S}$ with respect to $\Hr$ is
also well defined and we clearly have $\widehat{S^*}\subset
\widehat{S}^*$.

\begin{lemma}\label{prop33}
  If $S-z:\Gr\rightarrow \Gr^*$ is bijective for some $z\in\C$, then
  $\widehat{S}$ is a closed densely defined operator, we have
  $\widehat{S}^*=\widehat{S^*}$ and $z\in\rho(\what S)$.  Moreover,
  the domains $\Dc(\widehat{S})$ and $\Dc(\widehat{S}^*)$ are dense
  subspaces of $\Gr$.
\end{lemma}
\proof Clearly we can assume $z=0$. From the bijectivity of
$S:\Gr\rarrow\Gr^*$ and the inverse mapping theorem it follows that
$S$ and $S^*$ are homeomorphisms of $\Gr$ onto $\Gr^*$. Since $\Hr$
is dense in $\Gr^*$, we see that $\Dc(\what S)$ and
$\Dc(\what{S^*})$ are dense in $\Gr$, hence in $\Hr$. Since
$\what{S^*}\subset\what S^*$, the operator $\what S^*$ is also
densely defined in $\Hr$. Thus $\what S$ is densely defined and
closable. We now show that it is closed.  Consider a sequence of
elements $u_n\in\Dc(\what S)$ such that $u_n\rarrow u$ and $\what
Su_n\rarrow v$ in $\Hr$.  Then $Su_n\rarrow v$ in $\Gr^*$ hence,
$S^{-1}$ being continuous, $u_n\rarrow S^{-1}v$ in $\Gr$, so in
$\Hr$. Hence $u=S^{-1}v\in\Dc(\what S)$ and $\what Su=v$.

We have proved that $\what S$ is densely defined and closed and
clearly $0\in\rho(\what S)$. Then we also have $0\in\rho(\what
S^*)$, so $\what S^*:\Dc(\what S^*)\rarrow\Hr$ is bijective. Since
$\what{S^*}:\Dc(\what{ S^*})\rarrow\Hr$ is also bijective and $\what
S^*$ is an extension of $\what{S^*}$, we get $\what{S^*}=\what S^*$.
\qed

A standard example of operator satisfying the condition required in
Lemma \ref{prop33} is a {\em coercive} operator, i.e. such that
$\Re\langle u, Su\rangle\geq\mu\|u\|_{\Gr}^2-\nu\|u\|^2_\Hr$ for
some strictly positive constants $\mu,\nu$ and all $u\in\Gr$.
Indeed, replacing $S$ by $S+\nu$, we may assume $\Re\la
u,Su\ra\geq\mu\|u\|_\Gr^2$.  Since $S^*$ verifies the same estimate,
this clearly gives $\|Su\|_{\Gr^*}\geq\mu\|u\|_\Gr$ and
$\|S^*u\|_{\Gr^*}\geq\mu\|u\|_\Gr$ for all $u\in\Gr$. Thus $S$ and
$S^*$ are injective operators with closed range, which implies that
they are bijective.

\medskip

If $A$ is a self-adjoint operator on $\Hr$ then there is a natural
Gelfand triplet associated to it, namely
$\Dc(|A|^{1/2})\subset\Hr\subset\Dc(|A|^{1/2})^*$.  Then $A$ extends
to a continuous operator $A_0:\Dc(|A|^{1/2})\rightarrow
\Dc(|A|^{1/2})^*$ which fulfills the conditions of Lemma
\ref{prop33} and one has $\what A_0=A$. In our applications it is
interesting to know whether there are other Gelfand triplets
$\Gr\subset\Hr\subset\Gr^*$ with $\Dc(A)\subset\Gr$ and such that
$A$ extends to a continuous operator $\Gr\rightarrow\Gr^*$.  For not
semibounded operators, e.g.\ for Dirac operators, many other
possibilities exist such that $\Gr$ is not comparable to
$\Dc(|A|^{1/2})$.  But if $A$ is semibounded, then the class of
spaces $\Gr$ is rather restricted, as the next lemma shows.

\begin{lemma}\label{lem34}
Assume that $A$ is a bounded from below self-adjoint operator on $\Hr$
and such that $\Dc(A)\subset\Gr$ densely. Then $A$ extends to a
continuous operator $\widetilde{A}:\Gr\rightarrow\Gr^*$ if and
only if \ $\Gr\subset \Dc(|A|^{1/2})$ and in this case
$\widetilde{A}=A_0|_{\Gr}$.
\end{lemma}
\proof We prove only the nontrivial implication of the lemma.  So
let us assume that $A$ extends to some $\widetilde{A}\in \Bc(\Gr,
\Gr^*)$.  Replacing $A$ by $A+\lambda$ with $\lambda$ a large enough
number, we can assume that $A\geq 1$. For $u\in \Dc(A)$ we have
\begin{equation*}
\|A^{1/2}u\|_\Hr=\sqrt{\langle u, Au\rangle}
=\sqrt{\langle u, \widetilde{A}u\rangle}\leq C\|u\|_\Gr,
\end{equation*}
where $C^2=\|\widetilde{A}\|_{\Gr\rightarrow\Gr^*}$. Since $\Dc(A)$
is dense in $\Gr$, it follows that the inclusion map
$\Dc(A)\rightarrow \Dc(A^{1/2})$ extends to a continuous linear map
$J:\Gr\rightarrow \Dc(A^{1/2})$. If $u\in \Gr$ then there is a
sequence $\{u_n\}$ in $\Dc(A)$ such that $u_n\rightarrow u$ in
$\Gr$. Then $J(u_n)\rightarrow J(u)$ in $\Dc(A^{1/2})$. Since $\Gr$
and $\Dc(A^{1/2})$ are continuously embedded in $\Hr$ we shall have
$u_n\rightarrow u$ in $\Hr$ and $u_n=J(u_n)\rightarrow J(u)$ in
$\Hr$, hence $J(u)=u$ for all $u\in\Gr$. In other terms, $\Gr\subset
\Dc(A^{1/2})$.  \qed

We note that, under the conditions of the lemma, the inclusions
$\Dc(A)\subset\Gr$ and $\Gr\subset \Dc(|A|^{1/2})$ are continuous
(by the closed graph theorem), so we have a scale
\begin{equation*}
\Dc(A)\subset\Gr \subset \Dc(|A|^{1/2})\subset \Hr
\subset \Dc(|A|^{1/2})^*\subset\Gr^*\subset \Dc(A)^*
\end{equation*}
with continuous and dense embeddings (because $\Dc(A)$ is dense in
$\Dc(|A|^{1/2})$).


\end{document}